\overfullrule=0pt
\def\multi#1{\vbox{\baselineskip=0pt\halign{\hfil$\scriptstyle\vphantom{(_)}##$\hfil\cr#1\crcr}}}
%%%%%%%%%%%%%%%%%%%%%%%%%%%%%%%%%%%%%%%%%%%%%%%%%%%%%
%
%  General Macros and setup
%
%%%%%%%%%%%%%%%%%%%%%%%%%%%%%%%%%%%%%%%%%%%%%%%%%%%%%
\def\sqbox{\copy\squ\hskip -.7pt}
\def\ptt
{{
\hbox{$\vbox{\offinterlineskip
\hbox{\sqbox\sqbox}
\vskip -.4pt
\hbox{\sqbox\sqbox}}$}
}}

\def\potr
{{
\hbox{$\vbox{\offinterlineskip
\hbox{\sqbox}
\vskip -.4pt
\hbox{\sqbox\sqbox\sqbox}}$}
}}
\def\poot
{{
\hbox{$\vbox{\offinterlineskip
\hbox{\sqbox}
\vskip -.6pt
\hbox{\sqbox}
\vskip -.7pt
\hbox{\sqbox\sqbox}}$}
}}
\def\picture #1 by #2 (#3){
  \vbox {\hbox {\lower #2 \hbox
    {\special{picture #3}}} % this is the low-level interface
    \hrule width #1 height 0pt depth 0pt
    \vfill
    }
  }

\def\scaledpicture #1 by #2 (#3 scaled #4){{
  \dimen0=#1 \dimen1=#2
  \divide\dimen0 by 1000 \multiply\dimen0 by #4
  \divide\dimen1 by 1000 \multiply\dimen1 by #4
  \picture \dimen0 by \dimen1 (#3 scaled #4)}
  }

%    Un peu de place aprs

\font\Sc=cmcsc10
\font\Ch=cmbx10
%\font\Ch=msbm9
\font\gros=cmbx10 scaled\magstep2
\font\super=cmbx12 scaled\magstep2
\font\small=cmr8

\font\ita=cmsl10 
\font\bol=cmbx10

\font\Sc=cmcsc10
\font\super=cmbx12 scaled\magstep2

\hrule height 0pt
\abovedisplayskip=3mm
\belowdisplayskip=3mm
\abovedisplayshortskip=0mm
\belowdisplayshortskip=2mm
\normalbaselineskip=12pt
\parskip=2mm plus 0.5mm minus 0.5mm
\normalbaselines
\def\heute{\ifcase\month\or January\or February\or March\or April\or May
\or June\or July\or August\or September\or October\or November\or December\fi
        \space\number\day, \number\year}
\def\abstract#1{\bigskip\bigskip\centerline{\bf Abstract} \medskip
   {\parshape=1 .5in 5in \small\baselineskip=10pt #1\smallskip}}
\def\nl{\bigskip\noindent}
\newcount\sectionno
\sectionno = 0
\def\section#1
{\global\advance\sectionno by 1
 \goodbreak\par\bigskip\bigskip\centerline{\Sc #1}
\par\bigskip
}

\newcount\refctr  \refctr = 0
\def\enumerate#1 #2\endenumerate{\advance\refctr by 1
  \expandafter\edef\csname #1\endcsname{\the\refctr}%
  \if !#2! \else\enumerate #2\endenumerate\fi}
\newcount\numerobib \numerobib=0
\def\bibitem{\advance\numerobib by 1
    \par\noindent\item{[\number\numerobib]}}
\def\N{\hbox{\Ch N}}

\newcount\elemno
\elemno=0
\def\numero{\global\advance\elemno by 1
            {\the\elemno}.}
\def\TITRECO#1{\goodbreak\nl {\bf #1}\ \numero\ \ \ \ }
\def\Rem{\TITRECO{Remark}}

\outer\def\Thm#1#2\par{\bigbreak
         \noindent{{\bf Theorem}\ #1\quad}{\sl#2}\par
          \ifdim\lastskip<\medskipamount \removelastskip\penalty55\medskip\fi}
\outer\def\Prop#1#2\par{\bigbreak
         \noindent{{\bf Proposition}\ #1\quad}{\sl#2}\par
          \ifdim\lastskip<\medskipamount \removelastskip\penalty55\medskip\fi}
\outer\def\Cor#1#2\par{\bigbreak
         \noindent{{\bf Corollairy}\ #1\quad}{\sl#2}\par
          \ifdim\lastskip<\medskipamount \removelastskip\penalty55\medskip\fi}
\outer\def\Lem#1#2\par{\bigbreak
         \noindent{{\bf Lemma}\ #1\quad}{\sl#2}\par
          \ifdim\lastskip<\medskipamount \removelastskip\penalty55\medskip\fi}
\def\Proof{\nl {\it Proof\/}.\ \ \ }
\outer\def\Conj#1#2\par{\bigbreak
         \noindent{{\bf Conjecture\ #1\quad}}{\sl#2}\par
          \ifdim\lastskip<\medskipamount \removelastskip\penalty55\medskip\fi}

\newbox\squ  % box character for ends of proofs
\def\sqbox{\copy\squ\hskip -.2pt}

\def\ligne#1{\ifcase#1
\or\hbox{\sqbox}
\or\hbox{\sqbox\sqbox}
\or\hbox{\sqbox\sqbox\sqbox}
\or\hbox{\sqbox\sqbox\sqbox\sqbox}
\or\hbox{\sqbox\sqbox\sqbox\sqbox\sqbox}
\or\hbox{\sqbox\sqbox\sqbox\sqbox\sqbox\sqbox}
\or\hbox{\sqbox\sqbox\sqbox\sqbox\sqbox\sqbox\sqbox}
\or\hbox{\sqbox\sqbox\sqbox\sqbox\sqbox\sqbox\sqbox\sqbox}\fi}

\def\un#1{\ligne#1}

\def\deux#1#2{{\hbox{$\vbox{\offinterlineskip
\ligne#2
\vskip -.2pt
\ligne#1}$}}}

\def\trois#1#2#3{{\hbox{$\vbox{\offinterlineskip
\ligne#3
\vskip -.2pt
\ligne#2
\vskip -.2pt
\ligne#1}$}}}

\def\quatre#1#2#3#4{{\hbox{$\vbox{\offinterlineskip
\ligne#4
\vskip -.2pt
\ligne#3
\vskip -.2pt
\ligne#2
\vskip -.2pt
\ligne#1}$}}}

\nopagenumbers
\pageno=1
%\headline={CRM \heute, {\it Montr\'eal-San Diego} \hfill Preliminary Version.\hfill \folio}
%\headline={\vsmall Bergeron-Garsia\hfill Science Fiction and Macdonald Polynomials\hfill \folio}

\footline{\hfill\folio\hfill}

%%%%%%%%%%%%%%%%%%%%%%%%%%%%%%%%%%%%%%%%%%%%%%%%%%%%%
%
%  Macros of notation for the paper
%
%%%%%%%%%%%%%%%%%%%%%%%%%%%%%%%%%%%%%%%%%%%%%%%%%%%%%
\def\noc {\supseteq}
\def \DA {\downarrow}
\def \da {\downarrow}
\def  \lar   {\leftarrow}
\def  \rar   {\rightarrow}
\def\CF{{\cal F}}

\def \flip {{\bf flip}}

\def \TH {{\tilde H}}
\def \xon {x_1,\ldots ,x_n}
\def \yon {y_1, \ldots ,y_n}
\def \LLL {\langle}
\def \RRR {\rangle}

\def \om {\omega}

\def\TH{\widetilde{H}}

\def \CH{{\cal H}}
\def \CB {{\cal B}}
\def \CL {{\cal L}}
\def \CJ {{\cal J}}

\def\BM{{\rm {\bf M}}}
\def\TK{\widetilde{K}}

\def\part{\vdash}
\def\charac{{\raise 2pt\hbox{$\chi$}}}

\def\scalar#1{\langle #1\rangle}
\def\frac#1#2{{#1\over #2}}

\def\del{\partial}
\def\flp{{\bf flip}}
\def\char{{\rm char}\,}
\def\dim{{\rm dim}\,}
\def\eee{\varepsilon}
\def\aa{\alpha}
\def\da{\downarrow}
\def\teee{{\tilde{\eee}}}
\def\PD{D}
\def\hf{\hfill}

\def\om {\omega}
\def\la {\lambda}
\def \RA {\rightarrow}

\def\xon {x_1,x_2,\ldots ,x_n}

\def \sas {\vskip .06truein}
\def\sa{{\vskip .125truein}}

\def\sap{{\vskip .25truein}}

\def \eee {\epsilon}
\def\aaa {\alpha}
\def\bbb {\beta}

\def\ggg {\gamma}
\def\aa {\alpha}

\def\con {\subseteq}
\def \ses {\enskip = \enskip}
\def \sps {\enskip + \enskip}
\def \sms {\enskip -\enskip}

\def \scs {\ssp , \ssp}
\def \ess {\enskip}
\def \ssp {\hskip .25em}
\def \bigsp {\hskip .2truein}
\def \part {\vdash}
\def\P#1{{\Phi^{(#1)}}}

\def \BS {{\bf \Sigma}}
\def \tx {{\tilde x}}
\def \tmu {{\tilde \mu}}
\def \BA {{\bf A}} \def \BB{{\bf B}}\def \BC{{\bf C}}
\def \GGG {\Gamma}

\newdimen\Squaresize \Squaresize=14pt
\newdimen\Thickness \Thickness=0.5pt

\def\Square#1{\hbox{\vrule width \Thickness
   \vbox to \Squaresize{\hrule height \Thickness\vss
      \hbox to \Squaresize{\hss#1\hss}
   \vss\hrule height\Thickness}
\unskip\vrule width \Thickness}
\kern-\Thickness}

\def\Vsquare#1{\vbox{\Square{$#1$}}\kern-\Thickness}

\def\Young#1{
\vbox{\smallskip\offinterlineskip
\halign{&\Vsquare{##}\cr #1}}}

\newdimen\squaresize \squaresize=3pt
\newdimen\thickness \thickness=0.2pt

\def\square#1{\hbox{\vrule width \thickness
   \vbox to \squaresize{\hrule height \thickness\vss
      \hbox to \squaresize{\hss#1\hss}
   \vss\hrule height\thickness}
\unskip\vrule width \thickness}
\kern-\thickness}

\def\vsquare#1{\vbox{\square{$#1$}}\kern-\thickness}

\def\thisbox#1{\kern-.09ex\fbox{#1}}
\def\downbox#1{\lower1.200em\hbox{#1}}
\enumerate
butler chang BBGHM bergeronhamel bergerongarsia garsiahaiman garsiahaimanorbit
garsiahaimanpieri garsiahaimancatalan garsiahaimankostka garsiaprocesi grasiaatesler garsiaremmel
garsiateslerhaiman haiman kirillovnoumi knop knopquantum lapointevinet
macdonald macdonaldbook reiner sahi youngsubst
\endenumerate
\raggedbottom
%%%%%%%%%%%%%%%%%%%%%%%%%%%%%%%%%%%%%%%%%%%%%%%%%%%%%%%%%%%%%%%%%%%%%%%
%
%  Text starts here
%
%%%%%%%%%%%%%%%%%%%%%%%%%%%%%%%%%%%%%%%%%%%%%%%%%%%%%%%%%%%%%%%%%%%%%%%%%
%\hsize 6truein
%\vsize=8truein
 
\null\bigskip
\centerline{\super Science Fiction}
\medskip\centerline{\gros and }\medskip
\centerline{\super Macdonald's Polynomials}
\bigskip
\centerline{\Sc F. Bergeron\footnote{$\ ^{(\dag)}$}{With support from NSERC} and A. M. Garsia
\footnote{$\ ^{(\dag\dag)}$} {With support from NSF}} 

\abstract{
This work studies the remarkable relationships that hold among certain
$m$-tuples of the Garsia-Haiman modules $ \BM_\mu$ 
and  corresponding elements of the Macdonald basis. 
We recall that in [\garsiahaimankostka], $\BM_\mu$ is defined for a partition 
$\mu\part n$, as the linear span of derivatives of a
certain   bihomogeneous polynomial $  \Delta_ \mu(x,y)$ in the variables
$  x_1,x_2,\ldots ,x_n$, $  y_1,y_2,\ldots ,y_n$.   It has been conjectured in [\garsiahaiman],
[\garsiahaimankostka] that
$ \BM_\mu$ has $  n!$ dimensions and that its bigraded Frobenius characteristic
is given by the symmetric polynomial $  \TH_\mu(x;q,t)=\sum_{\lambda\part n }
S_\lambda(X)\TK_{\lambda\mu}(q,t)$ where  the $  \TK_{\lambda\mu}(q,t)$ are related to the
Macdonald $   q,t$-Kostka coefficients $   K_{\lambda\mu}(q,t)$ by the identity 
$   \TK_{\lambda\mu}(q,t)=K_{\lambda\mu}(q,1/t)t^{n(\mu )}$ with $   n(\mu)$ the 
$   x$-degree of $   \Delta_ \mu(x;y)$. Using this conjectured relation  we can
translate observed or proved properties of the modules $\BM_\mu$ into
identities for Macdonald polynomials. 
Computer data has suggested that as
$   \nu$ varies among the immediate predecessors of a partition $   \mu$, the spaces
$   \BM_\nu$ behave like a boolean lattice. The same appears to holds true when
$   \nu$ varies among the immediate successors of $   \mu$. 
Combining this property with a number of observed facts and some 
``heuristics'' we have been led to formulate a number of remarkable conjectures 
about the Macdonald polynomials. In particular we obtain
a representation theoretical interpretation for some of the symmetries that can be found in the
computed tables of $q,t$-Kostka coefficients.
The expression ``Science Fiction'' here refers to a package of ``heuristics'' that we use
to describe relations amongst the modules $\BM_\nu$. 
These heuristics are purely speculative assertions that are used as a convenient guide 
to the construction of identities relating the corresponding bigraded Frobenius characteristics.
Nevertheless computer experimentation reveals that these assertions are ``generically''
correct. Moreover, the evidence in support of the symmetric function identities we have derived from them 
are overwhelming. In particular, we show that various independent consequences of our heuristics 
lead to the same final identities.  
}

\section{Introduction}
Throughout this writing $\mu$ will be a partition of $n$ (denoted $\mu\vdash n$) and $\mu'$ will denote its
conjugate. We shall use the French convention here and, 
given that the parts of $\mu$ are $\mu_1\geq \mu_2\geq \cdots\geq \mu_k>0$, we let  
the corresponding Ferrer's diagram have  $\mu_i$ lattice squares in the $i^{th}$ row
(counting from the bottom up), hence the diagram of $4311$ is:

  $$ \Young{ \cr
            \cr
            & & \cr
            & & & \cr}
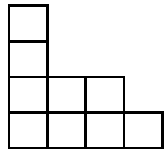
$$

\setbox\squ=\hbox{\vrule width.2pt 
   \vbox{\hrule height.2pt width.3em\kern.7ex 
         \hrule height.2pt}%
   \vrule width.2pt depth0pt }
\noindent
As in Macdonald [\macdonald], for $\mu=(\mu_1\geq\mu_2\geq \cdots \geq
\mu_k>0)$, we let  
$$
n(\mu):= \sum_{i=1}^k (i-1) \mu_i,
\eqno{\rm I}.1
$$
We also let the
coordinates $(i,j)$ of a cell $c\in\mu$ respectively measure the height of $c$ and the
position of $c$ in its row. 

We shall also adopt the Macdonald convention of calling 
the {\it arm, leg, co-arm\/} and {\it co-leg\/} of a lattice square $s$  the parameters
$a(s)$, $l(s)$, $a'(s)$ and $l'(s)$ giving the number of cells of $\mu$ that are
respectively {\it strictly\/} {\Sc east, north, west} and { \Sc south} of $s$ in $\mu$.
Recall that Macdonald in [\macdonald] established the existence of a symmetric function basis
$\{P_\mu(x;q,t)\}_\mu$ uniquely characterized by the following conditions
  $$\eqalign{
    &{\rm a)} \quad P_\lambda  = S_\lambda + \sum_{\mu < \lambda} S_\mu  \xi_{\mu \lambda} (q,t), \cr
    &{\rm b)} \quad \scalar{P_\lambda,P_\mu} _{q,t} = 0 \quad {\rm for}\quad  \lambda \neq \mu,\cr}
    \eqno{\rm I}.2 
$$
where $\scalar{,}_{q,t}$ denotes the scalar product of symmetric polynomials defined by setting
for the power basis $\{ p_\mu \}$ 
$$
\scalar{p_{\mu},p_{\lambda}}_{q,t} = \cases { z_\mu  p_\mu 
  \Bigl[{1-q \over 1-t}\Bigr]
   & if $ \mu=\lambda $, \cr\cr
  0 & otherwise. \cr}
$$
Here, as customary, $z_\mu$ denotes the integer that makes
$n!/z_\mu$ the number of permutations with cycle structure $\mu$.
We have also used  $\lambda$-ring notation in I.3.We recall
that $\lambda$-substitution in
a symmetric function is the linear and multiplicative extension of a substitution
explicitly defined for the power sums. Where, for any given expression $E$ and any
integer $k\geq 1$, we let $p_k[E]$ be the expression obtained from $E$
by replacing all the variables occurring in $E$ by their $k^{\rm th}$-powers.  For example, writing
$X$ for $x_1+x_2+x_3+\ldots$, the substitution of ${X\over 1-t}$
in a given symmetric function $f$, henceforth denoted $f[{X\over 1-t}]$, corresponds to the
linear and multiplicative extention of
$$
p_k[X]=x_1^k+x_2^k+x_3^k+\ldots\  \mapsto\ 
            p_k[{\hbox{$\textstyle X\over 1-t$}}] = {1\over 1-t^k}\ess p_k[X].
$$
 There
are a number of outstanding conjectures concerning the polynomials $ P_\lambda$ (see
[\macdonaldbook]). We are dealing  here with those involving  the so called {\sl
integral forms\/} $J_\mu(x;q,t)$ and their associated Macdonald-Kostka
coefficients $K_{\lambda\mu}(q,t)$. We shall use the same notation as in
[\macdonaldbook]. In particular
$\{ Q_\lambda(x;q,t)\}$ denotes the basis dual to $\{ P_\lambda(x;q,t)\}$ with respect to 
the scalar product $\scalar{,}_{q,t}$. Clearly, I.2 b)  gives
$$ 
Q_\lambda(x;q,t) = d_\lambda (q,t)   P_\lambda(x;q,t) ,
  \eqno{\rm I}.3 
$$
for a suitable rational function $d_\lambda (q,t)$. However in [\macdonaldbook], it is shown that
$$ 
d_\lambda (q,t) = {h_\lambda (q,t)\over h'_\lambda (q,t)} 
\eqno{\rm I}.4
$$
with
$$
h_\lambda (q,t):= \prod _{s\in \lambda }(1-q^{a_\lambda (s)} t^{l_\lambda (s) +1}),\qquad
    h'_\lambda (q,t):= \prod _{s\in \lambda }(1-q^{a_\lambda (s)+1 } t^{l_\lambda (s) }), 
$$
where $s$ denotes a generic lattice square and $a_\lambda (s)$, $l_\lambda (s)$  respectively
denote the {\it arm\/} and the {\it leg\/} of $s$ in the Ferrers' diagram of
$\lambda $. 

From I.3 and I.4 it follows that we may set, as Macdonald does in [\macdonaldbook]:
$$ 
J_\mu (x;q,t) = h_\mu(q,t) P_\mu(x;q,t) = h'_\mu(q,t) Q_\mu(x;q,t)\,.
       \eqno{\rm I}.5 
$$
The coefficients $K_{\lambda\mu}(q,t)$ can now be defined through an expansion which, in
$\lambda$-ring notation, may be written as
$$
 J_\mu (x;q,t) = \sum_\lambda S_\lambda [X (1-t)]  K_{\lambda\mu}(q,t). \eqno{\rm I}.6 
$$
Macdonald conjectured that these coefficients are polynomials in $q$ and $t$ with
non-negative integer coefficients. There are now a number of proofs of integrality
([\grasiaatesler], [\garsiaremmel], [\kirillovnoumi], [\knop], [\sahi]) but the positivity has still
to be demonstrated.

We define
 $$ H_\mu (x;q,t):= \sum_\lambda  S_\lambda  K_{\lambda\mu}(q,t) = J_\mu[X/(1-t);q,t],
$$
and set
 $$ \TH_\mu (x;q,t):= H_\mu (x;q,1/t) t^{n(\mu)}. 
  \eqno{\rm I}.7 
$$
Note that we also have
 $$ \TH_\mu (x;q,t):= \sum_\lambda  S_\lambda  \widetilde{K}_{\lambda\mu}(q,t),
\eqno{\rm I}.8 
$$
with 
  $$\widetilde{K}_{\lambda\mu}(q,t)=t^{n(\mu)}K_{\lambda\mu}(q,1/t).$$

We shall set 
  $$ B_\mu(q,t):= \sum_{(i,j)\in \mu} t^{i-1}q^{j-1},\qquad 
  T_\mu(q,t):= \prod_{(i,j)\in \mu} t^{i-1}q^{j-1}=q^{n(\mu')}t^{n(\mu)}.
  \eqno{\rm I}.9 $$
The pairs ${(i-1,j-1)}$ occurring in the above sum are briefly referred to as the
{\it biexponents\/} of $\mu$ and $B_\mu(q,t)$ itself will be called the {\sl
biexponent generator\/} of
$\mu$.

 Now let $(p_1,q_1),\ldots , (p_n,q_n)$ denote the set of biexponents arranged in
lexicographic order and set
  $$ \Delta_ \mu(x,y):=\Delta_ \mu (x_1,\ldots ,x_n;y_1,\ldots ,y_n)=
   \det \| x_i^{p_j}y_i^{q_j}\|_{i,j=1..n}. \eqno{\rm I}.10 $$ 
For example
$$ \Delta_{\trois321}=\det\!\pmatrix{ 1 & x_{1} & {x_{1}}^2 & y_{1} & x_{1}\,y_{1}
 & {y_{1}}^2 \cr 1 & x_{2} & {x_{2}}^2 & y_{2} & x_{2}\,y_{2} & {y_{2}
}^2 \cr 1 & x_{3} & {x_{3}}^2 & y_{3} & x_{3}\,y_{3} & {y_{3}}^2 \cr 1
 & x_{4} & {x_{4}}^2 & y_{4} & x_{4}\,y_{4} & {y_{4}}^2 \cr 1 & x_{5}
 & {x_{5}}^2 & y_{5} & x_{5}\,y_{5} & {y_{5}}^2 \cr 1 & x_{6} & {x_{6}
}^2 & y_{6} & x_{6}\,y_{6} & {y_{6}}^2 \cr}\eqno{\rm I}.11$$
This given we let $\BM_\mu$ be the collection of polynomials in the variables 
$x_1,\ldots ,x_n;y_1,\ldots ,y_n $ obtained by taking the linear span of all the
partial derivatives of $\Delta_ \mu$. In symbols we may write
  $$ \BM_\mu = {\cal L}[\del_x^p\del_y^q \Delta_ \mu (x;y)]
   \eqno{\rm I}.12 $$
where $\del_x^p=\del_{x_1}^{p_1}\cdots \del_{x_n}^{p_n}$,
$\del_y^q=\del_{y_1}^{q_1}\cdots \del_{y_n}^{q_n}$ .
\sas

The natural action of a permutation $\sigma=(\sigma_1,\ldots , \sigma_n)$ on a polynomial 
$P(x_1,\ldots ,x_n;y_1, \ldots ,y_n)$ is the so called {\it diagonal
action\/} which is defined by setting
 $$ \sigma P(x_1,\ldots,x_n ;y_1,\ldots,y_n ):= P(x_{\sigma_1},\ldots,x_{\sigma_n};y_{\sigma_1}, \ldots
  ,y_{\sigma_n}).$$
Since $\sigma \Delta_ \mu= \pm \Delta_ \mu$ according to the sign of $\sigma$, the space 
$\BM_\mu$ necessarily remains invariant under this action.  
\sas

Note that, since $\Delta_ \mu$ is bihomogeneous of degree $n(\mu)$ in $x$ and $n(\mu ')$ in
$y$, we also have the direct sum decomposition 
$$
\BM_\mu=\bigoplus _{h=0}^{n(\mu )}\bigoplus _{k=0}^{n(\mu ')}  \CH_{h,k}(\BM_\mu),
$$
where $ \CH_{h,k}(\BM_\mu)$ denotes the subspace of $\BM_\mu$ 
spanned by its bihomogeneous elements of degree $h$ in $x$ and degree $k$ in $y$.
Since the diagonal action clearly preserves bidegree, each of the 
subspaces $\CH_{h,k}(\BM_\mu)$ is also $S_n$-invariant.
Thus we see that $\BM_\mu$ has the structure of a bigraded module. The generating
function of the characters of its bihomogeneous components, which we shall refer to
as the {\it bigraded character\/} of $\BM_\mu$, may be written in the form
  $$\charac^\mu (q,t) = \sum_{h=0}^{n(\mu )}\sum _{k=0}^{n(\mu ' )} t^h q^k \char
           \CH_{h,k}(\BM_\mu).$$
The bivariate Frobenius characteristic of $\BM_\mu$ is
$$
{\cal F}(\BM_\mu) :={1\over n!} \sum_{\sigma\in S_n} \charac^{\mu} (\sigma;q,t)
p_{\lambda(\sigma)}(x) .
\eqno{\rm I}.13 
$$
where $\charac^\mu (\sigma;q,t)$ denotes the value of $\charac^\mu (q,t)$ at $\sigma$
and $p_{\lambda(\sigma)}(x)$ is the power symmetric function indexed by the shape of
$\sigma$. It will be convenient to set 
$$
{\cal F}(\BM_\mu)\ses C_\mu(x;q,t)\ses
\sum_{\lambda\part n}S_\lambda(x) 
C_{\lambda\mu}(q,t)\ess .
\eqno{\rm I}.14
$$
Where, $C_{\lambda\mu}(q,t)$ is the polynomial which gives the  occurrences of the
irreducible character $\chi^\la$ in the various bihomogeneus components of $\BM_\mu$.
More precisely the coefficient of $t^hq^k$ in $C_{\lambda\mu}(q,t)$  gives the
multiplicity of $\chi^\la$ in the submodule $\CH_{h,k}(\BM_\mu)$.
It has been conjectured in [\garsiahaiman] and
supported in [\garsiahaimanpieri], [\garsiahaimankostka] and [\reiner] (by verification of special
cases and several  representation theoretical considerations),
that
$$
C_\mu(x;q,t)\ses \TH_\mu(x;q,t)
\eqno{\rm I}.15
$$
This forces the equality 
$$
C_{\lambda\mu}(q,t)=\TK_{\lambda\mu}(q,t)
\eqno{\rm I}.16
$$
which in particular
implies the Macdonald conjecture that the $K_{\lambda\mu}(q,t)$ are polynomials
with non-negative integer coefficients. We shall briefly refer to I.15 as the $C=\TH$ conjecture.
\sas

Macdonald in [\macdonaldbook] derives a number of properties of the $K_{\lambda\mu}(q,t)$, in particular he
shows that for any partition $\mu$ 
$$ 
K_{\lambda\mu}(1,1)= f_\lambda ,
$$
where as customary $f_\lambda$ denotes the number of standard tableaux of shape $\lambda$.
Thus the validity of the $C=\TH$ conjecture requires that  $\BM_\mu$ should
yield a bigraded version of the left regular representation. Then a fortiori we should have
$$
\dim \BM_\mu = n!
\eqno{\rm I}.17
$$
This has come to be referred to as the {\it $n$-factorial Conjecture} (see [\garsiahaimankostka]).
In the nearly six years since this deceptively simple conjecture was
formulated, an overwhelming amount of theoretical and experimental evidence in its support has
been gathered. But perhaps the most surprising development in this context is a recent
work of M. Haiman [\haiman] where it is shown that the  $n$-factorial conjecture for a given $\mu$
is all that is needed to establish the identity $C_\mu(x;q,t)=\TH_\mu(x;q,t)$
for that same $\mu$. 
It is easy to see that, for $\mu=1^n$ and $\mu=(n)$,
$\Delta _\mu$ reduces to the Vandermonde determinant in $x$ and $y$ respectively. 
In these cases (I.15) is a classical result (see [\garsiahaimanorbit]).
\sas

Experimental evidence has revealed that the  modules $\BM_\nu$ corresponding
to  partitions $\nu$ lying  immediately below (resp. above) a given partition $\mu$
have remarkable intersection properties. Our purpose here is to explore various
results and conjectures that can be obtained  by combining these findings
with the $C=\TH$ conjecture. 
\sas

Let us begin by observing that a space $\BM$ spanned by the partial derivatives of a single 
bihomogeneous polynomial
$\Delta(x,y)=\Delta(x_1,\ldots,x_n;y_1,\ldots,y_n)$ of bidegree $(a,b)$, in symbols
$$
\BM={\cal L}[\del_x^p\del_y^q\,\Delta(x,y)],
$$
has a bidegree complementing automorphism, called {\bf flip}, defined by setting, for every $P\in\BM$,
$$
{\bf flip}\, P:=P(\del_x,\del_y)\,\Delta(x,y)=P(\del)\,\Delta(x,y),
$$
where here and  after, if $P$ is a polynomial in $x_1,\ldots,x_n;y_1,\ldots,y_n$, $P(\del)$ is to represent the differential
operator obtained by replacing in $P$, $x_i$ by $\del_{x_i}$ and $y_i$ by $\del_{y_i}$. 

Note that, if $\Delta$ is alternating
then
$\BM$ is invariant under the diagonal action of
$S_n$, and moreover the effect of flip will be to sign-twist each of the irreducible submodules.  Hence, 
if $\Phi(x;q,t)$ is the  bivariate Frobenius characteristic of $\BM$, then we necessarily have that
$$
\Phi(x;q,t)=t^aq^b\omega\, \Phi(x;q^{-1},t^{-1}).\eqno{\rm I}.18
$$
Here, as usual, $\omega$ denotes the involution that sends $S_\lambda$ to $S_{\lambda'}$. 
More generally,
for any bigraded submodule ${\bf M_1}$ of $\BM$, the subspace 
$$
\flp\, {\bf M_1}:={\cal L}[\flp\, P\ |\  P\in{\bf M_1}],\eqno{\rm I}.19
$$
is also a bigraded submodule of $\BM$. Moreover, 
if  $\Psi(x;q,t)$ is the bivariate Frobenius characteristic of ${\bf M}$, then the  bivariate
Frobenius characteristic of $\flp\,{\bf M_1}$ is given by the formula
$$
{\cal F}(\flp\,{\bf M_1})=t^aq^b\, \da \Psi(x;q,t),\eqno{\rm I}.20
$$
where we set
$$
\da \Psi(x;q,t):=\omega \Psi(x;q^{-1},t^{-1}).\eqno{\rm I}.21
$$

Here and after, we shall assume that our
partition $\mu$ has $m$ corners, and that the set of predecessors of $\mu$
is 
$$
\pi(\mu)=\{\aaa^{(1)},\aaa^{(2)},\ldots, \aaa^{(m)}\},
$$ 
with $\aaa^{(i)}$ the partition obtained by removing the $i^{\rm th}$-corner  
(as we encounter it from left to right). 
Our computer experimentations indicate that the space ${\bf V}_\mu$,  sum  of the $\BM_{\aaa^{(i)}}$'s, in symbols
$$
{\bf V}_\mu:=\bigvee_{i=1}^m \BM_{\aaa^{(i)}},\eqno{\rm I}.22
$$
has a  basis ${\cal B}={\cal B}(\mu)$ , of bihomogeneous polynomials, 
with the properties (i)--(v) given below.

\itemitem{(i)}{\it For each $1\leq i\leq m$, we have a subset ${\cal B}_i\subset {\cal B}$ 
such that\/}
$$
\BM_{\aaa^{(i)}}={\cal L}[{\cal B}_i], \eqno{\rm I}.23
$$

\noindent For a word $\eee=\eee_1 \eee_2\cdots \eee_m$, with $\eee_i=0,1$, set
$$
{\cal B}^\eee:=\bigcap_{i=1}^m {\cal B}_i^{\eee_i},
$$
where, for a given subset $S\subset {\cal B}$, we let 
$$ 
S^1:=S\,,\qquad S^0:= {\cal B}\setminus S\, .
$$
This given,

\itemitem{(ii)} {\it For each $\eee$ the space

$$
\BM^\eee_\mu:={\cal L}[{\cal B}^\eee_\mu],\eqno{\rm I}.24
$$
is $S_n$-invariant, and therefore has the structure of a bigraded $S_n$-module, whose bivariate Frobenius
characteristic we denote $\Phi^\eee=\Phi^\eee_\mu$, in symbols\/}
$$
\Phi^\eee_\mu(x;q,t):={\cal F}(\BM^\eee_\mu).\eqno{\rm I}.25
$$

\noindent Note that if we set 
$$
\BM_{\aaa^{(i)}}^1:={\cal L}[{\cal B}_i]=\BM_{\aaa^{(i)}}\,,
\qquad \BM_{\aaa^{(i)}}^0:= {\cal L}[{\cal B}\setminus {\cal B}_i]\, ,
$$
then we necessarily have
$$
\BM_\mu^\eee=\bigcap_{i=1}^m \BM_{\aaa^{(i)}}^{\eee_i}. \eqno{\rm I}.26
$$

\noindent Let now, ${\bf flip}_i$ denote the flip operation within the space $\BM_{\aaa^{(i)}}$, and let  $\da_i$
denote its effect on the bivariate Frobenius characteristics of the bigraded submodules of $\BM_{\aaa^{(i)}}$.
According to (I.20) and (I.21), if $\Psi(x;q,t)$ is such a characteristic, then we have
$$
\da_i \Psi(x;q,t)=T_i\, \da \Psi(x;q,t),\eqno{\rm I}.27
$$
where hereafter we simply use $T_i$ for $T_{\aa_i}(q,t)$.

Let $\tau_i$ denote the operator on words in $0,1$ which complements all but the $i^{\rm th}$ letter.
That is, setting $\teee_j:=1-\eee_j$, we have 
$$
\tau_i\bigl(  \eee_1\cdots\eee_{i-1}\eee_i\eee_{i+1}\cdots \eee_m\bigr) :=
\teee_1\cdots\teee_{i-1}\eee_i\teee_{i+1}\cdots \teee_m.\eqno{\rm I}.28
$$
This given, the most remarkable property of our flip operations can be stated as follows
 
\itemitem{(iii)} {\it  For all words $\eee=\eee_1\eee_2\cdots\eee_m$ such that $\eee_i=1$,
we have\/}
$$
(a)\ess\ess {\bf flip}_i \, \BM_\mu^{\eee}\cong
             \BM_\mu^{\tau_i\,\eee}.
\ess\ess\ess\ess {\rm (b)}\ess\ess\ess\ess 
 \BM_\mu^{\tau_i\,\eee}  \cap \ssp {\bf flip}_i \, \BM_\mu^{\eee}\ses \{0\}
\eqno{\rm I}.29
$$

\noindent Note that I.29 (a) can be  written more explicitly as
$$
\eqalign{
{\bf flip}_i \bigl( \BM_{\aaa^{(1)}}^{\eee_1}\cap \cdots\cap
\BM_{\aaa^{(1)}}^{\eee_{i-1}} &
\cap{\BM}_{\aaa^{(i)}}\cap
\BM_{\aaa^{(1)}}^{\eee_{i+1}}
\cap \cdots\cap\BM_{\alpha_{m}}^{\eee_{m}}\bigr)
\cong\cr
&\cong \ess\BM_{\aaa^{(1)}}^{\teee_1}\cap\cdots\cap 
\BM_{\aaa^{(1)}}^{\teee_{i-1}}\cap
{\BM}_{\aaa^{(i)}}\cap
\BM_{\aaa^{(1)}}^{\teee_{i+1}}\cap
\cdots\cap
\BM_{\alpha_{m}}^{\teee_{m}}\ess .\cr }
$$
Furthermore from (I.27) we derive that the bivariate Frobenius 
characteristics $\Phi_\mu^\eee(x;q,t)$ are related by the following identity
$$
\da_i \Phi_\mu^\eee=\Phi_\mu^{\tau_i\,\eee}.\eqno{\rm I}.30
$$
\sas

There is a natural scalar product $\LLL\scs \RRR$ in the space of polynomials
in $\xon;\yon$ which is defined by setting 
$$
\LLL P\scs Q\RRR\ses P(\del_x,\del_y)Q(x,y)\ssp |_{x=y=0}\ess .
\eqno{\rm I}.31
$$
It can be shown that the orthogonal complement of the space $\BM_\mu$ consists
of the ideal
$$
\CJ_\mu\ses \{\ssp P(x,y)\ssp :\ssp P(\del_x,\del_y)\ssp \Delta_\mu(x,y)=0\ssp \}\ess .
\eqno{\rm I}.32
$$
Our examples suggest thatwe should be able to construct the basis $\CB_\mu$ so
that
\sas

\itemitem{(iv)} {\ita If $\eee$ and $\eta$ are such that $\eee_i=1$ and $\eta_i=0$ then
$$
\BM_\mu^{\eta}\ess \con\ess \bigl(\ssp \BM_\mu^{\eee}\ssp \bigr)^\perp
\eqno{\rm I}.33
$$
where the symbol ``$\perp$'' is to represent the operation of taking orthogonal
complements with respect to the scalar product in I.31.}

%\itemitem{(iv)} {\ita For each $i=1,\ldots ,m$ the collection $\CB/\CB_i$ is
%contained in $\CJ_{\aa^{(i)}}$. In particular, when $\eee_i=0$ we must necessarily have 
%$$
%\BM_\mu^{\eee_i} \ssp \con \ssp \BM_i^\perp
%\eqno{\rm I}.33
%$$
\sa

\noindent Finally, one further property that is suggested by our computer experimentations 
is that the
intersection of any $k$ of the modules $\BM_{\aaa^{(i)}}$ always turns out to have 
dimension $1/k^{\rm th}$ of the common
dimension of the $\BM_{\aaa^{(i)}}$'s. More precisely:

\itemitem{(v)} {\it  For any $\mu\vdash n+1$, with $m$ parts, and for 
any $k$-subset $S$ of $\{1,2,\ldots,m\}$ 
\/}
$$
\dim \bigcap_{i\in S}\, \BM_{\aaa^{(i)}}={n!\over k}.
\eqno {\rm I}.34
$$

This given, here and after we shall work under the following basic 
\sas

\noindent{\bol Heuristic:}
{\ita  For every partition $\mu$, a basis ${\cal B}_\mu$ exists with
properties {\rm (i)--(v)}}  above. 

We shall refer to this assertion  as the ``Science fiction heuristic'' or
briefly ``SF''.  We do not venture calling this a {\ita conjecture} since it is only generically true
and  was constructed from data obtained from
experimentations with relatively small partitions.Nevertheless we shall see that, by combining SF
with the $C=\TH$ conjecture, we can predict several
surprising  identities for the polynomials $\TH_\mu(x;q,t)$ and the coefficients
$\TK_{\la\mu}(q,t)$. Remarkably, these identities can actually be established
outright in some cases and in others they can be verified numerically for significantly
large examples. We should note that tables of the $\TK_{\la\mu}(q,t)$ are now available up
to $n=12$.    
\sas

For example, note that the definition (I.25) gives that $\Phi_\mu^\eee$ is a bivariate Frobenius 
characteristic. Therefore it has a Schur function expansion
$$
\Phi_\mu^\eee(x;q,t)=\sum_{\lambda} S_\lambda\,\varphi_{\lambda\mu}^\eee(q,t),
\eqno{{\rm I}.35}
$$
with the $\varphi_{\lambda\mu}$ polynomials in $q,t$ with positive integer coefficients, a property 
of symmetric 
polynomials which here and after will by briefly referred to as ``Schur positivity''.
\sas

It develops further that the symmetric functions $\Phi_\mu^\eee$ depend on $\eee$ in a very simple 
manner. In fact, using our heuristic we will show that, for each $\mu\part n$ with $m$ corners,
we can construct $m$ Schur positive symmetric functions 
$$
\Phi^{(1)}(x;q,t)\scs \Phi^{(2)}(x;q,t)\scs \dots\scs \Phi^{(m)}(x;q,t)
\eqno{{\rm I}.36}
$$
yielding (for $\eee_1+\eee_2+\cdots+\eee_n=k$)
$$
\Phi_\mu^\eee\ses {\Phi^{(k)}(x;q,t)\over \prod_{\eee_i=0} T_i}
\eqno{{\rm I}.37}
$$
One of the most remarkable consequences of I.37 is that,
combined with the $C=\TH$ conjecture, it yields that the Macdonald polynomial $\TH_{\aaa^{(i)}}$
has the expansion 
$$
\TH_{\aaa^{(i)}}(x;q,t)\ses \sum_{k=1}^m\ssp \Phi^{(k)}\ess
e_{m-k}\Bigl[\ssp {1\over T_1}+{1\over T_2}+\cdots +{1\over T_m}\sms {1\over T_i}\ssp \Bigr ]\ess .
\eqno{{\rm I}.38}
$$
Equating coefficients of $S_\la$ we can then further derive that
$$
\TK_{\la,\aaa^{(i)}}(q,t)\ses \sum_{k=1}^m\ssp \Phi_{\la,k}^{(k)}(q,t)\ess
e_{m-k}\Bigl[\ssp {1\over T_1}+{1\over T_2}+\cdots +{1\over T_m}\sms {1\over T_i}\ssp \Bigr ]\ess .
\eqno{{\rm I}.39}
$$
where, (assuming SF), the $\Phi_{\la,k}(q,t)$ must be polynomials with positive integer coefficients.
In this manner, we obtain
a representation theoretical explanation for many
properties of the Kostka-Macdonald tables which have been previously observed by examination of data. 
In particular,  this brings to light that, when  $\lambda$ remains fixed, the dependence
of $\widetilde{K}_{\lambda\mu}(q,t)$ on $\mu$ is governed by a yet finer mechanism than the one that
was noticed by L. Butler 
\footnote{(\dag)} {Personal communication.}.

It is appropriate to review Butler's observations here. 
Let $\mu$ and $\nu$ be partitions of the same number. Following A. Young it is customary to say
that $\mu$ is obtained from $\nu$ by a ``Raising Operator'' if $\mu$ can be obtained by
lifting a number of cells of $\nu$ from lower to upper rows 
\footnote{(\dag\dag)}{This definition requires that diagrams of partitions be drawn by the english convention}. 
A raising operator will be called minimal if it
lifts a single cell up a single row or to the right a single column. It can be shown (see [\knop])
that the dominance partial order is the transitive closure of the relation $\mu=R\nu$
with $R$ a minimal raising operator. L. Butler [\butler] observed that the Macdonald coefficients
$\TK_{\la\mu}(q,t)$ change in a remarkably simple way when $\mu$ is obtained from $\nu$ by 
a minimal raising operator. Her observation which was originally made from tables 
computed by  hand by Macdonald [\macdonaldbook] (up to $n\leq 6$) is now confirmed by more extensive tables
obtained by computer (up to $n\leq 12$). Butler noticed that for any such pair of partitions
$\mu,\nu$ and any $\la$ there are two polynomials $\Phi_\la$ and $\Theta_\la$ with non-negative
integer coefficients such that
$$
\eqalign{
\TK_{\la\mu}(q,t)\ses \Phi_\la\sps \Theta_\la\cr
\TK_{\la\nu}(q,t)\ses \Phi_\la\sps \Theta_\la\ssp T_\nu/T_\mu ,\cr}
\eqno {\rm I}.40
$$
This has come to be referred to as the {\ita Butler conjecture}. 
To explore the full implications of this conjecture let us set for a moment $\psi_\la=\Theta_\la/ T_\mu$
and write these two identities in the more symmetric form 
$$
\eqalign{
\TK_{\la\mu}(q,t)\ses \Phi_\la\sps \Psi_\la\ssp T_\mu\cr
\TK_{\la\nu}(q,t)\ses \Phi_\la\sps \Psi_\la\ssp T_\nu ,\cr}
$$
Equivalently we may write
$$
\Phi_\la={T_\nu\ssp \TK_{\la\mu}\sms T_\mu\ssp \TK_{\la\nu} \over T_\nu \sms T_\mu}\ess ,
\bigsp
\Psi_\la={ \TK_{\la\mu}\sms  \TK_{\la\nu} \over T_\mu \sms T_\nu}\ess .
$$
Now it may be shown (see [\garsiahaimanorbit]) that for any $\la,\mu\part n$ we have 
$$
\TK_{\la\mu}(1/q,1/t)\ess T_\mu = \TK_{\la'\mu}(q,t)\ess , 
\bigsp
\TK_{\la\nu}(1/q,1/t)\ess T_\nu = \TK_{\la'\nu}(q,t)\ess . 
$$
These identities imply that $\Theta_\la$ is superfluous. In fact, 
combinining them with I.40 we can easily derive that
$$
\Psi_\la(q,t)\ses \Phi_{\la'}(1/q,1/t)\ess .
$$
Setting $\Phi(x;q,t)=\sum_\la\ssp \Phi_\la(q,t) \ssp S_\la(x)$ the identity in I.8
written for $\mu$ and $\nu$ yields that 
$$
\TH_\mu(x;q,t)= \Phi(x;q,t)+ T_\mu\ssp \da  \Phi(x;q,t)\ess ,\ess\ess\ess
\TH_\nu(x;q,t)= \Phi(x;q,t)+ T_\nu\ssp \da   \Phi(x;q,t)\ess .
\eqno {\rm I}.41
$$
Our findings in the case of the pair $\mu=31,\ \nu=22$ (which is treated in the next section)
reveals that the phenomenon observed by Butler should be the result of a beautiful mechanism which we venture
to state as follows.

\Conj{1.1}
If $\mu=R\nu$ with $R$ a minimal raising operator then the identities in I.41
hold true with $\Phi(x;q,t)$ the bigraded Frobenius characteristic of $\BM_\mu\wedge \BM_\nu$

On the validity of our $C=\TH$ conjecture, this property of the Macdonald polynomials 
would be an immediate consequence of
\sas

\Conj{1.2}
If $\mu=R\nu$ with $R$ a minimal raising operator then
$$
\eqalign{
\BM_\mu \ses \BM_\mu\wedge \BM_\nu\ssp \oplus \ssp \flip_\mu \ssp \BM_\mu\wedge \BM_\nu\ess ,\cr
\BM_\nu \ses \BM_\mu\wedge \BM_\nu\ssp \oplus \ssp \flip_\nu \ssp \BM_\nu\wedge \BM_\nu\ess .\cr}
\eqno {\rm I}.42
$$

It not difficult to see that both these conjectures may be derived from our heuristic
even in the more general case that $\mu$ and $\nu$ are any two predecessors of the
same partition.
\sas
 
Finally, we should mention that in recent joint work M. Haiman and C. Chang have been able to
give an Algebraic Geometrical setting to some of the identities we derive here. It develops that
in their setting,
some of the ingredients we use here in a purely heuristic way 
take a remarkable natural place within the theory of Hilbert schemes.  
This permitted them to establish the validity of our assumption for any two part partition and to
construct a general mechanism for proving its validity in full generality.
\sas

The contents  of this paper are divided into four sections. In the first section we  
work out some examples in full detail. We do this to show that our SF heuristic stems right out
of a very natural approach to the proof the
$n$-factorial conjecture. In section 2, we work out the case $\mu=321$ since it is
the first case which clearly displays all the facets of the SF heuristic.  
\sas

In section 3, using $C=\TH$, we derive 
formulas for all the characteristics $\Phi_\mu^\eee$ in terms of the polynomials $\TH_\mu(x;q,t)$
and derive from them some additional positivity properties of the Macdonald basis.
We also show in section 2 that (I.15) combined with (I.34) 
yields precise
dimension counts for all the submodules $\BM_\mu^\eee$. 
\sas

In section 4 we present  applications. In particular we derive there new and refined 
versions of Pieri formulas for the polynomials $\TH_\mu(x;q,t)$. The fact that
these formulas  are completely consistent with those originally obtained by Macdonald should
provide support for the validity of at least some weakened version of the SF  heuristic.

\section{1. A recursive approach to the $n!$-conjecture}
There is an elementary approach to proving the $n!$ conjecture 
which can be successfully used in various special cases.
To get accross the basic idea we shall systematically work out in this manner the
proof of the $n!$-conjecture for all $n\leq 5$. Of course, in doing this we shall
have to take for granted some of the results obtained in previous work. However,
all the auxiliary material we need here can  be found in [\garsiahaimanorbit], [\garsiahaimancatalan] or
[\garsiahaimankostka].
\sas

Our point of departure is a result (proved in [\garsiahaimankostka]) that for
any partition $\mu$ we have
$$
\dim \BM_\mu\ssp \leq \ssp n!
\eqno 1.1
$$
This means that to establish the $n!$-factorial conjecture for a given $\mu$ we 
need only exhibit $n!$ independent elements in $\BM_\mu$.
\sas

Note next that from the definition I.10 we immediately derive that
$$
\Delta_{\mu'}(x,y)\ses \Delta_{\mu}(y,x)
$$ 
In particular the Frobenius characteristics $C_\mu(x;q,t)$ 
are related by the identity
$$
C_{\mu'}(x;q,t)\ses  C_\mu(x;t,q)\ess .
$$ 
This is consistent with the $C=\TH$-conjecture since the duality result of
Macdonald implies that
$$
\TH_{\mu'}(x;q,t)\ses  \TH_\mu(x;t,q)\ess .
$$
Thus we need  establish the $n!$-conjecture only for one member of each pair of conjugate partitions.
\sas

Let us recall that the portion of 
$0$ $y$-degree
in any of our modules $\BM_\mu$ has dimension ${n \choose \mu}$. This may be easily derived
from the basic result in [\garsiaprocesi] and
the fact that this bihomogeneous subspace of $\BM_\mu$ is simply the linear span
of derivatives of the Garnir elements of shape $\mu$. Details of this derivation
can be found in [\bergerongarsia]. Thus in the particular case that
$\mu=21^{n-2}$ this dimension is $n!/2$. Since $\bf flip$ operation in any
$\BM_\mu$ yields 
a sign-twisted $S_n$-module isomorphism of the $0$ $y$-degree portion into the $n(\mu')$
$y$-degree portion, we see that for all these partitions (and their conjugates)
the $n!$-conjecture must necessarily hold true. Moreover the case of 
$\mu=1^n$ (single columns) or $\mu=(n)$ (single rows)  is classical since 
$\Delta_\mu$ then reduces to the Vandermonde determinant in the $x's$ and the $y's$ respectively.
In particular, for $n=4$ we are only left with $\mu=22$.
We shall start by dealing with this case. 

\sa
Expanding with respect to the last column we get
$$
\Delta_{22}\ses \det 
\pmatrix{
1& 1& 1& 1\cr
x_1& x_2& x_3& x_4\cr
y_1& y_2& y_3& y_4\cr
x_1y_1& x_2y_2& x_3y_3& x_4y_4\cr}\ses
\Phi_{00}+x_4\ssp \Phi_{10}+y_4\ssp \Phi_{01}+x_4y_4\Phi_{11}\ess .
\eqno 1.2
$$
So we may write
$$
\del_{x_4}\Delta_{22}= \Phi_{10}+y_4\Phi_{11} \scs\ess
\del_{y_4}\Delta_{22}= \Phi_{01}+x_4\Phi_{11} \scs\ess
\del_{x_4}\del_{y_4}\Delta_{22}= \Phi_{11} \ess .
\eqno 1.3
$$
Note that
$$
\Phi_{11}\ses \Delta_{21}\ses
\det \pmatrix{
1& 1& 1\cr
x_1& x_2& x_3\cr
y_1& y_2& y_3\cr}
\eqno 1.4
$$
Recall that if $P$ is a polynomial in $x_1,.., x_n;y_1,.., y_n$, 
then $P(\del)=P(\del_{x_1},..,\del_{x_n};\del_{y_1},..,\del_{y_n})$. This given, let $\CB$ be a collection of monomials in
$x_1,x_2,x_3;y_1,y_2,y_3$ chosen so that $\bigl\{ b(\del)\Delta_{21}\bigr\}_{ b\in \CB}$
is a basis for $\BM_{21}$
\footnote{$(\dag)$}
{We may take $\CB=\{1,x_1,x_2,y_1,y_2,x_1y_2\}$}.
 We claim that the polynomials in the union
$$
\CB_{22}=\bigl\{b(\del)\Delta_{22}\bigr\}_{b\in \CB}+
\bigl\{ b(\del)\del_{x_4}\Delta_{22}\bigr\}_{b\in \CB}+
\bigl\{ b(\del)\del_{y_4}\Delta_{22}\bigr\}_{b\in \CB}+
\bigl\{ b(\del)\del_{x_4}\del_{y_4}\Delta_{22}\bigr\}_{b\in \CB}
\eqno 1.5
$$
are linearly independent.  Note that since from 1.1 we get that $\dim\BM_{22}\leq 24$
we can then conclude that $\CB_{22}$ is a basis for $\BM_{22}$ and that the conjecture is true for
$\mu=22$. To this end, let there be $b_{00},b_{10},b_{01},b_{11}$ in the linear span of $ \CB$ 
giving
$$
b_{00}(\del)\Delta_{22}\sps
b_{10}(\del)\del_{x_4}\Delta_{22}\sps
b_{01}(\del)\del_{y_4}\Delta_{22}\sps
b_{11}(\del)\del_{x_4}\del_{y_4}\Delta_{22}\ses 0\ess .
$$
Using 1.2, 1.3 and 1.4 and extracting the coefficients of $1,x_4,y_4,x_4y_4$ 
we derive that we must have the simultaneous equations
$$
\matrix{
b_{00}(\del)\Phi_{00}   &+&   b_{10}(\del)\Phi_{10}      &+&       b_{01}(\del)\Phi_{01}      & + &      b_{11}(\del)\Delta_{21}&=&0 \cr
                      &&   b_{00}(\del)\Phi_{10}      &&                                & + &      b_{01}(\del)\Delta_{21}&=&0 \cr
                      &&                            &&       b_{00}(\del)\Phi_{01}      & + &      b_{10}(\del)\Delta_{21}&=&0 \cr
                      &&                            &&                                &  &      b_{00}(\del)\Delta_{21}&=&0 \cr}
$$
Note now that if $b_{00}=\sum_{b\in \CB}c_b\ssp b$ then the last equation may be written as
$$
\sum_{b\in \CB}c_b\ssp b(\del)\Delta_{21}\ses 0\ess ,
$$
but this contraddicts the choice of $\CB$ unless the coefficients $c_b$ are all equal to zero. Substituting $b_{00}=0$
in the third equation reduces it to $ b_{10}(\del)\Delta_{21}=0$ which again implies that the coefficients of $b_{10}$
must all vanish as well. Substituting $b_{00}=0$ in the second equation yields the same for the coefficients of $b_{01}$.
Finally, setting $b_{00}=b_{10}=b_{01}=0$ in the first equation forces the vanishing of the coefficients of
$b_{11}$. This proves that $\CB_{22}$ is an independent set as asserted.  
\sas

We should note that our argument here establishes a bit more than the validity of the $n!$
conjecture for $\mu=22$. 
Let us recall that if $\Theta$ is the Frobenius image of an $S_n$ character $\chi$ then the
partial derivative\footnote{$\ ^{(\dag)}$}
{Here $p_1$ denotes the first power symmetric polynomial and the symbol $\del_{p_1}\Theta$ is
to mean the partial derivative of $\Theta$ as a polynomial in the power symmetric functions.} $\del_{p_1}\Theta$  
yields the Frobenius image of the restriction of
$\chi$ to $S_{n-1}$. In particular, the polynomial
$$
G_\mu(q,t)=\del_{p_1}^n\ssp C_\mu(x;q,t)\ses \sum_{\la\part n}\ssp f_\la C_{\la\mu}(q,t)
$$  
must give the bigraded Hilbert series of $\CB_\mu$. Of course under the $C=\TH$ conjecture
this Hilbert series is given by the polynomial
$$
F_\mu(q,t)\ses \del_{p_1}^n\ssp \TH_\mu(x;q,t)\ses \sum_{\la\part n}\ssp f_\la \TK_{\la\mu}(q,t)
\ess .
\eqno 1.6
$$
We should also keep in mind  (see [\garsiahaimanorbit])  that the ``duality'' result for Macdonald
polynomials in our notation becomes
$$
\om\TH_\mu(x;1/q,1/t)t^{n(\mu)}q^{n(\mu')}\ses \TH_\mu(x;q,t)
\eqno 1.7
$$
where $\om$ is the involution which sends $S_\la$ into $S_{\la'}$. Thus 1.7
implies that we must also have
$$
F_\mu(1/q,1/t)t^{n(\mu)}q^{n(\mu')}\ses F_\mu(q,t)\ess .
\eqno 1.8
$$
The analogous identities
$$
\om C_\mu(x;1/q,1/t)t^{n(\mu)}q^{n(\mu')}\ses C_\mu(x;q,t)
\eqno 1.9
$$
$$
G_\mu(1/q,1/t)t^{n(\mu)}q^{n(\mu')}\ses G_\mu(q,t)\ess .
\eqno 1.10
$$
follow from I.18.
\sas

It is easy to see that if $\beta(q,t)$ gives the bidegree distribution of the monomials in
$\CB$ then $qt\beta(1/q,1/t)$ must give the Hilbert series of $\BM_{21}$. Thus from 1.10
we get that $\beta=G_{21}(q,t)$. From our construction 1.4 of the basis $\CB_{22}$ we can then
immediately derive that the Hilbert series of $\BM_{22}$ must be given by formula
$$
G_{22}(q,t)=\beta(1/q,1/t)(1+1/t+1/q+1/tq)q^2t^2=(1+t+q+tq)G_{21}(q,t)\ess .
$$ 
A slightly more refined argument which takes account of the action of $S_3$ on the basis $\CB_{22}$,
yields that we must also have 
$$
\del_{p_1}C_{22}(x;q,t)\ses (1+t+q+tq)C_{21}(x;q,t)\ess .
\eqno 1.11
$$
\sas

It is not difficult to see that the argument we have given here can be generalized to the case of
arbitrary rectangular partitions $\mu=r^s$ and obtain that

\Thm{1.1}
If the $n!$ conjecture holds for the partition $r-1,r^{s-1}$ then it holds for $r^s$,
Moreover,
$$
\del_{p_1}C_{r^s}(x;q,t)\ses B_{r^s}(q,t)\ssp C_{r-1,r^{s-1}}(x;q,t)
\eqno 1.12
$$
$$
G_{r^s}(q,t)\ses B_{r^s}(q,t)\ssp G_{r-1,r^{s-1}}(q,t)
\eqno 1.13
$$

Continuing with our examples we see that for $n=5$ we only need to deal with the partitions
$32$ and  $311$. We start with $32$.
Expanding again with respect to the last column we get
$$
\Delta_{32}\ses \det 
\pmatrix{
1& 1& 1& 1&1\cr
y_1& y_2& y_3& y_4& y_5\cr
y_1^2& y_2^2& y_3^2& y_4^2& y_5^2\cr
x_1& x_2& x_3& x_4& x_5\cr
x_1y_1& x_2y_2& x_3y_3& x_4y_4& x_5y_5\cr
}\ses
\Phi_{00}+y_5\ssp \Phi_{01}+y_5^2\ssp \Phi_{02}+x_5\Phi_{10}+x_5y_5\Phi_{11}\ess .
\eqno 1.14
$$
So we may write
$$
\matrix{
\Delta_{32}&=&
 b_{00}(\del)\Phi_{00} 
            &+& y_5\Phi_{01}  &+&  y_5^2\Delta_{22} & + &  x_5\Phi_{10}&+& x_5y_5\ssp \Delta_{211}\cr
\del_{y_5}\Delta_{32}&=&
            & & \ \  \Phi_{01}&+ & \hfill 2y_5 \Delta_{22} &   &            & + & \hfill x_5\Delta_{211}  \cr
\del_{y_5^2}\Delta_{32}&=&
            & &              & &  \hfill 2  \Delta_{22}     &   &               &   &             \cr
\del_{x_5}\Delta_{32}&=&
            & &              & &                 &  & \hfill \Phi_{10}     &  + &\hfill y_5  \Delta_{211} \cr
\del_{x_5}\del_{y_5}
            \Delta_{32}&=&
            & &              & &                 &   &            &   &   \hfill     \Delta_{211} \cr
}
\eqno 1.15
$$
Proceeding as we did for $\mu={22}$ we are to construct $5$ 
collections $\CB_{00},\CB_{01},\CB_{02},\CB_{10},\CB_{11}$
of polynomials in $x_1,x_2,x_3,x_4;y_1,y_2,y_3,y_4$ such that  
the polynomials in the union
$$
\eqalign{
\CB_{32} 
=\bigl\{b(\del)\Delta_{32}\bigr\}_{b\in \CB_{00}}+ 
\bigl\{ b(\del)&\del_{y_5} \Delta_{32}\bigr\}_{b\in \CB_{01}}+
\bigl\{ b(\del)\del_{y_5}^2\Delta_{32}\bigr\}_{b\in \CB_{02}}+\cr
&+\bigl\{ b(\del)\del_{x_5}\Delta_{32}\bigr\}_{b\in \CB_{10}} +
\bigl\{ b(\del)\del_{x_5}\del_{y_5}\Delta_{32}\bigr\}_{b\in \CB_{11}}\cr}
$$
are independent. If we succeed in choosing them so that $\CB_{32}$ has altogether $120$ elements
we will have proved the $n!$ conjecture for $\mu=32$. In  this case we shall not venture a
guess and determine the $\CB_{ij}$ from the equations that would have to hold if there was a
dependence between the elements of $\CB_{32}$. So let $b_{ij}$ be in the linear span of
$\CB_{ij}$ and  assume that we have
$$
b_{00}(\del)\Delta_{32}+  
b_{01}(\del)\del_{y_5}\Delta_{32}+  
b_{02}(\del)\del_{y_5}^2\Delta_{32}+  
b_{10}(\del)\del_{x_5}\Delta_{32}+  
b_{11}(\del)\del_{x_5}\del_{y_5}\Delta_{32}\ses 0
$$
Using the expansions in 1.15 and equating to zero the
coefficients of $1,y_5,y_5^2 ,x_5, x_5 y_5$ we get the
system of equations
$$
\matrix{
b_{00}(\del)\Phi_{00} 
&+  &
b_{01}(\del)\Phi_{01} 
&+  & 
2 b_{02}(\del)\Delta_{22} 
& + & 
b_{10}(\del)\Phi_{10} 
& + &
b_{11}(\del)\Delta_{31}&=0 
\cr 
& &
b_{00}(\del)\Phi_{01} 
& + & 2 b_{01}(\del)\Delta_{22}                      
&   &
& + & b_{10}(\del) \Delta_{31}                     
&=0 \cr
& &                     
& & b_{00}(\del)\Delta_{22}
&   &   
&  
& 
&=0  \cr
& &                     
& &                      
&   &
b_{00}(\del)\Phi_{10} 
&+ & 
b_{01}(\del)\Delta_{31}                      
&=0  \cr
& &                     
& &                      
&  
&                         
&  & b_{00}(\del)\Delta_{31}&=0  
\cr }
\eqno 1.16
$$
Note that in I.32 we have defined  $\CJ_\mu$ as the ideal of polynomials which kill $\Delta_\mu$.
With this notation the third and fifth equations
yield us that $b_{00}\in \CJ_{31}\wedge \CJ_{22}$ and
since  $(\CJ_{31}\wedge \CJ_{22})^\perp=\BM_{31}\vee \BM_{22}$, we see that
we can choose $\CB_{00}$ to be any bihomogenous basis of $\BM_{31}\vee \BM_{22}$. Indeed, if
we do so any linear combination $b_{00}$ of elements of $\CB_{00}$ which satisfies
the third and fifth equations in 1.16 would have to be orthogonal to itself and therefore must
have all its coefficients equal to zero. We shall visualize this choice by writing
$$
\CB_{00}\ses \potr \vee \ptt
\eqno 1.17
$$
This done,
setting $b_{00}=0$,
the fourth equation reduces to $b_{01}(\del)\Delta_{31}=0$ which gives $b_{01}\in
\CJ_{31}$. Thus using the same imagery we can set
$$
\CB_{01}\ses \potr
\eqno 1.18
$$
Substituting $b_{00}=b_{01}=0$  the second equation reduces to $b_{10}(\del)\Delta_{31}=0$. 
Thus we can again set
$$
\CB_{10}\ses \potr
\eqno 1.19
$$
Substituting $b_{00}=b_{01}=b_{10}=0$ the first equation reduces to 
$b_{02}(\del)\Delta_{22}+ b_{11}(\del)\Delta_{31}=0$. Here we have two possible choices. 
We can let $b_{11}$ vary freely in $\BM_{31}$ and force $b_{02}(\del)\Delta_{22}$
out of $\BM_{31}/\{0\}$. Alternatively we  can let $b_{02}$ vary freely 
in $\BM_{22}$ and force $b_{11}(\del)\Delta_{31}$ out of $\BM_{22}/\{0\}$.  
Let us take
$$
\CB_{11}\ses \potr
\eqno 1.20
$$
With this choice, we want $b_{02}(\del)\Delta_{22}\in\BM_{31}$ to imply $b_{02}=0$. We can
assure this if $b_{02}(\del)\Delta_{22}$ lies in the orthogonal complement of
$\BM_{31}\wedge \BM_{22}$ in $\BM_{22}$. In other words we want ${\bf flip}_{22}b_{02}$ 
\footnote{\dag}
{Here and after ${\bf flip}_\mu$ will denote flipping with respect to $\Delta_\mu$.}
to lie in 
$\bigl( \BM_{31}\wedge\BM_{22} \bigr)^\perp \wedge \BM_{22}$. 
To do this we must take $\CB_{02}$ to be a basis
of ${\bf flip}_{22}^{-1}\Bigl[\bigl(\BM_{31}\wedge \BM_{22}\bigr)^\perp\wedge \BM_{22}\Bigr]$. 
We represent this final choice by
writing
$$
\CB_{02}\ses {\bf flip}_{22}^{-1}\biggl[\Bigl(\potr \wedge \ptt\Bigr)^\perp \wedge \ptt\biggr]
$$
We have now assured that $\CB_{32}$ is an independent set. Let us count how many elements 
it has. We see that 1.18, 1.19 and 1.20 yield $3\times 24$ elements, moreover, 1.17 yields 
$2\times 24 - \dim  \BM_{31}\wedge\BM_{22}$ elements. Since $\dim \BM_{32}\leq 120$ the
independence of $\CB_{32}$ yields the inequality
$$
5\times 24 \sms \dim \poot \wedge \ptt\sps 
\dim \Bigl(\potr \wedge \ptt\Bigr)^\perp \wedge \ptt\ssp\ess \leq 120
$$   
or better
$$
\dim \Bigl(\potr \wedge \ptt\Bigr)^\perp \wedge \ptt 
\ssp\ess \leq\ess \ssp \dim \potr \wedge \ptt
\eqno 1.21
$$
In particular, $\CB_{32}$ would be a basis for $\BM_{32}$ and the $n!$ conjecture would then
be established for
$\mu=32$ if we prove that equality holds in 1.21. Amazingly it develops that this equality
holds true in the strongest possible sense. More precisely, we can prove that
$$
{\bf flip}_{22}\ssp
\bigr[ \BM_{31}\wedge \BM_{22}\bigr]\ses \bigl(\BM_{31}\wedge \BM_{22}\bigr)^\perp\wedge \BM_{22}
\ess .
\eqno 1.22
$$
Let us postpone the proof of this identity and proceed to examine all its implications. 
To begin with it shows that we can take $\CB_{02}$ to be any bihomogeneous basis of
$\BM_{31}\wedge\BM_{22}$. This done we can visualise the Hilbert series of $\BM_{32}$ by
writing
$$
G_{32}\ses  
\potr  \sps \ptt \sms  \potr \wedge \ptt\sps  q\ssp \potr  \sps t\ssp\potr  \sps tq \ssp \potr 
\sps q^2 \ssp \potr  \vee \ptt
\ess .
$$
To be precise, the left hand side of this relation is $G_{32}(1/q,1/t)t^2q^4$,
however we can write it that way because of 1.10. It will be convenient to extend the definition 
of ${\bf flip}_\mu$
to act on a Frobenius characteristic $H(x;q,t)$ and on a Hilbert series $F(q,t)$ by setting 
$$
{\bf flip}_\mu\ssp H(x;q,t)=\om H(x;1/q,1/t)\ssp t^{n(\mu)}q^{n(\mu')}
\ess ,\ess\ess\ess 
{\bf flip}_\mu\ssp F(q,t)=F(1/q,1/t)\ssp t^{n(\mu)}q^{n(\mu')}
\eqno 1.23
$$
This given, the more refined argument which takes account of the action of $S_4$ on the basis
$\CB_{32}$ yields the relation
$$
\del_{p_1}  C_{32}(x;q,t)= (1+t+q+tq) C_{31}(x;q,t)+  C_{22}(x;q,t)
+ (q^2-1)\ess \Phi
\eqno 1.24
$$
where for convenience we have denoted by $\Phi$ the Frobenius characteristic of $\potr
\wedge\ptt$. Here again the left-hand side should have been ${\bf flip}_{32}\del_{p_1} 
C_{32}$, but 1.9 makes it right the way it is. 
\sas
Note that if we follow the other alternative and choose $\CB_{02}=\ptt $, 
then we must take $\CB_{11}$ to be a basis of 
$$
{\bf flip}_{31}^{-1}\ssp \Bigl[ \Bigl(\BM_{31}\wedge \BM_{22}\Bigr)^\perp \wedge\BM_{31}  \Bigr]\ess .
$$
Since we can also show that
$$
{\bf flip}_{31}\ssp \Bigl(\BM_{31}\wedge \BM_{22}\Bigr)\ses 
 \Bigl(\BM_{31}\wedge \BM_{22}\Bigr)^\perp \wedge\BM_{31}
 \ess ,
\eqno 1.25
$$
we can again set
$$
\CB_{11}\ses \potr  \wedge \ptt\ess .
$$ 
These choices yield 
$$
\del_{p_1}  C_{32}(x;q,t)= (1+t+q) C_{31}(x;q,t)+ (1+q^2) C_{22}(x;q,t)
+ (tq-1)\ess \Phi\ess .
\eqno 1.26
$$
Subtracting 1.24 from 1.26 we deduce that these two expansions of
$\del_{p_1}  C_{32}(x;q,t)$ are one and the same if and only if
$$
\Phi\ses {t\ssp C_{31}- q\ssp C_{22}\over t-q}\ess .
\eqno 1.27
$$
It develops that this is a consequence of 1.22, and 1.25. Indeed   
these two relations yield that $\BM_{31}$ and $\BM_{22}$ decompose as follows
(as bigraded $S_4$-modules):
$$
\eqalign{
\BM_{31}&\ses\BM_{31}\wedge\BM_{22}\oplus {\bf flip}_{31}\BM_{31}\wedge\BM_{22}\ssp\ess ,\cr
\BM_{22}&\ses\BM_{31}\wedge\BM_{22}\oplus {\bf flip}_{22}\BM_{31}\wedge\BM_{22}\ssp\ess
.\cr}
$$
In particular we must also have that
$$
\eqalign{
C_{31}(x;q,t)&\ses\Phi+ \flip_{31}\ssp \Phi\ses\Phi+ tq^3 \downarrow\Phi\ess ,\cr
C_{22}(x;q,t)&\ses\Phi+ \flip_{22}\ssp \Phi\ses\Phi+ t^2q^2 \downarrow\Phi\ess .\cr}
\eqno 1.28
$$
Now it is easily seen that these relations are equivalent to 1.27.
The reader should recognize at this point that 1.22 and 1.25 are but two special instances of I.29.
In fact, in the notation of the introduction they may be written as
$$
{\bf flip}_1\ssp \BM_{32}^{11}= \BM_{32}^{10}
\ess\ess\ess\ess\ess {\rm and}\ess\ess\ess\ess
{\bf flip}_2\ssp \BM_{32}^{11}= \BM_{32}^{01} \ess .
$$

It goes without saying that we can verify by computer the validity of 1.22, 1.25 and 1.27.
We can also easily check that the symmetric
function given by 
$$
{t\ssp \TH_{31}- q\ssp \TH_{22}\over t-q}
$$
is indeed Schur positive and that the $C=\TH$ conjecture is
in fact true for $32$, $31$ and $22$. Nevertheless, it is more illuminating to
show how all of these identities may be proved using representation theory. 
\sas

Recalling that
$$
\Delta_{22}=\det  
\pmatrix{
1& 1& 1& 1\cr
y_1& y_2& y_3& y_4\cr
x_1& x_2& x_3& x_4\cr
x_1y_1& x_2y_2& x_3y_3& x_4y_4\cr}
, \bigsp
\Delta_{31}\ses \det 
\pmatrix{
1& 1& 1& 1\cr
y_1& y_2& y_3& y_4\cr
y_1^2& y_2^2& y_3^2& y_4^2\cr
x_1& x_2& x_3& x_4\cr
}\ess ,
$$
we get
$$
-\del_{x_4}\del_{y_4}\Delta_{22}\ses {1\over 2}\ssp \del_{y_4}^2 \Delta_{31}\ses \det 
\pmatrix{
1& 1& 1\cr
x_1& x_2& x_3\cr
y_1& y_2& y_3\cr
}
$$
Now from a well known result of A. Young (see Theorem 5.8 of [\youngsubst]) we derive that
the action of $S_4$ on this element
generates a $3$-dimensional irreducible representation with character 
$\chi^{211}$ and weight $tq$ which is shared by $\BM_{31}$ and $\BM_{22}$. Similarly, we see
that  its partial derivatives 
$$
\del_{x_3} 
\det  \pmatrix{
1& 1& 1\cr
x_1& x_2& x_3\cr
y_1& y_2& y_3\cr}= 
-\det  \pmatrix{
1& 1\cr
y_1& y_2\cr}\ssp ,
\ess\ess\ess
\del_{y_3} \det  \pmatrix{
1& 1& 1\cr
x_1& x_2& x_3\cr
y_1& y_2& y_3\cr}=
\det  \pmatrix{
1& 1\cr
x_1& x_2\cr}\ess ,
$$
generate common $3$-dimensional representations with character $\chi^{31}$ 
and weights $q$ and $t$ respectively.
Now, we also have
$$
{1\over 2}(\del_{x_4}\del_{y_3}-\del_{x_3}\del_{y_4})\Delta_{31}\ses 2 \del_{x_3}\del_{x_4}\Delta_{22}\ses 
\det  \pmatrix{
1& 1\cr
y_1& y_2\cr}
\det  \pmatrix{
1& 1\cr
y_3& y_4\cr}\ess,
$$
and this generates a common $2$ dimensional representation with character $\chi^{22}$ and weight $q^2$.
Including the trivial, we have thus identified a submodule of $\BM_{31}\wedge\BM_{22}$
with bigraded Frobenius characteristic
$$
\Phi\ses S_{4}+(t+q)\ssp S_{31}\sps tq\ssp S_{211}\sps q^2\ssp S_{22}\ess .
$$
Using I.19 and I.20 and the definition in 1.23 we then immediately derive that $\BM_{31}$ and
$\BM_{22}$ contain submodules with respective Frobenius characteristics
$$
\eqalign{
\flip_{31}\ssp\Phi&\ses tq^3\ssp S_{1^4}+(q^3+tq^2)\ssp S_{211}\sps q^2\ssp S_{31}\sps tq\ssp S_{22}\ess ,\cr
\bigsp
\flip_{22}\ssp \Phi&\ses t^2q^2S_{1^4}+(tq^2+t^2q)\ssp S_{211}\sps tq\ssp S_{31}\sps t^2\ssp S_{22}\ess .\cr }
\eqno 1.29
$$
Since, $\Phi$ and $\flip_{31}\ssp\Phi$ have no common terms and each accounts for $12$ dimensions
we must conclude that together they must give the Frobenius characteristic of $\BM_{31}$. 
A similar reasoning applies to $\Phi$ and $\flip_{22}\ssp\Phi$ and we can write
$$
C_{31}\ses \Phi\sps \flip_{31}\ssp\Phi\ess ,
\bigsp
C_{22}\ses \Phi\sps \flip_{22}\ssp\Phi\ess ,
\eqno 1.30
$$
Comparing the expressions in 1.29, we see that the possibility still remains that $\BM_{31}$ and $\BM_{22}$
may have in common an irreducible submodule with character $\chi^{211}$ and weight $tq^2$. 
However, this submodule is generated in $\BM_{22}$ by the action of $S_4$ on the polynomial
$$
P\ses \det\pmatrix{1 &1\cr \del_{x_1}&\del_{x_2}\cr}\ssp \Delta_{22}
$$
If this polynomial were to belong to $\BM_{31}$ it would be killed by any element that kills $\Delta_{31}$. 
In particular
we should have $\del_{x_4} \del_{y_4}P=0$. But we see that
$$
\del_{x_4} \del_{y_4}P\ses \det\pmatrix{1 &1\cr \del_{x_1}&\del_{x_2}\cr}\ssp
 \det \pmatrix{1&1&1\cr x_1& x_2& x_3\cr y_1& y_2& y_3\cr}\ses y_1+y_2-2y_3\neq 0\ess .
$$
This completes the proof of 1.22 and 1.25 and shows that $\Phi$ is none other than the
Frobenius characteristic of the intersection of $\BM_{31}$ and $\BM_{22}$.\
\sa

\Rem{1.1}{It is interesting to see what 1.24 reduces to when we express its right hand side
entirely in terms of $\Phi$. To this end we use the identities in 1.29 to replace
$C_{31}$ and $C_{22}$ and obtain (using the notation introduced in I.21)
$$
\eqalign{
\del_{p_1}\ssp C_{32} &\ses (1+t+q+tq)(\Phi+{\bf flip}_{31}\Phi)\sps (\Phi+{\bf flip}_{22}\Phi)
 \sps (q^2-1)\Phi\cr
&\ses
(1+t+q+tq+q^2)\Phi\sps (1+t+q+tq)tq^3\da\Phi
\sps t^2q^2\da\Phi\cr
&\ses
(1+t+q+tq+q^2)\Phi\sps t^2q^4(1+{1\over t}+{1 \over q}+{1\over tq}+{1\over q^2})\da \Phi
\ess .\cr}
$$
and this may be rewritten in the very suggestive form
$$
\del_{p_1}\ssp C_{32}\ses B_{32}(q,t)\ssp \Phi \sps {\bf flip}_{32}\bigl(B_{32}(q,t)\ssp \Phi\bigr) 
\ess .
\eqno 1.31
$$}

What we have discovered in this particular example holds true in full generality.
\sas

\Prop{1.1}Let $\mu$ be a two-corner partition of $n+1$ and let $\aaa^{(1)}$ and $\aaa^{(2)}$ be the partitions obtained 
by removing one of the corners. Then 
$$
\eqalign{
\dim \Bigl((\BM_{\aaa^{(1)}}\wedge\BM_{\aaa^{(2)}})^\perp \wedge \BM_{\aaa^{(1)}}\Bigr) \ess 
&\leq\ess
\dim  (\BM_{\aaa^{(1)}}\wedge\BM_{\aaa^{(2)}} )
\ess ,\cr
\dim \Bigl( (\BM_{\aaa^{(1)}}\wedge\BM_{\aaa^{(2)}} )^\perp \wedge \BM_{\aaa^{(2)}}\Bigr) 
\ess &\leq\ess
\dim  (\BM_{\aaa^{(1)}}\wedge\BM_{\aaa^{(2)}} )
\ess .\cr
}
\eqno 1.32 
$$
If the $n!$ conjecture holds for $\aaa^{(1)}$ and $\aaa^{(2)}$ and
$$
\dim  (\BM_{\aaa^{(1)}}\wedge\BM_{\aaa^{(2)}} )\ses {n!\over 2}
\eqno 1.33
$$
then equalities hold in 1.32 and the $(n+1)!$-conjecture holds for $\mu$ as well.
Moreover, if
$$
\eqalign{
  (\BM_{\aaa^{(1)}}\wedge\BM_{\aaa^{(2)}})^\perp \wedge \BM_{\aaa^{(1)}}  \ses 
& 
{\bf flip}_{\aaa^{(1)}}  (\BM_{\aaa^{(1)}}\wedge\BM_{\aaa^{(2)}} )
\ess ,\cr
 (\BM_{\aaa^{(1)}}\wedge\BM_{\aaa^{(2)}} )^\perp \wedge \BM_{\aaa^{(2)}} 
\ses 
& 
{\bf flip}_{\aaa^{(2)}}  (\BM_{\aaa^{(1)}}\wedge\BM_{\aaa^{(2)}} )
\ess .\cr
}
\eqno 1.34 
$$
Then, if we let $\Phi$ denote the bivariate Frobenius characteristic of $\BM_{\aaa^{(1)}}\wedge\BM_{\aaa^{(2)}}$,
we also have the identity
$$
\del_{p_1} C_\mu(x;q,t)\ses B_\mu(q,t)\ssp \Phi \sps {\bf flip_\mu}\bigl( B_\mu(q,t)\ssp \Phi\bigr)
\eqno 1.35 
$$

\noindent The proof of this will be given in section 4.
\sas

This given, we see that the validity of Conjecture I.2 (restated for the case of $\la,\mu$ predecessors
of a two-corner partition) depends only on the verification of the equality in 1.33

\section{2. The three corner case}
We will better understand the formalism we have used to state properties (i)--(iv) of the heuristic,
if we study in detail a sufficiently rich special case. To this end, we  consider
the case of a partition $\mu$ with 3 corners. Letting $\aaa^{(1)},\aaa^{(2)},\aaa^{(3)}$
denote the predecessors of $\mu$, we see that property (i) states that the space 
$$
{\bf V}_\mu=\BM_{\aaa^{(1)}}\vee \BM_{\aaa^{(2)}}\vee \BM_{\alpha^{(3)}}
$$
has a basis ${\cal B}$, with 3 subsets ${\cal B}_1$, ${\cal B}_2$ and ${\cal B}_3$ such that
$$
\BM_{\aaa^{(1)}}={\cal L}[{\cal B}_1],\qquad 
\BM_{\aaa^{(2)}}={\cal L}[{\cal B}_2],\qquad
\BM_{\alpha^{(3)}}={\cal L}[{\cal B}_3]\ess .
$$
Moreover, each of these subsets breaks up as disjoint unions
$$\eqalign{ {\cal B}_1&= {\cal B}^{100}+{\cal B}^{110}+{\cal B}^{101}+{\cal B}^{111},\cr
              {\cal B}_2&= {\cal B}^{010}+{\cal B}^{110}+{\cal B}^{011}+{\cal B}^{111},\cr
              {\cal B}_3&= {\cal B}^{001}+{\cal B}^{011}+{\cal B}^{011}+{\cal B}^{111},\cr
}
$$
where, for example,
$$
\eqalign{ 
{\cal B}^{100}= {\cal B}_1\cap \overline{{\cal B}}_2 \cap \overline{{\cal B}}_3,\cr
{\cal B}^{110}= {\cal B}_1\cap {\cal B}_2 \cap \overline{{\cal B}}_3,\cr
{\cal B}^{101}= {\cal B}_1\cap \overline{{\cal B}}_2 \cap {\cal B}_3,\cr
{\cal B}^{111}= {\cal B}_1\cap {\cal B}_2 \cap {\cal B}_3,\cr
}
$$
with $\overline{{\cal B}}_i={\cal B}\setminus {\cal B}_i$.

%All of this information is best visualized by the following Venn diagram.
%\Lafig{3Bcircles}{93mm}{88mm}{500}{1}{Basis decomposition for 3 corners}

%\noindent
%By abuse of notation, we use the same imagery to represent the various subspaces $\BM^\eee_\mu$. This given, we
%may visualize parts (ii) and (iii) of the heuristic by the diagram below.

%\Lafig{3Mcircles}{93mm}{88mm}{600}{2}{Module decomposition and flips for 3 corners}

%\noindent
%Here, each of the regions corresponds to an $S_n$ invariant subspace, and the arrows represent the effect of
%flip maps. We should add that arrows whithin the circle labelled $\BM_{\aaa^{(i)}}$, represent
%the effect of ${\bf flip}_i$. We can thus visualize the geometric significance of (I.29), which
In this case, (I.29) asserts that
  $$\matrix{{\rm a)}\ \ {\bf flip}_1 \BM^{100}\cong \BM^{111},\qquad\quad &
              \hbox{\rm a')}\ \ {\bf flip}_1 \BM^{111}\cong \BM^{100},\cr
            {\rm b)}\ \ {\bf flip}_1 \BM^{110}\cong \BM^{101},\hfill &
              \hbox{\rm b')}\ \ {\bf flip}_1 \BM^{101}\cong \BM^{110},\cr
            {\rm c)}\ \ {\bf flip}_2 \BM^{110}\cong \BM^{011},\hfill &
              \hbox{\rm c')}\ \ {\bf flip}_2 \BM^{011}\cong \BM^{110},\cr
            {\rm d)}\ \ {\bf flip}_2 \BM^{010}\cong \BM^{111},\hfill &
              \hbox{\rm d')}\ \ {\bf flip}_2 \BM^{111}\cong \BM^{010},\cr
            {\rm e)}\ \ {\bf flip}_3 \BM^{101}\cong \BM^{011},\hfill &
              \hbox{\rm e')}\ \ {\bf flip}_3 \BM^{011}\cong \BM^{101},\cr
            {\rm f)}\ \ {\bf flip}_3 \BM^{001}\cong \BM^{111},\hfill &
              \hbox{\rm f')}\ \ {\bf flip}_3 \BM^{111}\cong \BM^{001}.\cr}\eqno{2.1}$$
Of course these relations yield identities involving the corresponding bivariate characteristics as indicated in
(I.28). For example, b'), c) and e') yield the following identities 
  $$\eqalign{{\rm 1)}\ \ T_1\da \Phi^{101}=\Phi^{110},\cr
             {\rm 2)}\ \ T_2\da \Phi^{110}=\Phi^{011},\cr
             {\rm 3)}\ \ T_3\da \Phi^{011}=\Phi^{101}.\cr}\eqno{2.2}$$
Note that combining 1) and 2) above, and using the involutory nature of $\da$, we derive that
  $$\Phi^{011}=T_2\da T_1\da \Phi^{101}={T_2\over T_1}\, \Phi^{101}.$$
Similarly, 2) and 3) of (2.2) yield
$$
\Phi^{101}=T_3\da T_2\da \Phi^{110}={T_3\over T_2}\, \Phi^{110},
$$
and these two relations give that
$$
T_1\,\Phi^{011} = T_2\,\Phi^{101}= T_3\,\Phi^{110}.
\eqno{2.3}
$$
Thus, if we denote $\Phi^{(2]}$ this common expression, we can write
$$
\Phi^{011}={\Phi^{(2)}\over T_1},\qquad
\Phi^{101}={\Phi^{(2)}\over T_2},\qquad  
\Phi^{110}={\Phi^{(2)}\over T_3}.
\eqno{2.4}
$$  
In the same manner, using a'), d') and f') of (2.1), we get
  $$T_1\da \Phi^{111}=\Phi^{100},\qquad
    T_2\da \Phi^{111}=\Phi^{010},\qquad    
    T_3\da \Phi^{111}=\Phi^{001}.$$ 
Thus, if we set 
   $$\Phi^{(1)}:=T_1T_2T_3\, \da \Phi^{111},$$
we obtain that
   $$\Phi^{100}={\Phi^{(1)}\over T_2T_3},\qquad
    \Phi^{010}={\Phi^{(1)}\over T_1T_3},\qquad  
    \Phi^{001}={\Phi^{(1)}\over T_1T_2}.\eqno{2.5}$$  
This leads us to discover the remarkable fact (mentioned in the introduction)
that each $\Phi^\eee_\mu$, up to a factor, depends
only on the sum of the $\eee_i$'s. Postponing the proof of the general result to the next section, it will
be good to see what can be further derived from these relations.

First of all, we should note that part (ii)  of the SF heuristic implies that the
$\BM_{\aaa^{(i)}}$'s have the following decomposition into $S_n$-invariant submodules\footnote{$^{\rm \dag}$}{(see
Figure 2)}.
$$
\eqalign{    \BM_{\aaa^{(1)}}=\BM^{100}\oplus \BM^{101}\oplus \BM^{110}\oplus \BM^{111}\ess ,\cr
             \BM_{\aaa^{(2)}}=\BM^{010}\oplus \BM^{110}\oplus \BM^{011}\oplus \BM^{111}\ess ,\cr
             \BM_{\aaa^{(3)}}=\BM^{001}\oplus \BM^{101}\oplus \BM^{011}\oplus \BM^{111}\ess .\cr}
\eqno{2.6}
$$
It is clear that in general we have the following decomposition
   $$\BM_{\aaa^{(i)}}=\bigoplus_{\eee\atop {\eee_i=1}} \BM_\mu^{\eee}.$$
Of course, all these relations yield corresponding relations for the associated bivariate Frobenius
characteristics. For instance, those in (2.6) combined with (2.4) and (2.5) give
  $$\eqalign{{\cal F}(\BM_{\aaa^{(1)}})=
               {1\over T_2T_3}\Phi^{(1)}\,+\, {1\over T_2}\Phi^{(2)}\,+\,{1\over T_3}\Phi^{(2)}\,+\,
                        \Phi^{(3)}\ess ,\cr
             {\cal F}(\BM_{\aaa^{(2)}})=
               {1\over T_1T_3}\Phi^{(1)}\,+\,{1\over T_3}\Phi^{(2)}\,+\,{1\over T_1}\Phi^{(2)}\,+\,
                        \Phi^{(3)}\ess ,\cr
             {\cal F}(\BM_{\aaa^{(3)}})=
               {1\over T_1T_2}\Phi^{(1)}\,+\,{1\over T_2}\Phi^{(2)}\,+\,{1\over T_1}\Phi^{(2)}\,+\,
                       \Phi^{(3)}\ess ,\cr}
$$
where for consistency we have set
$$
\Phi^{111}:=\Phi^{(3)}\ess .
$$
Note that the $C=\TH$ conjecture then yields the expansions
$$
\eqalign{
\TH_{\aaa^{(1)}}=
               {1\over T_2T_3}\Phi^{(1)}\,+\, {1\over T_2}\Phi^{(2)}\,+\,{1\over T_3}\Phi^{(2)}\,+ 
                        \Phi^{(3)}\ess ,\cr
\TH_{\aaa^{(2)}}=
               {1\over T_1T_3}\Phi^{(1)}\,+\,{1\over T_3}\Phi^{(2)}\,+\,{1\over T_1}\Phi^{(2)}\,+ 
                         \Phi^{(3)}\ess  ,\cr
\TH_{\aaa^{(3)}}=
               {1\over T_1T_2}\Phi^{(1)}\,+\,{1\over T_2}\Phi^{(2)}\,+\,{1\over T_1}\Phi^{(2)}\,+  
                       \Phi^{(3)}\ess .\cr}
\eqno{2.7}
$$
which is I.38 for $m=3$ 
\sas

Note further that the non singularity of the flip maps combined with the relations in (2.1) implies that
$$
\dim \BM^{100}=\dim \BM^{010}=\dim \BM^{001} =\dim \BM^{111}\ess ,
\qquad\qquad
\dim \BM^{110}=\dim \BM^{101}=\dim \BM^{011}\ess . 
$$   
Calling $d_1$ the first common dimension and $d_2$ the second,
any  of the equations in (2.6) gives that
  $$n!=2\,d_1+2\,d_2.$$
If we add to this, part (v) of the heuristic which gives that $\dim  \BM^{111}=n!/3$, we get
  $$d_1={n!\over 3},\qquad {\rm and}\qquad d_2={n!\over 6}.
\eqno{2.8}
$$
Since $\BM_{\aaa^{(1)}}\cap \BM_{\aaa^{(2)}}=\BM^{110}+\BM^{111}$, we deduce from (2.8) that 
    $$\dim  \BM_{\aaa^{(1)}}\cap \BM_{\aaa^{(2)}}={n!\over 2},$$
which shows that at least in this case part (v) of the heuristic is self-consistent.
\sas

Leaving all this aside for a moment, we shall next narrow our study to the case $\mu=321$.
We begin by showing here, that the heuristic enables us to push through 
the recursive approach to proving  the $n!$ conjecture. 
We shall then complete the proof by the explicit construction of a basis $\CB$
which satisfies properties (i)-(iv) required by the heuristic. In particular,
since here
$$
\aaa^{(1)}\ses \deux32
\ess\ess\ess ,\ess\ess\ess
\aaa^{(2)}\ses \trois311
\ess\ess\ess ,\ess\ess\ess
\aaa^{(3)}\ses \trois221\ess ,
$$
we see that $\CB_1,\CB_2,\CB_3$ must respectively give bases for 
$$
\BM_{32}=\CL[\del_x^p\del_y^q\Delta_\deux32]
 \ess ,\ess\ess\ess
\BM_{311}=\CL[\del_x^p\del_y^q\Delta_\trois311] 
 \ess ,\ess\ess\ess
\BM_{221}=\CL[\del_x^p\del_y^q\Delta_\trois221]\ess .
$$
Let us also keep in mind that according to I.33 we must have
$$
\CB\backslash \CB_1\con\BM_{32}^\perp
\ess ,\ess\ess\ess
\CB\backslash \CB_2\con\BM_{311}^\perp
\ess ,\ess\ess\ess
\CB\backslash \CB_3\con\BM_{221}^\perp\ess .
\eqno 2.9
$$
\sa
  
Note now that in this case, $\Delta_\mu$ is given by the
determinant in (I.11) which may be written as
$$
\Delta_{\trois321}=\PD_{00}+y\,\PD_{01}+y^2\PD_{02}+x\,\PD_{10}+xy\,\PD_{11}+x^2\PD_{20},
\eqno{2.10}
$$ 
where for simplicity we have set $x_6=x$ and $y_6=y$, and $\PD_{ij}$ represents the appropriately signed
complementary minor of the element $x^iy^j$. 
Using 1.1 again, our strategy will be to construct a collection ${\cal B}_{321}\subset \BM_{321}$ of the form
  $$
  \eqalign{
{\cal B}_{321}  =\bigl\{b(\del)\,&\Delta_{321}\bigr\}_{b\in {\cal B}_{00}}+ 
                 \bigl\{ b(\del)\,\del_{y}\,\Delta_{321}\bigr\}_{b\in {\cal B}_{01}}+
                 \bigl\{ b(\del)\,\del_{y}^2\,\Delta_{321}\bigr\}_{b\in {\cal B}_{02}}+\cr &
                 \bigl\{b(\del)\,\del_{x}\,\Delta_{321}\bigr\}_{b\in {\cal B}_{10}}+
                 \bigl\{b(\del)\,\del_{x}\del_{y}\,\Delta_{321}\bigr\}_{b\in {\cal B}_{11}}+
                 \bigl\{ b(\del)\,\del_{x}^2\,\Delta_{321}\bigr\}_{b\in {\cal B}_{20}}\cr}
\eqno 2.11
 $$
where the ${\cal B}_{ij}$ are collections of polynomials in the variables $x_1,x_2,\ldots, x_5; y_1,y_2,\ldots,y_5$, to
be determined from the equations that would have to assure that
 
  \itemitem{1)}  ${\cal B}_{321}$ is an independent set,
  \itemitem{2)} the cardinality of ${\cal B}_{321}$ is $6!$,

\noindent To construct these equations, we assume that a basis $\CB$ with properties (i)-(v)
has been constructed, then determine how the subsets ${\cal B}_{ij}\con \CB$
must be chosen so that no matter how we pick $b_{ij}\in{\cal L}[{\cal B}_{ij}]$ the equality
$$
\eqalign{
  b_{00}(\del)\,\Delta_{321}+  
  b_{01}(\del)\, \del_{y}\,\Delta_{321}&+  
  b_{02}(\del)\,\del_{y}^2\Delta_{321}+\cr  
 &+ b_{10}(\del)\,\del_{x}\,\Delta_{321}+  
  b_{11}(\del)\,\del_{x}\del_{y}\,\Delta_{321}+  
  b_{20}(\del)\,\del_{x}^2\Delta_{321}\,=\, 0\cr
}
\eqno{2.12}
$$
forces all the $b_{ij}$ to vanish identically.

Note that from  2.10  we derive 
$$
\matrix{
\Delta_{321}\hf &=& \PD_{00} &+& \hf y\,\PD_{01}&-& \hf y^2\Delta_{\trois221}  &+&  x\,\PD_{10} & - &  \hf xy\,\Delta_{\trois311} & + & 
              x^2\Delta_{\deux32}\cr
\del_{y}
 \Delta_{321}\hf &=&          & & \hf  \PD_{01}&-& \hf 2 y\,\Delta_{\trois221}  & &         & - & \hf x\,\Delta_{\trois311}&&             
\cr
\del_{y}^2
 \Delta_{321}\hf &=&          & &            &-&\hf 2 \,\Delta_{\trois221}  & &            &   &             &   &              \cr
\del_{x}
 \Delta_{321}\hf &=&          & &            & &              & &\hf  \PD_{10} & - &\hf  y\,\Delta_{\trois311} & + &  2
            x\,\Delta_{\deux32}\cr
\del_{x}\del_{y}
 \Delta_{321}\hf &=&          & &            & &              & &            & - &\hf   \Delta_{\trois311} &   &            \cr
\del_{x}^2
 \Delta_{321}\hf &=&          & &            & &              & &            &   &             &   &\hf 2 \,\Delta_{\deux32}\cr
}
\eqno{2.13}
$$
where we have used the fact that $\PD_{02}=-\Delta_{32}$, $\PD_{11}=-\Delta_{311}$ and $\PD_{20}=\Delta_{221}$. Thus
equating to $0$ the coefficients of $x^iy^j$ in (2.12), we can replace it by the following system of 6 equations
$$
\matrix{
         {\rm a)}\  b_{00}(\del)\PD_{00}    &\hf +    b_{01}(\del)\PD_{01}    &\hf -2   b_{02}(\del)\Delta_{\trois221} &
         \hf+  b_{10}(\del)\PD_{10}    &\hf -    b_{11}(\del)\Delta_{\trois311}& + \hf 2 b_{20}(\del)\Delta_{\deux32}&     =0 \cr
         {\rm b)}\hf                    &\hf     b_{00}(\del)\PD_{01}    &\hf-   2b_{01}(\del)\Delta_{\trois221} &
                                    &\hf-    b_{10}(\del)\Delta_{\trois311}&                              &     =0 \cr
         {\rm c)}\hf                   &                               &\hf     b_{00}(\del)\Delta_{\trois221} &
                                    &                              &                              &     =0 \cr
         {\rm d)}\hf                    &                              &                              &
         \hf b_{00}(\del)\PD_{10}    &\hf-    b_{01}(\del)\Delta_{\trois311}&\hf +   2b_{10}(\del)\Delta_{\deux32}&     =0 \cr
         {\rm e)}\hf                                                 &                              &
            &                        &\hf   b_{00}(\del)\Delta_{\trois311}&        &                           =0 \cr
          {\rm f)}\hf                  &                              &                              &
                                    &                              &\hf    b_{00}(\del)\Delta_{\deux32}&     =0 \cr
}
\eqno{2.14}
$$
Our guiding principle in looking for such  sets ${\cal B}_{i,j}$, will be to choose them as ``large'' as possible in a manner
that still forces equations in  2.14  to have as unique solution $b_{i,j}=0$ for all $i,j$. 
To this end, note that whatever our choice of $\CB_{00}$ is in $\BM_{\aaa^{(1)}}\vee\BM_{\aaa^{(2)}}\vee
\BM_{\aaa^{(3)}}$   equations c), e) and f) force $b_{00}$ to be orthogonal to itself and therefore
identically zero. Hence our maximal choice is
$$
{\cal B}_{00}={\cal B}_1\cup {\cal B}_2\cup {\cal B}_3\ess .
$$ 
To proceed, we need to adopt a convention. For a polynomial
$P\in \CL[\CB ]$ we shall let $P^{\eee_1\eee_2\eee_3}$ denote the projection of
$P$ into $\BM^{\eee_1\eee_2\eee_3}$. Putting it in another way, if
$$
P\ses \sum_{b\in \CB }\ssp c_b \ssp b
$$
then
$$
P^{\eee_1\eee_2\eee_3}\ses \sum_{b\in \CB^{\eee_1\eee_2\eee_3} }\ssp c_b \ssp b
$$
Thus, whatever $\CB_{10}$ and  $\CB_{01}$ are chosen to be in $\CB_1\cup\CB_2\cup\CB_3$, 
our elements $b_{10}$ and $b_{10}$ will have the decompositions
$$
\eqalign{
& b_{10}
\ses
b_{10}^{100}\sps
b_{10}^{110}\sps
b_{10}^{010}\sps 
b_{10}^{111}\sps
b_{10}^{101}\sps
b_{10}^{011}\sps
b_{10}^{001}\ess , \cr
& b_{01}
\ses
b_{01}^{100}\sps
b_{01}^{110}\sps
b_{01}^{010}\sps 
b_{01}^{111}\sps
b_{01}^{101}\sps
b_{01}^{011}\sps
b_{01}^{001}\ess , \cr
}\eqno 2.15
$$
Note that since the spaces $M^{\eee_1\eee_2\eee_3}$ are orthogonal to each other  
$b_{01}$ or $b_{10}$ will vanish if and only if each of their components
$b_{01}^{\eee_1\eee_2\eee_3}$ and $b_{10}^{\eee_1\eee_2\eee_3}$ 
separately  vanishes. Now, without any assumptions regarding the nature
of the determinants $D_{ij}$ occurring in 2.14, our only method for
showing the vanishing of one of these components is to force some of its
flips to vanish.  Under these conditions, the only way we can show $b_{10}^{001}$
to vanish is to force  
$$
b_{10}^{001}(\del)\Delta_{\trois221}=0\ess .
$$
However, we see no way to extract such a result out of 2.14 since $b_{10}$ 
only acts on $\Delta_{\trois311}$ and $\Delta_{\deux32}$. In conclusion, we would
be unable to force the vanishing of $b_{10}$ if we allow $b_{10}^{001}\neq 0$.
This given, the most we can include in  $\CB_{10}$ is
all of $\CB_1\cup\CB_2$. An analogous reasoning yields that the most we can include in
$\CB_{01}$ is all of $\CB_2\cup\CB_3$. This leaves us with
$$
\eqalign{
& b_{10}
\ses
b_{10}^{100}\sps
b_{10}^{110}\sps
b_{10}^{010}\sps 
b_{10}^{111}\sps
b_{10}^{101}\sps
b_{10}^{011}\ess , \cr
& b_{01}
\ses
b_{01}^{110}\sps
b_{01}^{010}\sps
b_{01}^{111}\sps 
b_{01}^{101}\sps
b_{01}^{011}\sps
b_{01}^{001}\ess . \cr
}\eqno 2.16
$$

Note next that setting $b_{00}=0$ in 2.14 reduces the $2^{nd}$ and $4^{th}$ equations to
$$
\eqalign{
 &{\rm b)}\ess\ess\ess             2b_{01}(\del)\ssp \Delta_{\trois221} \sps  \ssp b_{10}(\del)\ssp\Delta_{\trois311} \ssp\ses 0 \ess ,\cr
 &{\rm d)}\ess              - b_{01}(\del)\ssp \Delta_{\trois311}  \sps    2b_{10}(\del)\ssp \Delta_{\deux32} \ses 0 \ess .\cr
}
\eqno 2.17
$$
In other words we must have 
$$
\eqalign{
 &{\rm b)}\ess\ess\ess 2\ssp{\bf flip \  }_{\trois221}\ssp b_{01}\ses \ssp  - {\bf flip \  }_{\trois311}\ssp b_{10} \ess ,\cr
 &{\rm d)}\ess\ssp\ssp\ess\ess{\bf flip \  }_{\trois311}\ssp b_{01}\ses 2\ssp {\bf flip \  }_{\deux32}\ssp   b_{10}\ess .\cr
}
$$
Projecting these equations in each of our spaces $\BM^{\eee_1\eee_2\eee_3}$
and using the relations in 2.1 as equalities we obtain seven pairs of 
equations whose consequences are easily
visualised from the two tables given below.
$$
\matrix{
{\rm by}\ess{\bf flip\ssp}_{\trois221}  &    &{\rm  Spaces} &     &{\rm by}\ess {\bf flip\ssp}_{\trois311} \cr               
                                  &    &              &     &                                  \cr               
             \{0\}                   &    &   \BM^{100}  &     &              \{0\}                  \cr 
                                  &    &              &     &                                  \cr               
             \{0\}                   &    &   \BM^{110}  &\lar &      b_{10}^{011}                \cr 
                                  &    &              &     &                                  \cr               
             \{0\}                   &    &   \BM^{010}  &\lar &      b_{10}^{111}                \cr 
                                  &    &              &     &                                  \cr               
    b_{01}^{011}                  &\rar&   \BM^{101}  &     &              \{0\}                  \cr
                                  &    &              &     &                                  \cr               
    b_{01}^{001}                  &\rar&   \BM^{111}  &\lar &      b_{10}^{010}                \cr 
                                  &    &              &     &                                  \cr               
    b_{01}^{101}                  &\rar&   \BM^{011}  &\lar &      b_{10}^{110}                \cr 
                                  &    &              &     &                                  \cr               
    b_{01}^{111}                  &\rar&   \BM^{001}  &     &              \{0\}                  \cr 
}
\ess\ess\ess ,\ess\ess\ess
\matrix{
{\rm by}\ess{\bf flip\ssp}_{\trois311}  &    &{\rm  Spaces} &     &{\rm by}\ess {\bf flip\ssp}_{\deux32} \cr               
                                  &    &              &     &                                  \cr               
             \{0\}                   &    &   \BM^{100}  &\lar &      b_{10}^{111}                \cr 
                                  &    &              &     &                                  \cr               
    b_{01}^{011}                  &\rar&   \BM^{110}  &\lar &      b_{10}^{101}                \cr 
                                  &    &              &     &                                  \cr               
    b_{01}^{111}                  &\rar&   \BM^{010}  &     &              \{0\}                  \cr
                                  &    &              &     &                                  \cr               
             \{0\}                   &    &   \BM^{101}  &\lar &      b_{10}^{110}                \cr 
                                  &    &              &     &                                  \cr               
    b_{01}^{010}                  &\rar&   \BM^{111}  &\lar &      b_{10}^{100}                \cr 
                                  &    &              &     &                                  \cr               
    b_{01}^{110}                  &\rar&   \BM^{011}  &     &              \{0\}                  \cr 
                                  &    &              &     &                                  \cr               
             \{0\}                   &    &   \BM^{001}  &     &              \{0\}                  \cr 
}
$$
For instance, the second row of the table in the left states that ${\bf flip \  }_{\trois311}$
sends $b_{10}^{011}$ into $\BM^{110}$ and at the same time there is no component  of $b_{01}$ 
that is sent into  $\BM^{110}$ by ${\bf flip \  }_{\trois221}\ssp $. Thus, as long as 
$b_{10}$ and $b_{01}$ are are given by 2.16, equation b) forces $b_{10}^{011}=0$.
On the other hand, the second row of the table on the right, simply says that ${\bf flip \  }_{\trois311}$
sends $b_{01}^{011}$ into $\BM^{110}$ and at the same time ${\bf flip \  }_{\deux32}$
sends $b_{10}^{101}$ into the same space. So we can't deduce any implications at this stage from
this particullar row. Proceeding in this manner we immediately derive from equation b)
and the table on the left that $b_{10}^{011}=b_{10}^{111}=b_{01}^{011}=b_{01}^{111}=0$.
Similarly, equation d) and the table on the right imply that
$b_{10}^{111}=b_{10}^{110}=b_{01}^{111}=b_{01}^{110}=0$. Less trivially, if we feed the 
information we get from one table into the other, we see that the vanishing of $b_{10}^{110}$
forced by the fourth row of the table on the right, makes no component of 
$\ssp{\bf flip\ssp }_{\trois311}\ssp b_{10}$
available to match the image of $b_{01}^{101}$ under ${\bf flip\ }_{\trois221}\ssp $. So from
the sixth row of the table on the left we deduce that $b_{01}^{101}=0$. A similar reasoning
based on the second row of the table on the right yields that $b_{10}^{101}=0$.
This given, without committing ourselves to any particular choices of $\CB_{01}$
and $\CB_{10}$, (other than 2.15), our equations reduce the decompositions in 2.15 to 
$$
\eqalign{
& b_{10}
\ses
b_{10}^{100}\sps
b_{10}^{010}\ess , \cr
& b_{01}
\ses
b_{01}^{010}\sps
b_{01}^{001}\ess , \cr
}
$$
Setting all the remaining components equal to zero reduces both our tables to a single row:
$$
\matrix{
{\rm by}\ess{\bf flip\ssp}_{\trois221}  &    &{\rm  Spaces} &     &{\rm by}\ess {\bf flip\ssp}_{\trois311} \cr               
                                  &    &              &     &                                  \cr               
    b_{01}^{001}                  &\rar&   \BM^{111}  &\lar &      b_{10}^{010}                \cr 
}
\ess\ess\ess ,\ess\ess\ess
\matrix{
{\rm by}\ess{\bf flip\ssp}_{\trois311}  &    &{\rm  Spaces} &     &{\rm by}\ess {\bf flip\ssp}_{\deux32} \cr               
                                  &    &              &     &                                  \cr               
    b_{01}^{010}                  &\rar&   \BM^{111}  &\lar &      b_{10}^{100}                \cr 
}
$$
Now the non vanishing of $b_{10}^{100}\scs b_{10}^{010}$ and $ b_{01}^{010}\scs  b_{01}^{001}$ is perfectly
consistent with the equations in 2.17. So it is pretty clear that we can't choose  at the same time
$
\CB_{10}= \CB_1\cup\CB_2
$
and
$
\CB_{01}= \CB_2\cup\CB_3\ess .
$
On the other hand if we take the subset $\CB^{001}+\CB^{010}$ out of $\CB_{01}$
or the subset  $\CB^{100}+\CB^{010}$ out of $\CB_{10}$ then equation b) and d) will do the
rest and force $b_{01}=b_{10}=0$. It turns out that either choice leads to the construction of a basis
for $\BM_{\trois321}$. To be definite we shall take
$$
\CB_{10}= \CB_1\cup\CB_2
\ess\ess\ess {\rm and } \ess\ess\ess
\CB_{01}= \CB^{110}+\CB^{111}+\CB^{011}+\CB^{101}
\eqno 2.18
$$  
Our choices so far force $b_{00}=b_{10}=b_{01}=0$ and the system 2.14 now reduces 
to the single equation
$$
 -\ssp 2  b_{02}(\del)\Delta_{\trois221} \sms  b_{11}(\del)\Delta_{\trois311}\sps 2 b_{20}(\del)\Delta_{\deux32} \ses 0\ess .
\eqno 2.19
$$
Remarkably, we can chose $\CB_{02}\scs \CB_{11}$ and $\CB_{20}$ in an optimal way and still guarantee
that this single equation will force $b_{02}=b_{11}=b_{20}=0$. Note that at this stage, the maximal choices
for $\CB_{02}\scs \CB_{11}$ and $\CB_{20}$ are $\CB_3\scs \CB_2$ and $\CB_1$ respectively. 
Adopting the {\ita greedy algorithm} strategy we chose
$$
\CB_{20}\ses \CB_1\ess .
\eqno 2.20
$$
In order that the second term in 2.19 does not produce a component that could be canceled by the third
term we need  $b_{11}(\del)\Delta_{\trois311}$ to fall out of $\BM_{\deux32}$. A look at Figure 2.
suggests that $\CB_{11}$ should then lie in $\BM^{110}\hskip -.03truein\oplus\BM^{111}$. The  maximal choice
is then
$$
\CB_{11}\ses \CB^{110}\ssp +\ess \CB^{111}\ess .
$$
This done we must assure that $ b_{02}(\del)\Delta_{\trois221}$ falls out of $\BM_{\trois311}$ and $\BM_{\deux32}$.
Again from Figure 2. we get that we must take
$$
\CB_{02}\ses  \CB^{111}\ess .
$$
With this final choice we have assured that the system in 2.17 forces all $b_{ij}$ to vanish
and consequently establish the independence of the system given in  2.11. We are left with 
checking that its cardinality is indeed $6!$. Now our choices have been
$$
\matrix{
Basis  &&  Cardinality \cr
&&\cr
\CB_{00}=\CB_1\cup\CB_2\cup\CB_3=\CB_1 + \CB_2\cap (\CB\backslash \CB_1) +\CB^{001}
 &&  5!\sps {5!\over 2}\sps {5!\over 3}\cr
&&\cr
\CB_{10}=\CB_1+\CB^{010}+\CB^{011}  && 5!\sps {5!\over 3}\sps {5!\over 6}\cr
&&\cr
\CB_{01}=\CB^{110}+\CB^{111}+\CB^{101}+\CB^{011}  && 3  \times  {5!\over 6}\sps {5!\over 3}\cr
&&\cr
\CB_{20}=\CB_1   && 5! \cr
&&\cr
\CB_{11}= \CB^{110}+\CB^{111}  && {5!\over 2}\cr
&&\cr
\CB_{02}=\CB^{111}  &&   {5!\over 3} \cr
} 
$$
Here all these cardinalities result from the assumption that the $n!$ conjecture 
holds true for $\trois221\scs \trois311 \scs \deux32$ and that the dimensions
of the spaces $\BM^{\eee_1\eee_2\eee_3}$ are as given by 2.8. 
As we can see the count is indeed $6\times 5!$ as desired. Thus, to prove 
the $n!$ conjecture for the partition $(3,2,1)$ and obtain a basis for the
module $\BM_{321}$ we need only exhibit the collection $\CB$ with properties (i)-(iv)
of the SF heuristic.

Using the dimensions given
in 2.8 we can easily see that $\CB$ will consist of $5!\times 13/6=260$ elements.
This given, we shall limit ourselves here to exhibiting the basic ingredients from
which $\CB$ may be easily constructed. Note first that we need only exhibit $\CB^{111}$
and $\CB^{110}$ since (via 2.1) all the other $\CB^{\eee_1\eee_2\eee_3}$ may be recovered by
successive flips. Secondly, note that if we obtain a complete decomposition of 
the modules $\BM^{111}$ and $\BM^{110}$ into irreducible constituents, then we
need only exhibit a collection of cyclic elements generating these constituents.
In fact, if the cyclic elements are given as Young idempotents acting on 
polynomials, then $\BM^{111}$ and $\BM^{110}$ may be constructed by selecting from
the $S_5$-orbits of these cyclic elements those corresponding to standard tableaux
idempotents. In the following tables we list such a  set of cyclic elements as
sums of products of determinants. This description is a special case of
a construction of basis for these intersections given by F. Bergeron and S. Hamel 
in [\bergeronhamel]. We recall that if $a_1<a_2<\cdots<a_m$ is a
sequence of integers then Young uses the symbol $[a_1,a_2,\ldots ,a_m]$ to denote the
formal sum of all elements of the symmetric group $S_{\{a_1,a_2,\ldots ,a_m\}}$.
Likewise Young uses the symbol $[a_1,a_2,\ldots ,a_m]'$ to denote the 
sum of the elements of $S_{\{a_1,a_2,\ldots ,a_m\}}$ multiplied their sign.
In this notation, the Young idempotent correspoding to the standard tableau 
  $$ T=\matrix{\Young{6 \cr
            4& 5 \cr
            1& 2 &3 \cr}\cr}$$
is simply given by the group algebra expression
$$
[1,4,6]'[2,5]'\  [1,2,3][4,5]  
$$ 
Young calls $[1,2,3][4,5]$ the $\underline {\rm row\ssp   group}$ of $T$ and $[1,4,6]'[2,5]'$
the $\underline {\rm signed \ssp column \ssp  group}$. Young in [\youngsubst] makes the following 
important observation. Namely, if $\gamma$ is any element of the group algebra of  $S_n$
and $V$ is a vector space spanned by some symbols on which $S_n$ acts then, for any $v\in V$,
the action of $S_n$ on the submodule generated by the element $\gamma\ssp v$  is identical
to the action of $S_n$ on the left ideal generated by the group algebra element 
$\gamma\ssp \Sigma_v$ where $\Sigma_v$ denotes the formal sum of the elements
of $S_n$ that stabilize $v$. This observation should be kept in mind when reading
our table below. For instance, suppose we are to study the representation
resulting from the action of $S_6$  on the polynomial 
$$
P(x;y)= 
{\rm det} \left[\matrix{1 & 1 & 1\cr y_1& y_2 &y_3\cr x_1& x_2 &x_3\cr }\right]
{\rm det}\ssp \left[ \matrix{1 & 1 \cr y_4& y_5 \cr }\right]
\ess .
$$ 
We first notice that it may be written in the form
$$
P(x;y)=[1,2,3]'[4,5]'\ssp y_2\, x_3 \,y_5\ess .
$$
Then observe that, since the stabilizer of $y_2\, x_3\, y_5$ in $S_6$ is $[1,4,6]\ [2,5]$, 
we may also write
$$
P(x;y)={1\over 12} \ess [1,2,3]'[4,5]'[1,4,6][2,5]\ssp y_2 x_3 y_5\ess .
\eqno 2.21
$$
This makes it evident that $P(x,y)$ generates an irreducible representation of $S_6$
with character $\chi^{3,2,1}$. Note further that, since $P$ is a bihomogeneous polynomial
of $x$-degree $1$ and $y$-degree $2$, the submodule generated by $P$ would contribute
the term $tq^2S_{321}$ to the bivariate Frobenius characteristic of any bigraded
$S_6$-module that contains $P$. It will be convenient to refer to $tq^2$ as the 
{\ita weight} of this representation and to the partition $321$ as its {\ita shape}.
Finally, we also immediately read from 2.21 that a basis for the submodule
generated by $P$ is given by $P$ and the $15$ polynomials obtained by replacing in 2.21
the idempotent $[1,2,3]'[4,5]'[1,4,6][2,5]$ by the idempotents corresponding
to the other $14$ standard tableaux of shape $321$. Another example may be helpful
in checking some of the entries in the tables below. Suppose that we are to 
study the $S_5$-module generated by the action of $S_5$ on the polynomial
$$
Q(x,y)\ses \det \pmatrix{1 & 1 & 1\cr y_1^2& y_2^2 &y_3^2\cr x_1& x_2 &x_3\cr }
\sms 2\ssp y_5 \ssp \det\pmatrix{1 & 1 & 1\cr y_1& y_2 &y_3\cr x_1& x_2 &x_3\cr }
$$ 
Then it is easily seen that we have
$$
[4,5]'Q(x,y)\ses -2\ssp\det\pmatrix{1 & 1 & 1\cr y_1& y_2 &y_3\cr x_1& x_2 &x_3\cr }
\det\pmatrix{1 & 1 \cr y_4& y_5 \cr} 
\ses -{1\over 2}\ssp [1,2,3]'[4,5]'[1,4][2,5]\ssp y_2x_3y_5 
$$
as well as
$$
\eqalign{
{1\over 2}\ssp [4,5]Q(x,y) &\ses
\det\pmatrix{1 & 1 & 1\cr y_1^2& y_2^2 &y_3^2\cr x_1& x_2 &x_3\cr }
\sms \ssp \det\pmatrix{1 & 1 & 1\cr y_1& y_2 &y_3\cr x_1& x_2 &x_3\cr }
\ess (y_4+y_5) \cr\cr 
& \ses
[1,2,3]'\bigl(\ssp y_2^2x_3\sms y_2x_3\ssp( y_1+y_4+y_5)\ssp \bigr)\cr\cr 
&\ses
{1\over 6}\ssp [1,2,3]'[1,4,5]\bigl(\ssp y_2^2x_3\sms y_2x_3\ssp( y_1+y_4+y_5)\ssp \bigr)\ess .\cr} 
\eqno 2.22
$$
\noindent
We thus see that $Q$ generates a bigraded $S_5$-module with bivariate Frobenius 
characteristic
$$
t\ssp q^2 S_{2,2,1}\sps t\ssp q^2  S_{3,1,1}\ess .
\eqno 2.23
$$
This given, it should not  be difficult to  check that the following two tables
give all the cyclic elements needed to generate our bases $\CB^{111}$ and $\CB^{110}$.
\sap

\vbox{
\centerline {\bol TABLE FOR $\CB^{111}$} 
$$
\matrix
{
 {\bf\bigsp Cyclic\ess element}\ess\bigsp &\bigsp{\bf Weight}
\bigsp &   {\bf Shape }
\cr
\cr
\left\{
\matrix
{
\ess {1\over 2}(\del_{x_5}+\del_{x_4})\del_{y_4} \del_{y_5}\Delta_{\deux32}
{{\ \atop \ }  \atop {\ \atop \ } }\cr  
{1\over 4}(\del_{y_4}\del_{x_5}^2+\del_{y_5}\del_{x_4}^2)\Delta_{\trois311} \cr
{1\over 2}\del_{x_1}\del_{x_2}\del_{x_3}\Delta_{\trois221} \cr
}
\right\}\ses
{\rm det}\left[ \matrix{1 & 1& 1\cr y_1 & y_2& y_3\cr x_1 & x_2 & x_3 \cr}\right]
{\rm det}\left[ \matrix{1 & 1\cr y_4 & y_5\cr}\right]
&q^2t
& \trois221
\cr
\cr
\cr\left\{
\matrix
{
\ess -{1\over 2}\del_{y_1}\del_{y_2}\del_{y_3}
\Delta_{\deux32}
{{\ \atop \ }  \atop {\ \atop \ } }\cr  
-{1\over 4}(\del_{x_4}\del_{y_5}^2+\del_{x_5}\del_{y_4}^2)\Delta_{\trois311} \cr
{1\over 2}(\del_{y_5}+\del_{y_4})\del_{x_4} \del_{x_5}\Delta_{\trois221} \cr
}
\right\}\ses
{\rm det}\left[ \matrix{1 & 1& 1\cr y_1 & y_2& y_3\cr x_1 & x_2 & x_3 \cr}\right]
{\rm det}\left[ \matrix{1 & 1\cr x_4 & x_5\cr}\right]
&q t^2
& \trois221
\cr
\cr
\cr
\left\{
\matrix
{
{1\over 2}\del_{y_3}\del_{y_4}\del_{y_5}^2\Delta_{\deux32} \cr
{1\over 4}(\del_{y_3}\del_{x_4}+\del_{x_3}\del_{y_4})\del_{y_5}^2\Delta_{\trois311} \cr
{1\over 2}\del_{y_3}\del_{y_4}\del_{x_5}^2\Delta_{\trois221} \cr
}
\right\}=
{\rm det}\left[ \matrix{1 & 1\cr x_1 & x_2  \cr}\right] 
{\rm det}\left[ \matrix{1 & 1\cr x_3 & x_4  \cr}\right]
&t^2
& \deux32
\cr
\cr
\cr
\cr
}
$$
}

\vbox{
\centerline {\bol TABLE FOR $\CB^{111}$ {\ita continued}} 

$$
\matrix
{
 {\bf\bigsp Cyclic\ess element}\ess\bigsp &\bigsp{\bf Weight}
\bigsp &   {\bf Shape }
\cr
\cr
\left\{
\matrix
{
{1\over 2}(\del_{y_1}\del_{y_2}+\del_{y_3}\del_{y_4}) \del_{x_5}\del_{y_5}\Delta_{\deux32} \cr
{1\over 4}(\del_{x_1}\del_{x_2}+\del_{x_3}\del_{x_4})\del_{y_5}^2\Delta_{\trois311} \cr
-{1\over 2}(\del_{x_1}\del_{x_2}+\del_{x_3}\del_{x_4}) \del_{x_5}\del_{y_5}\Delta_{\trois221} \cr
}
\right\}\ses
[1,2]'[3,4]'(y_2x_4+x_2y_4)
&qt
& \deux32
\cr
\cr
\cr
\left\{
\matrix
{
{1\over 2}\del_{x_3}\del_{x_4}\del_{y_5}^2\Delta_{\deux32} \cr
{1\over 4}(\del_{y_3}\del_{x_4}+\del_{x_3}\del_{y_4})\del_{x_5}^2\Delta_{\trois311} \cr
{1\over 2}\del_{x_3}\del_{x_4}\del_{x_5}^2\Delta_{\trois221} \cr
}
\right\}=
{\rm det}\left[ \matrix{1 & 1\cr y_1 & y_2  \cr}\right] 
{\rm det}\left[ \matrix{1 & 1\cr y_3 & y_4  \cr}\right]
&q^2
& \deux32
\cr\cr\cr
\left\{\matrix{
{1\over 2}\del_{x_4}\del_{y_4}\del_{y_5}^2\Delta_{\deux32} \cr
{1\over 4}\del_{x_5}^2\del_{y_4}^2\Delta_{\trois311} \cr
{1\over 2}\del_{x_4}\del_{y_4}\del_{x_5}^2\Delta_{\trois221} \cr
}\right\}=
{\rm det}\left[ \matrix{1 & 1& 1\cr y_1 & y_2& y_3\cr x_1 & x_2 & x_3 \cr}\right]
&qt
& \trois311
\cr
\cr
\cr
{\rm det}\left[ \matrix{1 & 1\cr x_1 & x_2  \cr}\right]
&t 
& \deux41
\cr
\cr
\cr
{\rm det}\left[ \matrix{1 & 1\cr y_1 & y_2\cr}\right]
& q
& \deux41
\cr
\cr
\cr
1
& 1
& \un5
\cr
\cr
}
$$
}
\sa
\vbox{
\centerline {\bol TABLE FOR $\CB^{011}$} 

$$
\matrix{
{\bf\bigsp Cyclic\ess element}\ess\bigsp &\bigsp{\bf Weight}\bigsp &   {\bf Shape }
\cr
\cr
\left\{
\matrix
{
[3,4,5]'
(\del_{y_3}\del_{x_4}\del_{x_5}\del_{y_5})\ssp \deux32\cr
\cr
{1\over 2}
[3,4,5]'
(\del_{y_3}\del_{x_4}\del_{x_5}^2)\ess\trois311\cr
}
\right\}
\ses [1,2]' 
(-6 y_2^2 +4( y_1 + y_3+ y_4 + y_5) y_2)
&
q^2
&
\deux41 
\cr
\cr
\cr\cr
\left\{
\matrix
{
{1\over 2}
(\del_{y_4}-\del_{x_5})\del_{x_4}\del_{y_5}\ssp \deux32 \cr
{1\over 4}
(\del_{x_4}^2\del_{y_5}-\del_{x_5}^2 \del_{y_4})\ssp \trois311
\cr
}
\right\}
\ses
[1,2,3]'\ssp \bigl(-y_2^2x_3+2(y_1+y_4+y_5)y_2x_3\bigr)
&
q^2 t
&
\trois311
\cr\cr\cr
\left\{
\matrix
{
\ess \del_{x_5}\del_{y_5}
\Delta_{\deux32}
{{\ \atop \ }  \atop {\ \atop \ } }\cr  
{1\over 2}\del_{x_5}^2\Delta_{\trois311} \cr
}
\right\}\ses
{\rm det}
\left[ 
\matrix{1 & 1& 1&1\cr y_1 & y_2& y_3&y_4\cr y_1^2 & y_2^2& y_3^2&y_4^2\cr x_1 & x_2 & x_3&x_4 \cr}
\right]
&q^3 t
& \quatre2111
\cr
\cr
\cr
\left\{
\matrix
{
\del_{x_4}\del_{x_5}\del_{y_5}\ssp \deux32 \cr
{1\over 2}
\del_{x_4}\del_{x_5}\del_{y_5}\ssp \trois311 \cr
}
\right\}\ses
\left[
\matrix
{ 
1&1&1\cr
y_1&y_2&y_3\cr
y_1^2&y_2^2&y_3^2\cr
}
\right]
&
q^3
&
\trois311
}
$$
}

\sas

At this point it is good to see to what extent the results of this section
are consistent with the explicit calculations made by Macdonald in the original
paper [\macdonald]. We are referring here to the tables of the coefficients $K_{\la\mu}(q,t)$
which are also given in [\macdonaldbook]. After making the changes of scale
$$
K_{\la\mu}(q,t)\ess \RA \ess\TK_{\la\mu}(q,t)=T_\mu\ssp K_{\la\mu}(q,1/t)\ess ,
$$
and using the expansion in I.8 we can easily compute all of the polynomials
$\TH_\mu(x;q,t)$ for all $|\mu|\leq 6$. This given, let us subtract the the second equation
in 2.7 multiplied by $1/T_2$ from the first multiplied by $1/T_1$, obtaining
$$
{1\over T_1}\TH_{\aaa^{(1)}}\sms{1\over T_2}\ssp \TH_{\aaa^{(2)}}
\ses 
\Bigl({1\over T_1} \sms {1\over T_2}\Bigr)
\Bigl({1\over T_3} \Phi^{(2)} \sps \Phi^{(3)} \ssp \Bigr)\ess  .
$$ 
Which gives
$$
\eqalign{
\Bigl({1\over T_1}\TH_{\aaa^{(1)}}\sms{1\over T_2}\ssp \TH_{\aaa^{(2)}}\Bigr)
\Big/
\Bigl({1\over T_1} \sms {1\over T_2}\Bigr)  
&\ses
 {1\over T_3} \Phi^{(2)} \sps \Phi^{(3)}  \cr
&\ses  \Phi_\mu^{110}+\Phi_\mu^{111}
\ses {\cal F}\bigl( \BM_{\aaa^{(1)}}\cap \BM_{\aaa^{(2)}}\bigr)\ess .
}
\eqno 2.24
$$ 
Similarly, the second and third  equations in 2.7 give
$$
\eqalign{
\Bigl({1\over T_2}\TH_{\aaa^{(2)}}\sms{1\over T_3}\ssp \TH_{\aaa^{(3)}}\Bigr)
\Big/
\Bigl({1\over T_2} \sms {1\over T_3}\Bigr)  
&\ses
 {1\over T_1} \Phi_\mu^{(2)} \sps \Phi_\mu^{(3)}\cr  
&\ses \Phi_\mu^{011}+\Phi_\mu^{111}
\ses {\cal F}\bigl( \BM_{\aaa^{(2)}}\cap \BM_{\aaa^{(3)}}\bigr)\ess .\cr
}\eqno 2.25
$$
Subtracting 2.25 multiplied by $1/T_3$ from 2.24  multiplied by $1/T_1$ 
and dividing the result by  
$ \bigl({1\over T_1} -{1\over T_3}\bigr)$ 
gives
$$
 \Phi^{(3)}=
\ses{
{1\over T_1}\Bigl({1\over T_1}\TH_{\aaa^{(1)}}-{1\over T_2}\ssp \TH_{\aaa^{(2)}}\Bigr)
\Big/
\Bigl({1\over T_1} - {1\over T_2}\Bigr) 
\sms
{1\over T_3}\Bigl({1\over T_2}\TH_{\aaa^{(2)}}-{1\over T_3}\ssp \TH_{\aaa^{(3)}}\Bigr)
\Big/
\Bigl({1\over T_2} - {1\over T_3}\Bigr)
\over 
\Bigl({1\over T_1} \sms{1\over T_3}\Bigr)
}
\ess .
\eqno 2.26
$$ 
Using Macdonald's tables to compute the left-hand side of 2.26
for $\ssp \aaa^{(1)}= \deux32\ssp $, $\ssp \aaa^{(2)}= \trois311\ssp $
and $\ssp \aaa^{(3)}= \trois221\ssp\ssp $ gives
$$
\eqalign{
&
{{1\over q^2}\Bigl({1\over q}\TH_{\deux32}\sms{1\over t }\ssp \TH_{\trois311}\Bigr)
\Big/
\Bigl({1\over q} \sms {1\over t }\Bigr) 
\sms
{1\over t^2 }\Bigl({1\over q}\TH_{\trois311}\sms{1\over t }\ssp \TH_{\trois221}\Bigr)
\Big/
\Bigl({1\over q} \sms {1\over t }\Bigr) 
\over {1\over q^2}\sms {1\over t^2 }
}
\ses\cr
&\ess\ess\ess\ess\ess\ess\ess\ess\ess\ess\ess\ess\ess\ess\ess\ess
\ses \Phi_{321}^{111} \ses S_4+(t+q)S_{41}+(t^2+tq+q^2)S_{32}+tqS_{311}+(qt^2+tq^2)S_{221}
} 
$$
which is easily seen to be in perfect agreement with the Frobenius characteristic
of 
$$
\bigl( \BM_{\deux32}\cap \BM_{\trois311}\cap\BM_{\trois221}\bigr)
$$
as may be put together from our tables giving $\CB^{111}$. Similarly, from 2.24
we derive that
$$
\Phi_{321}^{110}=
\Bigl({1\over q}\TH_{\deux32}\sms{1\over t }\ssp \TH_{\trois311}\Bigr) 
\Big/
\Bigl({1\over q} \sms {1\over t }\Bigr)
 \sms
\Phi_{321}^{111}
\ses q^2S_{41}+q^2(t+q)\ssp S_{311}+q^3t\ssp S_{2111}\ess .
$$
Which is again easily seen to agree with our table for $\CB^{110}$.
\sas

We shall see  in the next section that divided difference formulas as in 2.24, 2.25
and 2.26 are but special  cases of a general identity giving the Frobenius
characteristic of the intersection of any subset of the modules
$\BM_{\aaa^{(1)}},\BM_{\aaa^{(2)}},\ldots ,\BM_{\aaa^{(m)}}$.

\section{3. General identities}
We begin by extending the arguments that gave 2.4 and 2.5 
to the general case. Here and in the following we  adopt the
same notational conventions we made in the introduction. We  work
with a fixed partition $\mu$ with $m$ corners, with predecessors
$\aaa^{(1)},\aaa^{(2)},\ldots ,\aaa^{(m)}$ ordered as we indicated in the 
introduction. We shall assume that the SF heuristic is valid  for $\mu$
and that a basis $\CB$ has been constructed with the required properties.

\Prop{3.1}
For every $1\leq k\leq m$ there exists a Schur positive function $\Phi_\mu^{(k)}$
such that for all words $\eee=\eee_1\eee_2\cdots\eee_m$ with $\sum_{i=1}^m\eee_i=k$
we have
$$
\Phi_\mu^{\eee}\ses {{\Phi^{(k)} }\over \prod_{\eee_i=0} T_i }
\eqno 3.1
$$

\Proof For a moment let us set
$$
\Psi_\mu^{\eee}\ses T^{\tilde \eee}\ssp \Phi_\mu^{\eee}
\eqno 3.2
$$
with
$$
T^{\tilde \eee}\ses \prod_{i=1}^m T_i^{\teee_i }\ses \prod_{\eee_i=0} T_i \ess .
\eqno 3.3
$$
This given, to show 3.1 we need only establish that $\Psi_\mu^{\eee}$ depends only
on $\sum_{i=1}^m\eee_i$. Note that if $\eee=1111\cdots 1=1^m$ then we may take 
$\Phi^{(m)}:= \Phi_\mu^{1^m}$. In all other cases $\eee$ will contain at least a $0$
and a $1$. Choose a pair $i,j$ such that $\eee_i=1$ and $\eee_j=0$.
Note that two applications of I.30 yield that
$$
\da_j\da_i\ssp \Phi_\mu^{\eee} \ses\Phi_\mu^{\tau_j\tau_i\eee}\ess .
\eqno 3.4
$$
However, from the definition I.28 we immediately derive that (if $i<j$)
$$
\tau_j\tau_i\eee\ses 
\eee_1\cdots \eee_{i-1}{\bf 0}\ssp\eee_{i+1}\cdots \eee_{j-1}{\bf 1}\ssp\eee_{j+1}\cdots \eee_m
\ses (i,j) \ssp \eee\ess ,
\eqno 3.5
$$
where $(i,j)$ denotes the transposition that interchanges the $i^{th}$ and $j^{th}$
letters of a word. On the other hand I.27 gives
$$
\da_j\da_i\ssp \Phi_\mu^{\eee} \ses\da_j(T_i\ssp \da \Phi_\mu^{\eee})
\ses {T_j\over T_i}\ess\Phi_\mu^{\eee}\ess . 
\eqno 3.6
$$
Combining 3.4, 3.5 and 3.6 we get 
$$
\Phi_\mu^{(i,j)\eee}\ses {T_j\over T_i}\ess \Phi_\mu^{\eee}\ess .
\eqno 3.7
$$
Since $\teee_i=0$ and $\teee_j=1$ we can write
$$
T^{(i,j)\teee}\ses {T_i\over T_j}\ess T^{\teee}\ess .
$$
Thus multiplying both sides of 3.7 by $T^{(i,j)\teee}$ we finally obtain
$$
\Psi_\mu^{(i,j)\eee}\ses T^{(i,j)\teee}\ess \Phi_\mu^{(i,j)\eee}\ses
T^{\teee}\ess \Phi_\mu^{\eee}\ses \Psi_\mu^{\eee} \ess .
$$
Thus $\Psi_\mu^{\eee}$ remains unchanged when we arbitrarily permute
the letters of $\eee$, proving the proposition.
\sa

To proceed we need some further notational convenions. To begin with, note that
we can associate to the word $\eee=\eee_1\eee_2\cdots \eee_m$ the subset
$S=\{1\leq i\leq m:\eee_i=1\ssp \}$. Conversely, for $S\con\{1,2,\ldots ,m\}$ 
we let $\eee(S)$ be the word $\eee$ such that $\eee_i=1$ iff $i\in S$. This given, it
will be convenient to set
$$
\Phi_\mu^{=S}:=  \Phi_\mu^{\eee(S)}\ess .
$$
We should note that this notation is consistent with the fact that
$$
\Phi_\mu^{\eee(S)}\ses {\cal F} 
\bigl(\bigcap_{i\in S}\BM_{\aaa^{(i)}}\cap \bigcap_{j\notin S}\BM_{\aaa^{(j)}}^\perp\bigr)
\ess .
$$
In the same vein we let
$$
\Phi_\mu^{\supseteq S}\ses {\cal F} 
\bigl(\bigcap_{i\in S}\BM_{\aaa^{(i)}}\bigr)\ess . 
$$
Note that since
$$
\bigcap_{i\in S}\BM_{\aaa^{(i)}}\ses \bigoplus_{T\noc S}\ssp \BM^{\eee(T)}\ess ,
$$
we deduce that
$$
\Phi_\mu^{\supseteq S}\ses \sum_{T\noc S}\ssp \Phi_\mu^{=T}\ess .
\eqno 3.8
$$

\Prop{3.2}
For every subset $S\con \{1,2,\ldots ,m\}$ we have
$$
\Phi_\mu^{\supseteq S}\ses \sum_{k=0}^m\ssp \P k
e_{m-k}\Bigl[\ssp
{1\over T_1}+{1\over T_2}+\cdots +{1\over T_m}\sms 
\sum_{i\in S}\ssp{1\over T_i}\ssp \Bigr ]\ess . 
\eqno 3.9
$$

\Proof Using formula 3.1, the identity in 3.8 reduces to
$$
\Phi_\mu^{\supseteq S}\ses \sum_{k=0}^m\ssp \P k \sum_{T\noc S\ssp \&  |T|=k}\ssp  
\prod_{i\notin T}{1\over T_i}
$$
and this may also be rewritten as in 3.9.
\sas

It will be convenient to denote by $\BS_\mu$ the vector space spanned by the 
symmetric functions $\TH_{\aaa^{(i)}}$ for $i=1,\ldots,m\ssp .\ess $  In symbols
$$
\BS_\mu\ses \CL[\TH_{\aaa^{(1)}},\TH_{\aaa^{(2)}}\,\ldots ,\TH_{\aaa^{(m)}}]
\eqno 3.10
$$
We will make extensive use here of the operator $\nabla$, acting
on symmetric polynomials, which gives
$$
\nabla\ssp \TH_\la\ses T_\la\ssp \TH_\la
\bigsp (\ssp \forall \ess \la\ssp)\ess .
\eqno 3.11
$$
Since it can be shown (see [\garsiahaimanorbit]) that the symmetric polynomials $\TH_\la$ ($\la\vdash n$) form a basis (of the space
of hogeneous symmetric functions of degree $n$),
formula 3.11 may be used as the definition of $\nabla$. 

We  should note that, using the notation adopted in the introduction, we can write
$$
\nabla  \TH_{\aaa^{(i)}}\ses T_i\ess  \TH_{\aaa^{(i)}}\ess .
\eqno 3.12
$$
Note further that if $a_i$ and $l_i$ respectively denote the coarm and coleg of
the cell we must remove from $\mu$ to get $\aaa^{(i)}$ then setting
$$
x_i=t^{l_i}q^{a_i}\bigsp\bigsp (\ssp i=1,\ldots ,m\ssp )
\eqno 3.13
$$
we have
$$
T_i\ses T_\mu/x_i\ess .
\eqno 3.14
$$
In particular we see that the monomials $T_1,T_2,\ldots ,T_m$ are all distinct.
This enables us to write every element of $\BS_\mu$ as a polynomial in $\nabla$ applied
to a single element of $\BS$. In fact, any element 
$$
\Phi\ses \sum_{j=1}^m\ssp c_j\ess \TH_{\aaa^{(j)}}
$$
with all $c_i\neq 0$ may be used for this purpose. However, the most convenient one is
the element
$$
\Phi_\mu\ses  \sum_{j=1}^m\ssp \biggl( \prod_{s=1,s\neq j}^m\ssp
{1\over  \bigl(1-{T_j/ T_s}\bigr)}\biggr)\ess \TH_{\aaa^{(j)}} \ess ,
\eqno 3.15
$$
for, we have the following 
\sas
\Prop{3.3}
If
$$
\Phi\ses \sum_{j=1}^m\ssp c_j\ess \TH_{\aaa^{(j)}}
$$
then
$$
\Phi\ses P(\nabla)\ssp \Phi_\mu
$$
with
$$
P(\nabla)\ses \sum_{i=1}^m\ssp c_i\ssp \prod_{s=1,s\neq i}^m(1-\nabla/T_s)
$$

\Proof From 3.12 it follow that
$$
\prod_{s=1,s\neq i}^m\ssp \Bigl(1-{\nabla\over T_s}\Bigr)\ess
\Phi_\mu=  \biggl( \prod_{s=1,s\neq i}^m\ssp
{1\over  \bigl(1-{T_i/ T_s}\bigr)}\biggr)\ess\prod_{s=1,s\neq i}^m
\ssp \Bigl(1-{T_i\over T_s}\Bigr)\ess \TH_{\aaa^{(i)}}=
\ess  
\TH_{\aaa^{(i)}}\ssp , 
\eqno 3.16
$$
and the result follows by linearity.
\sas

Surprisingly, it develops that $\Phi_\mu$ is none other than $\P m $ itself, namely
the polynomial giving the bigraded Frobenius characteristic of the intersection
of the modules $\BM_{\aaa^{(i)}}$. More generally, we have the following remarkable
identities:

\Thm{3.1}
For  $k=1,..,m$  
$$
\P k \ses (-\nabla)^{m-k}\ssp \Phi_\mu\ess ,
\eqno 3.17
$$
and for any $S\con\{1,2,\ldots ,m\}$ 
$$
\Phi^{\noc S}\ses \prod_{i\notin S}\ssp \Bigl(\ssp 1-{\nabla\over T_i}\ssp \Bigr)\ess \Phi_\mu\ess .
\eqno 3.18
$$
We also have for $i=1,..,m$ 
$$
\eqalign{
&a)\ess\ess\TH_{\aaa^{(i)}}(x;q,t)\ses 
\prod_{s=1,s\neq i}^m\ssp \Bigl(1-{\nabla\over T_s}\Bigr)\ess
\Phi_\mu\ess ,
\cr
&b)\ess\ess\TH_{\aaa^{(i)}}(x;q,t)\ses \sum_{k=1}^m\ssp \P k \ess
e_{m-k}\Bigl[\ssp {1\over T_1}+{1\over T_2}+\cdots +{1\over T_m}\sms 
{1\over T_i}\ssp \Bigr ]\ess .\cr
}
\eqno 3.19
$$
Thus the $\P k $ form a basis for $\BS_\mu$.
They may be also be explicitely computed by means of the following formula:
$$
\P k \ses  \sum_{j=1}^m\ssp \biggl( \prod_{s=1,s\neq j}^m\ssp
{1\over  \bigl(1-{T_j/ T_s}\bigr)}\biggr)\ess (-T_j)^{m-k}\ess \TH_{\aaa^{(j)}} \ess .
\eqno 3.20
$$

\Proof Formula 3.19 a) is 3.16 and 3.19 b) is the special case $S=\{i\}$ of 3.9. 
This not only implies that 
the $\P k $ are a basis for $\BS_\mu$ but also that the matrix
$$
\Big\|
e_{m-k}\Bigl[\ssp {1\over T_1}+\cdots +{1\over T_m}\sms {1\over T_i}\ssp 
\Bigr ]\Big\|_{k,i=1}^m
\eqno 3.21
$$
is non-singular. Expanding the right-hand side of 3.19 b) as a polynomial in $\nabla$ we derive 
that
$$
\TH_{\aaa^{(i)}}(x;q,t)\ses \sum_{k=1}^m\ssp (-\nabla)^{m-k}\ssp \Phi_\mu \ess
e_{m-k}\Bigl[\ssp {1\over T_1}+{1\over T_2}+\cdots +{1\over T_m}\sms 
{1\over T_i}\ssp \Bigr ]\ess .
$$
Comparing this with 3.19 b) and using the non-singularity of the matrix in 3.21,
we derive that 3.17 must hold true for all $k=1,..,m$. This given, formula 3.20
follows immediately by applying $(-\nabla)^{m-k}$ to both sides of 3.15
and using 3.12. In the same vein we see that, using 3.17, 3.9 can be rewritten in the form
$$
\Phi_\mu^{\noc S}\ses \sum_{k=0}^{m-|S|}\ssp (-\nabla)^k\ssp 
e_k\Bigl[\sum_{i\notin S}{1\over T_i}\Bigr]\ess ,
$$
This gives 3.18 and completes our proof.
\sa

As a corollary of 3.18 we obtain a recursive way of computing $\Phi_\mu$ as well as
any of the $\Phi_\mu^{\noc S}$. 

\noindent
For given quantities $y_1,y_2,y_3,\ldots\ssp ; $
$A_1,A_2,A_3,\ldots $ we recursively set (for $k>2$)
$$
\Delta[y_1,\ldots, y_k;A_1,\ldots,A_k]\ses 
{
y_1\ssp \Delta[y_1,\ldots, y_{k-1};A_1,\ldots,A_{k-1}]\sms
\Delta[y_2,\ldots, y_{k};A_2,\ldots,A_{k}]\ssp y_k
\over 
y_1\sms y_k
}
$$
with
$$
\Delta[y_1,y_2;A_1,A_2]\ses 
{y_1\ssp A_1\sms A_{2}\ssp y_2
\over 
y_1\sms y_2
}
\bigsp\bigsp (\ess {\rm for} \ess\ess k=2\ess )
$$
\sas

\Prop{3.4}
For $S=\{1\leq i_1<i_2<\cdots <i_k\leq m\}$ we have
$$
\Phi_\mu^{\noc S}\ses 
\Delta\bigl[{1\over T_{i_1}},\ldots ,{1\over T_{i_k}};\TH_{\aaa^{i_1}},\ldots,\TH_{\aaa^{i_k}}\bigr]
\eqno 3.22
$$

\Proof Let $T$ be a subset of $\{1,2,\ldots ,m\}$ and let $i$ be smaller and $j$ bigger
than any element of $T$, then note that 
$$
\eqalign{
{1\over T_i}\prod_{s\notin T\cup\{i\}} \Bigl(1-{\nabla\over T_s}\Bigr)
- 
{1\over T_j}\prod_{s\notin T\cup\{j\}} \Bigl(1-{\nabla\over T_s}\Bigr)
&\ssp =\ssp \biggl({1\over T_i}\Bigl(1-{\nabla\over T_j}\Bigr)- 
{1\over T_j}\Bigl(1-{\nabla\over T_i}\Bigr)\biggr)
\prod_{s\notin T\cup\{i,j\}}  \Bigl(1-{\nabla\over T_s}\Bigr)
\cr
&\ssp =\ssp \biggl({1\over T_i}\sms {1\over T_j}\biggr)\prod_{s\notin T\cup\{i,j\}}  
\Bigl(1-{\nabla\over T_s}\Bigr)\ess .\cr
}
\eqno 3.23
$$
Applying both sides to $\Phi_\mu$ and using 3.18 gives
$$
\biggl(
{1\over T_i}\ess \Phi_\mu^{\noc T\cup\{i\}}\sms {1\over T_j}\ess \Phi_\mu^{\noc T\cup\{j\}}
\biggr)\big/\Bigl({1\over T_i}\sms{1\over T_j} \Bigr)
\ses
\Phi_\mu^{\noc T\cup\{i,j\}}
\eqno 3.24
$$
and 3.22 follows by induction on $k$ and our definition of the divided difference operator $\Delta$.
\sa

Of course all of these expressions and identities can be viewed solely as
part of the theory of Macdonald polynomials. However, the SF heuristic associates
a meaning to them which suggests manipulations, results and conjectures that are difficult
to predict without it. For instance, the case $T=\emptyset$ of 3.23 is
$$
\biggl(
{1\over T_i}\ess \TH_{\aaa^{(i)}}\sms {1\over T_j}\ess \TH_{\aaa^{(j)}}
\biggr)\big/\Bigl({1\over T_i}\sms{1\over T_j} \Bigr)\ess  
\ses
\Phi_\mu^{\noc \{i,j\}}
\eqno 3.25
$$
and the SF heuristic states that the divided difference on the left-hand
side must be Schur positive since it should give the Frobenius characteristic
of the intersection  of the modules $\BM_{\aaa^{(i)}}$ and $\BM_{\aaa^{(j)}}$.
Of course the further divided differences given by 3.23 should turn out to be
Schur positive as well. Now these particular facts can be verified
by computer much more extensively and for much larger partitions
than the partitions for which we can carry out the verification of the SF heuristic.
Although such experimental verifications give support to the SF heuristic,
more substantial support comes from results (predicted by SF) 
which we can actually prove in full generality from within the theory
of Macdonald polynomials. 
For instance, it can be shown (see [\bergeronhamel]) that for all $\lambda=j 1^{n-j}$
and many other infinite families of partitions, the coefficients of $S_\lambda$ in
$\Phi_\mu^{(k)}$ is (as expected) in $\N[q,t]$. In fact, Bergeron-Hamel derive many explicit
plethystic expressions for these coefficients from those given by Garsia-Tesler in [\grasiaatesler]. 
To emphasize this point we will carry out a calculation
which supports part (v) of SF. This is the assertion that the dimension of
of any $k$ of the modules $\BM_{\aaa^{(i)}}$ is $n!/k$. The verification 
of this property is based on the following.
\sas

\Prop{3.5}
 For any $S=\{1\leq i_1<i_2<\cdots <i_k\leq m\}$ we have 
$$
\Phi_\mu^{\noc S}\ses 
\sum_{\matrix{ i_s\in S}} \ssp 
\biggl(
\prod_{ \matrix {i_r\in S\cr i_r\neq i_s}} \ess {1\over 1- T_{i_s}/T_{i_r}}\ssp 
\biggr)\ess \TH_{\aaa^{(i_s)}}
\eqno 3.26
$$
In particular, under SF, the Hilbert series of the module $\bigcap_{s=1}^k\BM_{\aaa^{(i_s)}}$
is given by the rational expression
$$
F_{\noc S}(q,t)\ses \sum_{\matrix{ i_s\in S}} \ssp \biggl(
\prod_{ \matrix {i_r\in S\cr i_r\neq i_s}} \ess {1\over 1- T_{i_s}/T_{i_r}}\ssp 
\biggr)\ess \del_{p_1}^n\TH_{\aaa^{(i_s)}} 
\eqno 3.27
$$ 

\Proof Substituting in 3.18 $\Phi_\mu$ as given by 3.15 and using 3.12 we obtain
$$
\Phi_\mu^{\noc S}\ses \sum_{i\in S}\ssp 
\biggl(\prod_{\matrix{j=1\cr j\neq i}}^m {1\over 1-T_i/T_j}\ssp \biggr)\ess 
\prod_{r\notin S}\Bigl(1-T_i/T_r\Bigr) \ess \TH_{\aaa^{(i)}}
$$
making the appropriate cancellations gives 3.26 and 3.27 follows by taking the 
$n^{th}$ derivative of both sides of 3.26 with respect to $p_1$.
\sas

Of course, on the validity of SF, the expression on the right-hand side of 3.27 should
evaluate to a polynomial with positive integer coefficients and the dimension of the
module  $\bigcap_{s=1}^k\BM_{\aaa^{(i_s)}}$ could then be simply obtained  by evaluating
$F_{\noc S}(q,t)$ at  $t=q=1$. What we shall show here is that for any $\mu\part n+1$
and $S$ of cardinality $k$ we have
$$
\lim_{q\RA 1}\ssp F_{\noc S}(q,1) \ses  {n!\over k} 
\eqno 3.28
$$
Using one of the specializations computed by Macdonald in [\macdonaldbook] (see [\garsiahaimankostka]) it
can be shown  that for any partition $\bbb=(\bbb_1,\bbb_2,\ldots ,\bbb_l)$ we have
$$
\TH_\bbb(x;q,1)\ses \prod_{i=1}^l\ssp (q)_{\bbb_i}\ssp 
h_{\bbb_i}\bigl[{\textstyle {X\over 1-q}}\bigr]\ess .
\eqno 3.29
$$ 
where for any index $s$ we set $(q)_s=(1-q)(1-q^2)\cdots (1-q^s)$.
Assuming the row that contains the corner cell we remove from $\mu$ to
obtain $\aaa^{(i)}$ has length $a_i$ we deduce from 3.29 that
$$
\TH_{\aaa^{(i)}}(x;q,1)\ses \TH_\ggg(s;q,1)\ssp 
(q)_{a_i-1}h_{a_i-1} \bigl[{\textstyle {X\over 1-q}}\bigr]
\prod_{\matrix{j=1\cr j\neq i}}^m\ssp (q)_{a_j}
h_{a_j} \bigl[{\textstyle {X\over 1-q}}\bigr]
\eqno 3.30
$$ 
where $\ggg$ is the partition obtained by removing from $\mu$ all the rows
that contain corner cells.
\sas

Now it is easily seen that 
$$
T_{\aaa^{(i)}}/T_{\aaa^{(j)}}\big|_{t=1}\ses q^{a_j}/q^{a_i} \ess .
$$
Thus from 3.27 we derive
$$
F_{\noc S}(q,1)\ses \sum_{\matrix{ i_s\in S}} \ssp \biggl(
\prod_{ \matrix {i_r\in S\cr i_r\neq i_s}} \ess { q^{a_{i_s}}\over  q^{a_{i_s}}- q^{a_{i_r}}}\ssp 
\biggr)\ess \del_{p_1}^n\TH_{\aaa^{(i_s)}}(x;q,1)\ess . 
\eqno 3.31
$$
Note that if $a_1+a_2+\cdots +a_m-1=r$ then $|\ggg|=n-r$. 
Let us set for a moment
$$
F_\ggg(q,t)\ses \del_{p_1}^{n-r}\TH_\ggg(x;q,t)\ess .
$$ 
This given, applying $\del_{p_1}^n$ to both sides of 3.30
and using Leibnitz formula we obtain
$$
\eqalign{
\del_{p_1}^n\ssp \TH_{\aaa^{(i)}}(x;q,1) &\ses
{n!\over (n-r)!}\ssp F_\ggg(q,1)\ssp {[a_i-1]_q!\over (a_i-1)!}
\prod_{\matrix{j=1\cr j\neq i }}^m\ssp {[a_j]_q!\over a_j!} \cr
&\ses {n!\over (n-r)!} \ssp  F_\ggg(q,1) 
\ssp \Bigl(\prod_{j=1 }^m\ssp {[a_j]_q! \over a_j!}\ssp \Bigr)\ssp { (a_i)\over [a_i]_q }
\ess .
\cr
}
\eqno 3.32
$$
Where as customary, for any integer $l$ we set $[l]_q=(1-q^l)/(1-q)$ and
$[l]_q!=(q)_l/(1-q)^l$.

Setting $i=i_s$ in 3.32 and substituting in 3.31  we may write
$$
{
F_{\noc S}(q,t)\over
{n!\over (n-r)!} \ssp  F_\ggg(q,1) 
\ssp \Bigl(\prod_{j=1 }^m\ssp {[a_j]_q!\over a_j!}\ssp \Bigr)
}
\ses \sum_{\matrix{ i_s\in S}} \ssp \biggl(
\prod_{ \matrix {i_r\in S\cr i_r\neq i_s}} \ess { q^{a_{i_s}}\over  q^{a_{i_s}}- q^{a_{i_r}}}\ssp 
\biggr)\ess { (a_{i_s})  \over [a_{i_s}]_q}\ess . 
$$
Now we know (see [\garsiahaimankostka]) that $\lim_{q\RA 1}\ssp F_\ggg(q,1)=(n-r)! $. Moreover it is easy to
see that all the ratios $[a_j]_q!/a_j!$ tend to $1$ as $q\RA  1$. Thus we see that
the identity in 3.28 is a simple consequence of the following elementary fact:

\Lem{3.1}
For any distinct integers $y_1,y_2,\ldots ,y_k$ we have
$$
\lim_{q\RA 1}\ess (1-q)\ssp \sum_{s=1}^k\ssp 
\biggl(\prod_{1\leq r\leq k\atop  r\neq s} {q^{y_s}\over q^{y_s}-q^{y_r} }\biggr)
\ssp {y_s \over1-q^{y_s}} \ses {1\over k}
\eqno 3.33
$$

\Proof We give a sketch of the argument since the details are easily filled.
Our point of departure is the partial fraction decomposition
$$
\prod_{s=1}^k\ssp {1\over 1-t\ssp q^{y_s}} \ses
\sum_{s=1}^k\ssp \biggl(\prod_{1\leq r\leq k\atop  r\neq s} {q^{y_s}\over q^{y_s}-q^{y_r} }\biggr)
\ssp {1 \over1-t\ssp q^{y_s}}\ess .
\eqno 3.34
$$
Now the trick in proving 3.33 is to use the expansion
$$
{y_s\over 1-q^{y_s}}\ses {1\over  \log_e\ssp 1/q}\ssp \sum_{p\geq 1} {1\over p}\ssp  
(1-q^{y_s})^{p-1}
\eqno 3.35
$$
and noting that, since the term multiplying $y_s/(1-q^{y_s})$ in 3.33
has only a pole of order $k-1$
at $q=1$, in calculating our limit we need only the first $k$ terms of the series in 3.35.

Since 3.34 gives that for any integer $a\geq 0$
$$
\sum_{s=1}^k\ssp \biggl(\prod_{1\leq r\leq k\atop  r\neq s} {q^{y_s}\over q^{y_s}-q^{y_r} }\biggr)
\ssp q^{a y_s}\ses h_a[q^{y_1},q^{y_2},\ldots ,q^{y_k}]\ssp \RA 
\ssp \Bigl(\ssp \matrix{ a+k-1 \cr a }\ssp \Bigr)
\ess\ess\ess\ess (\ssp {\rm as} \ess\ess q\RA 1\ssp )
\ess ,
$$ 
we see that 3.33 is reduced to showing that
$$
\sum_{p=1}^k\ssp{1\over p}\ssp \sum_{a=0}^{p-1}\ssp (-1)^a\ssp 
\Bigl(\ssp \matrix{ p-1 \cr a }\ssp \Bigr)\Bigl(\ssp \matrix{ a+k-1 \cr a }\ssp \Bigr)
\ses {1\over k}\ess . 
$$
However, this turns out to be an immediate consequence of the classical 
Gauss summation formula.
\sas

Our proof of 3.28 is thus complete and we have verified that property
$(v)$ of the SF heuristic is consistent with Macdonald theory. 
\sas

We terminate this section with three curious identities which, on the validity
of SF, imply further Schur positivity results.
\sas

\Thm{3.2}
$$
\CF\ssp \bigvee_{i=1}^m\ssp \BM_{\aaa^{(i)}}
\ses \sum_{k=1}^m \Phi_\mu^{(k)}   \ssp 
e_{m-k}\Bigl[{1\over T_1}+{1\over T_2}+\cdots +{1\over T_m}\Bigr]
\ses {\nabla \ssp \Phi_\mu^{(1)}\over T_1T_2\cdots T_m}\ess ,
\eqno 3.36
$$
$$
\sum_{k=1}^m \P k \ssp  e_{m+1-k}\Bigl[{1\over T_1}+{1\over T_2}+\cdots +{1\over T_m}\Bigr]
\ses \nabla^{-1} \Phi_\mu\ess .
\eqno 3.37
$$

\Proof We clearly have the decomposition
$$
\bigvee_{i=1}^m\ssp \BM_{\aaa^{(i)}}\ses 
\bigoplus_{\matrix{S\con\{1,2,..,m\}\cr S\neq \emptyset}}
\BM^{\eee(S)}\ess .
$$
Thus
$$
\eqalign{
\CF\ssp \bigvee_{i=1}^m\ssp \BM_{\aaa^{(i)}}
&\ses \sum_{\matrix{S\con\{1,2,..,m\}\cr S\neq \emptyset}}
\ssp \Phi_\mu^{=S}\cr
&\ses\sum_{k=1}^m\ssp \P k 
\sum_{\matrix{S\con\{1,2,..,m\}\cr{|S|=k}\cr}}
\ssp \prod_{i\notin S}{1\over T_i}\ess .
\cr
}
$$
This shows the first equality in 3.36. For the second equality we note that we can write
$$
\eqalign{
\sum_{k=1}^m \P k \ssp e_{m-k}\Bigl[{1\over T_1}+{1\over T_2}+\cdots +{1\over T_m}\Bigr]
&\ses
\sum_{k=1}^m (-\nabla)^{m-k}\ssp \Phi_\mu \ssp  
\ssp e_{m-k}\Bigl[{1\over T_1}+{1\over T_2}+\cdots +{1\over T_m}\Bigr]\cr
&\ses
\biggl\{\ssp
\prod_{i=1}^m\ssp \Bigl(1-{\nabla\over T_i}\Bigr)
\sms (-1)^m\ssp {\nabla^m\over T_1T_2\cdots T_m }
\biggr\}\ess \Phi_\mu\cr
&\ses
(-1)^{m-1}\ssp {\nabla^m\over T_1T_2\cdots T_m }\ssp \Phi_\mu
\cr
}
$$
and the result follows from 3.17 for $k=1$.

This given, applying $\nabla^{-1}$ to the second equality in 3.36, and
using again the equalities in 3.17 we get:
$$
\sms \sum_{k=1}^{m-1} \Phi_\mu^{(k+1)}   \ssp 
e_{m-k}\Bigl[{1\over T_1}+{1\over T_2}+\cdots +{1\over T_m}\Bigr]
\sps \nabla^{-1}\ssp \Phi_\mu
\ses { \Phi_\mu^{(1)}\over T_1T_2\cdots T_m}\ess ,
$$
and 3.37 follows by rearranging terms and changing the index of summation. 
\sas

To state our next result we need further ingredients.
Let us recall that in 3.13 we have
set $x_i=t^{l_i}q^{a_i}$ with $l_i$ and $a_i$  denoting
the coleg and coarm of the cell we must remove from $\mu$ to  get $\aaa^{(i)}$. 
 We may refer to $x_i$ as the ``weight'' of the  $i^{th}$
``outer'' corner cell of $\mu$. We shall need here also the monomials 
$$
u_i\ses t^{l_{i+1}}q^{a_i}\bigsp\bigsp(\ssp i=1,2,\ldots,m-1\ssp )
\eqno 3.38
$$
giving the weights of the ``inner'' corner cells. We shall also set
$$
u_0\ses t^{l_1}/q\ess\scs\ess 
u_m\ses q^{a_m}/t\ess\scs\ess 
x_0\ses 1/tq
\eqno 3.39
$$

% Adriano,
%	Here's the illustration above Proposition 2.3, redone in TeX macros
%								Glenn

% coordinate system is ordinary Euclidean, (0,0) in lower left corner,
% x increases to right, y increases upwards
% coordinates are expressed in units \mylength

\newdimen\mylength
%\mylength=1pt
%\def\cornd{6}
\mylength=.3pt		% unit coordinate
\def\cornd{18}		% side length of corner squares in above units

\def\corndl{\cornd\mylength}

% place #3 at coordinates #1,#2
% #1,#2
\def\myput(#1,#2)#3{\raise#2\mylength\hbox to 0pt{\hskip #1\mylength #3\hss}}

% \hseg(#1,#2;#3) draws a line from (#1,#2) to (#1+#3,#2)
\def\hseg(#1,#2;#3){\myput(#1,#2){\vrule width#3\mylength depth.4pt}}

% \vseg(#1,#2;#3) draws a line from (#1,#2) to (#1,#2+#3)
\def\vseg(#1,#2;#3){\myput(#1,#2){\vrule height#3\mylength width.4pt depth.4pt}}

% \corn(#1,#2) draws a corner cell with given upper right coordinate
\def\corn(#1,#2){\myput(#1,#2){\vbox to0pt{\hbox to0pt{\hss%
	\vbox{\hrule width \corndl\hbox to\corndl{%
		\vrule height \corndl \hfil \vrule height \corndl}%
		\hrule width \corndl}%
	\kern-.4pt}\vss}}}
% \labUR(#1,#2;#3), \labLL(#1,#2;#3)
% put label #3 beyond upper right (\labUR) or lower left (\labLL) vertex of a
%  \corn(#1,#2) (but don't draw the \corn)
\def\labUR(#1,#2;#3){\myput(#1,#2){\hskip 2pt$#3$}}
\def\labLL(#1,#2;#3){\myput(#1,#2){\myput(-\cornd,-\cornd){%
	\vbox to0pt{\hbox to0pt{\hss $#3$\hskip 1pt}\vss}}}}

% \drawL(x1,x2,y1,y2;A,B)  draws

%              A
% y2  +-------+
%    O       O|
%   B         |
%             |
% y1          +
%     x1      x2
% where O is corner square, A, B are labels
\newcount\drawcalc
\def\drawL(#1,#2,#3,#4;#5,#6){%
%	\vseg(x2,y1;y2-y1)\corn(x2,y2)\labUR(x2,y2;A)%
	\drawcalc=-#3 \advance\drawcalc#4%
	\vseg(#2,#3;\drawcalc)\corn(#2,#4)\labUR(#2,#4;#5)%
%
%	\hseg(x1,y2;x2-x1)\corn(x1,y2)\labLL(x1,y2;B)%
	\drawcalc=-#1 \advance\drawcalc#2%
	\hseg(#1,#4;\drawcalc)\corn(#1,#4)\labLL(#1,#4;#6)%
}

To appreciate the geometric significance of these weights, in the figure below we illustrate
a $4$-corner case with 
outer corner cells labelled $A_1$, $A_2$, $A_3$, $A_4$   
and inner corner cells
labelled  $B_1$, $B_2$, $B_3$, $B_4$.
%\hskip 2truein
$$\hbox{%
	\hseg(0,0;300)\corn(300,0)\labLL(300,0;B_4)% bottom line
	\drawL(220,300,0,60;A_4,B_3)%
	\drawL(130,220,60,130;A_3,B_2)%
	\drawL(70,130,130,200;A_2,B_1)%
	\drawL(0,70,200,260;A_1,B_0)%
	\vseg(0,0;260)% left line
}
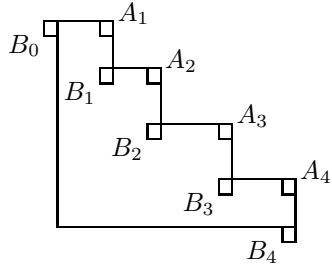
$$
\sa
\sas
%$$
%\cells
%$$
 
We recall (see [\garsiahaimanorbit]) that one of the Pieri rules given by Macdonald  
may be written in the form
$$ 
\del_{p_1} \TH_\mu\ses \sum_{\nu\RA \mu}\ssp c_{\mu\nu}(q,t)\ssp \TH_\nu
$$
where $\del_{p_1}$ denotes the Hall inner product adjoint of multiplication
by $p_1$ and $\nu\RA\mu$ is to indicate that the sum is to be carried out
over the partitions $\nu$ which immediately precede $\mu$ in the Young order.
In our present notation we can write this identity in the form
$$ 
\del_{p_1} \TH_\mu\ses \sum_{i=1}^m\ssp c_{\mu\aaa^{(i)}}(q,t)\ssp \TH_{\aaa^{(i)}}
\eqno 3.40
$$
It is well known and it is easy to show that if $\Psi$ is the Frobenius
characteristic of an $S_n$ module $\BM$ then $\del_{p_1}\Psi$ gives the Frobenius characteristic
of the  module $\BM$ {\ita restricted  to} $S_{n-1}$. Thus,
under the $C=\TH$ conjecture the polynomial $\del_{p_1} \TH_\mu$ should
give the Frobenius characteristic of $\BM_\mu$ as an $S_{n-1}$-module.
\sas

Our next and final task in this section is to compute the expansion 
of this characteristic in terms of our polynomials $\Phi_\mu^{(k)}$.
This computation became possible only recently due to the discovery in [\grasiaatesler] that 
the Macdonald Pieri coefficients $c_{\mu\nu}$ reduce to a remarkably simple
form when expressed in terms of the weights of the ``outer'' and ``inner''
corner cells of $\mu$. This result may be stated as follows.

\Prop{3.6}
For any partition $\mu$ with predecessors $\aaa^{(1)},\aaa^{(2)},\ldots,\aaa^{(m)}$
and corner weights $x_0,x_1,\ldots,x_m$, $u_0,u_1,\ldots,u_m$ given by 3.13 , 3.38
and 3.39 we have 
$$
\del_{p_1}\ssp \TH_\mu\ses {1\over M}\ssp 
\sum_{i=1}^m{1\over x_i}\ssp 
{
\prod_{s=0}^m (u_s-x_i)
\over 
\prod_{s=1,s\neq i}^m (x_s-x_i)
}\ess \TH_{\aaa^{(i)}} \ess ,
\eqno 3.41
$$
where for convenience we have set
$$
M=(1-1/t)(1-1/q)\ess .
\eqno 3.42
$$

\noindent The proof is given in [\grasiaatesler] (see Proposition 2.4 there).

\sas

Combining this result with our expansions of $\TH_{\aaa^{(i)}}$ we obtain
the following remarkable identities.

\Thm{3.3}
Under the same hypotheses  
$$
\del_{p_1}  \TH_\mu\ses {1\over M}\ssp 
{T_\mu\over \nabla}\ssp
\biggl\{\ssp 
\prod_{s=0}^m\Bigl(1-{\nabla \over T_\mu}\ssp u_s\Bigr)
-
\prod_{s=0}^m\Bigl(1-{\nabla \over T_\mu}\ssp x_s\Bigr)
\biggr\}\ssp \Phi_\mu\ess ,
\eqno 3.43
$$ 
 or equivalently:
$$
\del_{p_1} \TH_\mu\ses 
\sum_{k=1}^m\ssp{\Phi_\mu^{(k)}\over T_\mu^{m-k}}\ssp 
{e_{m+1-k}[x_0+\cdots+x_m]-e_{m+1-k}[u_0+\cdots+u_m]\over {\scriptstyle {(1-1/t)(1-1/q)}}}
\eqno 3.44
$$

\Proof It is easy to verify from the definitions 3.13 , 3.38 and 3.39 that we have
$$
x_0x_1\cdots x_m\ses u_0u_1\cdots u_m\ess .
\eqno 3.45
$$
Thus the expression
$$
{1\over z}\ssp \biggl\{\ssp 
\prod_{s=0}^m\Bigl(1-z\ssp u_s\Bigr)
-
\prod_{s=0}^m\Bigl(1-z\ssp x_s\Bigr)
\biggr\}
$$
defines a polynomial in $z$ of degree $m-1$. Using the Lagrange interpolation formula
at $z=x_1,x_2,\ldots ,x_m$ we then get that
$$
{1\over z}\ssp \biggl\{\ssp 
\prod_{s=0}^m\Bigl(1-z\ssp u_s\Bigr)
-
\prod_{s=0}^m\Bigl(1-z\ssp x_s\Bigr)
\biggr\}
\ses
\sum_{i=1}^m\ssp x_i
\ssp \prod_{s=0}^m\Bigl(1- \ssp u_s/x_i\Bigr)\ssp
{
\prod_{s=1, \neq i}^m\Bigl(1-z\ssp x_s\Bigr)
\over
\prod_{s=1, \neq i}^m\Bigl(1-  x_s/x_i\Bigr)
}
$$
and this may also  be rewritten as 
$$
{1\over z}\ssp \biggl\{\ssp 
\prod_{s=0}^m\Bigl(1-z\ssp u_s\Bigr)
-
\prod_{s=0}^m\Bigl(1-z\ssp x_s\Bigr)
\biggr\}
\ses
\sum_{i=1}^m\ssp {1\over x_i}
\ssp 
{
\prod_{s=0}^m (x_i- \ssp u_s  )
\over
\prod_{s=1, \neq i}^m (x_i-  x_s )
}
\ssp
\prod_{s=1, \neq i}^m\Bigl(1-z\ssp x_s\Bigr)\ess .
\eqno 3.46
$$
Setting $z={\nabla  \over T_\mu}$ and applying both sides to  $\Phi_\mu$ gives
$$
{T_\mu\over \nabla} 
\biggl\{\ssp 
\prod_{s=0}^m\Bigl(1-{\nabla \over T_\mu}\ssp u_s\Bigr)
-
\prod_{s=0}^m\Bigl(1-{\nabla \over T_\mu}\ssp x_s\Bigr)
\biggr\}\ssp \Phi_\mu
\ses
\sum_{i=1}^m\ssp {1\over x_i}
\ssp 
{
\prod_{s=0}^m (x_i- \ssp u_s  )
\over
\prod_{s=1, \neq i}^m (x_i-  x_s )
}
\ssp
\prod_{s=1, \neq i}^m\Bigl(1-{\nabla \over T_\mu}\ssp x_s\Bigr)
 \Phi_\mu\ess .
\eqno 3.47
$$ 
Since $x_s/T_\mu=1/T_s$ for $s=1,2,\ldots,m$ we see from 3.19 a) that
$$
\prod_{s=1, \neq i}^m\Bigl(1-{\nabla \over T_\mu}\ssp x_s\Bigr)
\Phi_\mu\ses \TH_{\aaa^{(i)}}\ess .
$$
Thus 3.42 follows by combining 3.41 and 3.47.

This given, note that expanding the products on the right-hand side of 3.43 we
get
$$
\del_{p_1}\TH_{\aaa^{(i)}}\ses{T_\mu\over \nabla}\ssp
\sum_{k=1}^m\Bigl(-{\nabla\over T_\mu}\Bigr)^{m+1-k}\Phi_\mu
\ess
{e_{m+1-k}[u_0+\cdots+u_m]-e_{m+1-k}[x_0+\cdots+x_m]\over {\scriptstyle {(1-1/t)(1-1/q)}}}
\ess ,
$$
and 3.43 follows from the identities in 3.17.

\Rem{3.1}
Routine manipulations starting from the definition in I.9 yield that
\footnote{$\dag$}{See Proposition 2.3 in [\grasiaatesler]}
$$
(1-1/t)(1-1/q)\ssp B_\mu(q,t)\ses x_0+x_1+\cdots+x_m\sms u_0-u_1-\cdots-u_m\ess . 
\eqno 3.48
$$
This given, we see from 3.44 that the coefficient of $\Phi_\mu^{(m)}=\Phi_\mu$ 
in $\del_{p_1}\TH_\mu$ is precisely by $B_\mu(q,t)$.

A simple calculation based on Property (iii) of SF (see eq. I.30) yield that
for any $k=1,2,..,m$ we have
$$
\da\ssp \Phi_\mu^{(k)}\ses {\Phi_\mu^{(m+1-k)}\over T_1T_2\cdots T_m} 
\eqno 3.49
$$
Note also that replacing $t,q$ by $1/t,1/q$ in 3.48, and using 3.45, we derive that
$$
\eqalign
{
(1-t)(1-q)\ssp B_\mu(1/q,1/t)&\ses 
{1\over x_0}+{1\over x_1}+\cdots+{1\over x_m}\sms {1\over u_0}-{1\over u_1}-\cdots-{1\over u_m}
\cr
&\ses
{e_m[x_0+\cdots+x_m]\sms e_m[u_0+\cdots+u_m]\over x_0 x_1 \cdots x_m}
\cr
}
$$
Since $x_0=1/tq$ and $x_i=T_\mu/T_i$ we deduce that
$$
{1\over T_\mu^{m-1}}\ssp{e_m[x_0+\cdots+x_m]\sms e_m[u_0+\cdots+u_m]\over(1-1/t)(1-1/q)}
\ses T_\mu\ssp B_\mu(1/q,1/t)\ssp {1\over  T_1T_2\cdots T_m}\ess .
$$ 
Thus using 3.49 we see that we may write 3.43 in the form
$$
\del_{p_1}\TH_\mu\ses B_\mu(q,t)\ssp \Phi_\mu
\sps \cdots \sps
T_\mu \da  B_\mu(q,t) \ssp \Phi_\mu\ess .
\eqno 3.50
$$
Note that in the two-corner case there are no further terms and under the $C=\TH$-conjecture
this identity yields
$$
\del_{p_1}C_\mu (x;q,t)\ses B_\mu(q,t)\ssp \Phi_\mu\sps 
T_\mu \da  B_\mu(q,t) \ssp \Phi_\mu\ess .
\eqno 3.51
$$
which is formula 1.35. Thus we obtain here the general case of an identity
we first encountered with $\mu=32$ (see 1.31).
\sas

We also see that the coefficients of $\P 1 $ and $\P m $ in 3.43
are polynomials in $q,t,1/q,1/t$ with positive integer coefficients. 
This is also true for the coefficients of the remaining $\P k $.
The basic result here may be stated as follows.

\Thm{3.4} 
If $x_0,x_1,\ldots , x_m$ and $u_0,u_1,\ldots , u_m$ are the corner weights 
of a partition then for any $k=1,2,\ldots ,m$ the expression
$$
B_\mu^{(k)}(q,t)\ses 
{
e_{k}[x_0+\cdots+x_m]-e_{k}[u_0+\cdots+u_m]
\over 
{\scriptstyle {(1-1/t)(1-1/q)}}
}
$$
evaluates to a polynomial in $q,t$ with positive integer coefficients.
\sas

\Proof We owe the following very simple induction argument to Glen Tesler.
Note first that for $m=1$ we have
$$
x_0=1/qt\scs\ess  x_1=t^{l_1}q^{a_1}\ess ;\ess\ess   u_0=t^{l_1}/q\scs\ess u_1=q^{a_1}/t\ess .
$$
Thus
$$
B_\mu^{(1)}(q,t)\ses 
{
1/qt + t^{l_1}q^{a_1}-t^{l_1}/q-q^{a_1}/t
\over
{ {(1-1/t)(1-1/q)}}
}\ses [l_1+1]_t\ssp [a_1+1]_q\ess .
$$
Let us then assume the result true for all partitions with $m-1$ corners. 
For a partition $\mu$ with $m$ corners and weights given by 3.13, 3.38 and 3.39
we write (see 3.42)
$$
 B_\mu^{(k)}=
{
e_{k}[x_0+X_{m-1}+x_m] - e_{k}[\tx_0+X_{m-1}+u_m]
\over M}
+ 
{e_{k}[\tx_0+X_{m-1}+u_m] - e_{k}[u_0+U_{m-1}+u_m]
\over M} 
$$
where for convenience we have set
$$
\tx_0= t^{l_m}/q\ess \ess ;  \ess\ess 
X_{m-1}=x_1+x_2+\cdots +x_{m-1}\scs \ess U_{m-1}=u_1+u_2+\cdots +u_{m-1}\ess .
\eqno 3.52
$$
Now we have
$$
\eqalign{
{e_{k}[x_0+X_{m-1}+x_m] - e_{k}[\tx_0+X_{m-1}+u_m]\over M}
&\ses \sum_{s=1}^2\ssp {  e_{s}[x_0+x_m]- e_{s}[\tx_0+u_m] \over M}
 \ssp
e_{k-s}[X_{m-1}]\cr
&\ses
{x_0+x_m - \tx_0-u_m  \over M }\ess e_{k-1}[X_{m-1}]
\cr
&\ses [l_m+1]_t[a_m+1]_q \ess  e_{k-1}[X_{m-1}]\ess ,
\cr
}  
$$
where the second equality is due to the relation $x_0x_m=\tx_0u_m$.
\sas

Similarly we deduce that
$$
\eqalign{
&
{e_{k}\bigl[\tx_0+X_{m-1}+u_m\bigr] - e_{k}\bigl[u_0+U_{m-1}+u_m\bigr]\over M}
\cr
&\hskip14.5cm{\rm 3.53}
\cr
&\ess\ess\ess\ess\ess\ess\ess\ess\ess\ess\ess\ess\ess
=
{e_{k}\bigl[\tx_0+X_{m-1}\bigr]-e_{k}[u_0+U_{m-1}\bigr]\over M}+  
{e_{k-1}\bigl[\tx_0+X_{m-1}\bigr]-e_{k-1}[u_0+U_{m-1}\bigr]\over M}\ess u_m
\cr
}
$$
Note that the auxiliary monomial $\tx_0$ we introduced in 3.52 may be viewed as the weight
of the cell $\tilde A_0$ that is immediately to the left of the diagram of $\mu$ and is in the same row
as the cell $A_m$ with weight $x_m$. Let $\tilde \mu$ denote the $m-1$-corner partition obtained by removing
from $\mu$ all the rows on or below the cells $\tilde A_0$ and $A_m$.
After drawing the diagrams of $\mu$ and $\tmu$, a minute of reflection should reveal that the monomials 
$$
{\tx_0\over t^{l_m}}\scs 
{x_1\over t^{l_m}}\scs \ldots ,
{ x_{m-1}\over t^{l_m}}
\ess\ess ;\ess\ess\ess
{u_0\over t^{l_m}}\scs 
{u_1\over t^{l_m}}\scs \ldots ,
{ u_{m-1}\over t^{l_m}}
$$  
are none other than the corner weights of the partition $\tilde \mu$. 

This given, we can rewrite 3.53 in the form
$$
{e_{k}\bigl[\tx_0+X_{m-1}+u_m\bigr] - e_{k}\bigl[u_0+U_{m-1}+u_m\bigr]\over M}
\ses t^{kl_m}
\ess B_\tmu^{(k)}(q,t)\sps
 u_m\ssp
t^{(k-1)l_m}
B_\tmu^{(k-1)}(q,t) \ess .
$$
In summary we have shown that
$$
B_\mu^{(k)}(q,t)\ses [l_m+1]_t[a_m+1]_q\ssp e_{k-1}[X_{m-1}]
\sps
t^{kl_m}
\ess B_\tmu^{(k)}(q,t)\sps
 u_m\ssp
t^{(k-1)l_m}
B_\tmu^{(k-1)}(q,t) \ess .
$$
which proves the positive integrality of $B_\mu^{(k)}(q,t)$ and completes our
induction argument. 
\sas

We shall see in the next section that the positivity of these coefficients expresses a remarkable combinatorial
and representation theoretical process intimately connected with the SF heuristic..  
\sa

\section{4. Dissecting $\BM_\mu$ as an $S_{n-1}$-module}
We have seen in section 3 that, at least in the case of the partition $\mu=321$, it is possible to use
the bases $\CB^{\eee_1\eee_2\eee_3}$ to construct a basis which decomposes $\BM_{321}$ into a direct sum
of six left regular representations of $S_5$. In summary for each cell $(i,j)\in 321$ we constructed a 
basis $\CB_{ij}$
as a union of bases $\CB^{\eee_1\eee_2\eee_3}$ then the basis for $\BM_{321}$ was taken to be of the form
$$
\CB_{321}\ses \bigcup_{(i,j)\in 321}  \CB_{ij}(\del)\ssp \del_{x_6}^i\del_{y_6}^j\Delta_{321}\ess .
\eqno 4.1 
$$
Note that our construction yields that the space $\BM_{ij}$ spanned by the basis $\CB_{ij}$
is a direct sum of the $S_5$-modules $\BM^{\eee_1\eee_2\eee_3}$ and thus is itself
an $S_5$-module. Using this fact, it is not difficult to derive that our construction
also yields the direct sum decomposition
$$
\BM_{321}\ssp \big|_{s_5}\ses \bigoplus_{(i,j)\in 321}\ssp {\bf flip}_{\Delta_{321}(ij)}\ess \BM_{ij}\ess ,
\eqno 4.2 
$$
where for convenience we have set
$$
\Delta_{321}(ij)\ses \del_{x_6}^i\del_{y_6}^j\Delta_{321}\ess .
\eqno 4.3 
$$
To be sure, the polynomial $\Delta_{321}(ij)$ is not one of our $\Delta_\mu$.  Nevertheless it is an $S_5$
alternant and this is all that is needed for the corresponding $\bf flip$ map to have all of the 
properties we need. In particular, letting $\Phi_{ij}$ denote the bigraded Frobenius characteristic of 
$\BM_{ij}$, we can deduce from 4.2 that the bigraded Frobenius characteristic
$C_{321}(x;q,t)$ satisfies the equation
$$
\del_{p_1}C_{321}(x;q,t)\ses \sum_{(i,j)\in 321}\ssp {T_{321}\over t^iq^j}\ssp \DA\Phi_{ij}\ess .
\eqno 4.4 
$$
Since for any partition $\mu$ we have $T_\mu \DA C_\mu(x,q,t)= C_\mu(x,q,t)$, applying $\DA$ to
both sides of 4.4 and multiplying by $T_{321}$ we derive that
$$
\del_{p_1}C_{321}(x;q,t)\ses \sum_{(i,j)\in 321}\ssp \Phi_{ij}\ssp t^iq^j  \ess .
\eqno 4.5 
$$
Looking back at the display of bases $\CB_{ij}$ given in section 2, and using the notation
we introduced in section 3, we can easily deduce that in this case
$$
\matrix{
\Phi_{00}\ses \Phi^{(3)}\sps {1\over T_{321}}\ssp \bigl( x_1+x_2+x_3\bigr)\Phi^{(2)}
 \sps {1\over T_{321}^2}\ssp \bigl( x_1x_2+x_1x_3+x_2x_3\bigr)\Phi^{(1)}\cr
\Phi_{01}\ses \Phi^{(3)}\sps {1\over T_{321}}\ssp \bigl( x_1+x_2+x_3\bigr)\Phi^{(2)}\hfill\cr
\Phi_{10}\ses \Phi^{(3)}\sps {1\over T_{321}}\ssp \bigl( x_1+x_2+x_3\bigr)\Phi^{(2)}
\sps {1\over T_{321}^2}\ssp \bigl(x_1x_3+x_2x_3\bigr)\Phi^{(1)}\hfill\cr
\Phi_{20}\ses \Phi^{(3)}\sps {1\over T_{321}}\ssp \bigl( x_2+x_3\bigr)\Phi^{(2)}
 \sps {1\over T_{321}^2}\ssp   x_2x_3\,\Phi^{(1)} \hfill\cr
\Phi_{11}\ses \Phi^{(3)}\sps {1\over T_{321}}\ssp x_3\,\Phi^{(2)} \hfill\cr
\Phi_{02}\ses \Phi^{(3)}\hfill\cr
}
\eqno 4.6
$$
To test the consistency of all our conjectures we should want to verify that substituting these 
expressions for the $\Phi_{ij}$ in 4.5 we do obtain the same final expression for $\del_{p_1}C_{321}(x;q,t)$
that can be computed by means of 3.44. Taking account that in this case we have  
$$
T_{321}=t^3q^3\ssp ; \ess\ess
\matrix{
 x_o= 1/tq\ssp , &x_1=t^2\ssp , & x_2=tq\ssp , & x_3=q^2\ssp ,\cr 
 u_o= t^2/q\ssp , &u_1=t\ssp , & u_2=q\ssp , & u_3=q^2/t\ssp ,\cr 
}
$$
it can be somewhat tediously verified that, in fact, 4.5 and 3.44 do turn out to be in complete agreement.
We shall not do so here 
since, very shortly,  we will prove a general result that includes this verification as a particular case.
%The assignment of modules $\BM_{ij}$ that yields the Frobenius characteristics $\Phi_{ij}$ can best
%be visualized by representing each $\BM_{ij}$ by the  appropriate union of pieces of the Venn diagram 
%for $\mu=321$ and placing the corresponding figure in the cell $(i,j)$. This results in the diagram
%given below.
%$$
%\diagram
%$$ 
 
Upon close examination of the process that led to these choices of the modules $\BM_{ij}$ for
the partition $\mu={321}$, Mark Haiman and Nantel Bergeron put together a scheme for constructing
an assignment of modules $\BM_{ij}$ to the cells any Ferrer's diagram.  They conjectured
that if $\mu$ is an $m$-corner partition of $n$, and the corresponding bases $\CB_{ij}$ are made
up of unions of the bases $\CB^{\eee_1\eee_2\cdots \eee_m}$,
then the collection
$$
\CB_\mu\ses \bigcup_{(i,j)\in \mu}  \CB_{ij}(\del)\ssp \del_{x_6}^i\del_{y_6}^j\Delta_\mu\ess .
\eqno 4.7 
$$
should turn out to be a basis for $\BM_\mu$. Moreover, if $\Phi_{ij}$ denotes 
the Frobenius characteristic of $\BM_{ij}$, then we should also have that
$$
\del_{p_1}C_{\mu}(x;q,t)\ses \sum_{(i,j)\in\mu}\ssp \Phi_{ij}\ssp t^iq^j  \ess .
\eqno 4.8 
$$
This is of course in complete analogy with 4.1 and 4.4. Although  we cannot prove any of this
at the present time, we are nevertheless in a position to show that their construction leads
to a formula for $\del_{p_1}C_{\mu}(x;q,t)$ which is in full generality consistent 
with Macdonald theory and the $C=\TH$ conjecture. To describe the Bergeron-Haiman Algorithm,
which here and after we shall simply refer as  refer as the $BH$ algorithm, we shall use the same
notation we introduced in section 3. The algorithm proceeds
one row at the time starting from the top row.  The first step is to assign to each cell
of the top row the module $\BM_1$. Inductively, the assignements for each row are  
obtained from the assignments for the previous row by exactly the same proceedure.
There are however two cases. Calling $A$ and $B$ the two successive rows, 
if they have the same length then the assignments for $B$ are the same as those for $A$. 
If $A$ has length $a$ and $B$ has length $b$ and $c=b-a\geq 1$ then let 
$$
\BA_1\scs \BA_2,\scs \ldots \scs \BA_a
$$  
be the assignments for the successive cells of row $A$ (from left to right) and set
$$
\BA_s'\ses 
\cases { 
\BA_s & for $1\leq s\leq a\ssp ,$ \cr
\vee _{i=1}^m \BM_i    & for $s\leq 0\ssp ,$ \cr
\{0\} & for $a+1\leq s\leq b\ssp .$
}
\eqno 4.9
$$
Then, the $s^{th}$ cell (from left to right) of row $B$ gets assigned the module 
$$
\BB_s\ses \BA_s'\ssp \vee\bigl( \BA_{s-c}'\cap \BC\ssp \big)
\eqno 4.10
$$ 
where $\BC=\BM_{(i)}$ if the cell at end of row $B$ is the $i^{th}$ corner cell of $\mu$.  
It is important to note that if 
$$
\BA_1\noc \BA_2\noc\cdots \noc \BA_a
$$
then 4.9 gives that 
$$
\BA_{1-c}'\noc \BA_{2-c}'\noc\cdots \noc \BA_{b-c}'
\eqno 4.11
$$
and 4.10 yields that  
$$
\BB_1\noc \BB_2\noc\cdots \noc \BB_b
$$
Since we start by assigning the same module to all the cells of the first row, we see that
this process will assign to each successive row a decreasing sequence of $S_{n-1}$-modules
each of which may be decomposed into a direct sum of the modules $\BM^{\eee_1\eee_2\cdots \eee_m}$.
\sas

Now, according to formula 4.5, if $A$ is the $j^{th}$ row of $\mu$ , 
the contribution to $\del_{p_1}\TH_\mu(x;q,t)$ coming from the Frobenius characteristics
assigned to the cells of row $A$ will be $t^{j-1}$ times the polynomial
$$
L_A\ses \CF \BA_1\sps q\ssp \CF \BA_2\sps \cdots \sps q^{a-1}\CF \BA_a
\eqno 4.12
$$ 
while the contribution of row $B$ will be $t^j$ times the polynomial
$$
L_B\ses \CF \BB_1\sps q\ssp \CF \BB_2\sps \cdots \sps q^{b-1}\CF \BB_b\ess .
\eqno 4.13
$$
We need to know how to compute $L_B$ starting from $L_A$. Note that if $\BM_1$ and $\BM_2$ are two
$S_n$-submodules of a given $S_n$-module $\BM$ then we necessarily have 
$$
\CF\ssp \bigl( \BM_1\vee\BM_2\bigr)\ses \CF\ssp\BM_1\sps \CF\ssp\BM_2 
\sms \CF\ssp\bigl( \BM_1\wedge\BM_2\bigr)
$$
This given, taking account of the inclusion $A'_{s-c}\noc A'_s$, we deduce from 4.10 that
$$
\CF \BB_s\ses \CF  \BA'_s\sps \CF\ssp \bigl(\BA'_{s-c} \cap \BC\ssp \bigr)
\sms \CF\ssp \bigl(\BA'_s \cap \BC\ssp \bigr)
$$
Summing from $1$ to $b$ and taking account of 4.9 we get
$$
\eqalign{
L_B &\ses \sum_{s=1}^aq^{s-1} \CF  \BA_s\sps \sum_{s=1}^cq^{s-1} \CF\ssp  \BC\sps 
\sum_{s=c+1}^bq^{s-1} \bigl(\BA_{s-c} \cap \BC\ssp \bigr)\sms \sum_{s=1}^aq^{s-1} \CF  \bigl(\BA_s\cap \BC\ssp \bigr)
\cr 
&\ses
L_A\sps {q^c-1\over q-1 }\ssp \CF\ssp \BC\sps (q^c-1)\ssp \sum_{s=1}^aq^{s-1} \CF  \bigl(\BA_s\cap \BC\ssp \bigr)\ess .
\cr
}
\eqno 4.14
$$
We have seen in section 3 that the Frobenius characteristic of any direct sum of the modules 
$\BM^{\eee_1\eee_2\cdots \eee_m}$ may be expressed as a polynomial in $\nabla$ applied to $\Phi_\mu$.
It develops that when we express  $L_A$ and $L_B$ in this manneer the terms corresponding to 
$\CF\ssp \bigl(\BA_s\cap \BC\ssp \bigr)$ take a particularly simple form. To this end let $\Pi_A$ and $\Pi_B$
denote the unique polynomials of degree $<m$ such that
$$
 L_A =\Pi_A(\nabla)\ssp \Phi_\mu \ess\ess\ess {\rm and}\ess\ess\ess L_B = \Pi_B(\nabla)\ssp \Phi_\mu\ess\ess\ess
\eqno 4.15
$$
To compute $\Pi_B$ from $\Pi_A$ according to 4.14, we need to make a precise assumption as to 
which corner cell is at the end of row $B$. So let us suppose that it is the $i^{th}$ corner cell 
so that $\BB=\BM_{\aaa^{(i)}}$. In that case, we see that the inductive process, expressed by the
recursion in 4.10, forces all the modules $\BB_s$ to be direct sums of
the submodules
$$
\BM^{\eee_1\eee_2\cdots \eee_{i}}\ses 
\BM_{\aaa^{(1)}}^{\eee_1}\wedge\BM_{\aaa^{(2)}}^{\eee_2}\wedge\cdots \wedge \BM_{\aaa^{(i)}}^{\eee_{i}}
$$
with the exponents $\eee_s$ taking the values $1,0$ according to the conventions we made at the introduction.
This is because we assign $\BM_{\aaa^{(1)}}$ to the top row,
and in constructing the assignments for the row that contains the $i^{th}$ corner we only intersect
with $\BM_{\aaa^{(i)}}$. This given the crucial step in solving the recursion in  4.14 is given by the following
result.

\Prop{4.1}
If $\BM$ is a direct sum of the modules $\BM^{\eee_1\eee_2\cdots \eee_{i-1}}$ then its Frobenius
characteristic is of the form
$$
\CF\ssp \BM\ses \Pi(\nabla )\ssp \Phi_\mu
$$
with $\Pi(x)$ a polynomial divisible by $\prod_{j\geq i}(1-x/T_j)$ and the Frobenius characteristic of $\BM\wedge\BM_{\aaa^{(i)}}$ 
is given by the formula
$$
\CF\ssp\bigl(\ssp \BM\wedge\BM_{\aaa^{(i)}}  \bigr)\ses \Pi'(\nabla )\ssp \Phi_\mu
\eqno 4.16
$$
with
$$
\Pi'(\nabla )\ses {\Pi(\nabla)\over 1-\nabla/T_i}  \ess .
\eqno 4.17
$$

\Proof Clearly it is sufficient to prove the assertion for $\BM=\BM^{\eee_1\eee_2\cdots \eee_{i-1}}$. But in this case we
have
$$
\BM^{\eee_1\eee_2\cdots \eee_{i-1}}\ses 
\bigoplus_{\eta_i=0}^1 \bigoplus_{\eta_{i+1}=0}^1\cdots \bigoplus_{\eta_{i+1}=0}^1
\BM^{\eee_1\eee_2\cdots \eee_{i-1}\eta_{i} \eta_{i+1}\cdots \eta_{m}}\ess .
$$
Thus, from 3.1 and 3.17 we derive that 
$$
\eqalign{
\CF\ssp \BM^{\eee_1\eee_2\cdots \eee_{i-1}}
&\ses 
\sum_{\eta_i=0}^1 \sum_{\eta_{i+1}=0}^1\cdots \sum_{\eta_{i+1}=0}^1
\ssp \prod_{j=1}^{i-1}\Bigl({-\nabla\over T_j}\Bigr)^{1-\eee_j}
\prod_{j=i}^{m}\Bigl({-\nabla\over T_j}\Bigr)^{1-\eta_j}
\ess\ssp \Phi_\mu
\cr
&\ses
  \prod_{j=1}^{i-1}\Bigl({-\nabla\over T_j}\Bigr)^{1-\eee_j}
\prod_{j=i}^{m}\Bigl(1-{\nabla\over T_j}\Bigr) 
\ess\ssp \Phi_\mu
\cr
}
\eqno 4.18
$$
This proves the divisibility by $\prod_{j\geq i}(1-x/T_j)$. Now, since
$$
\BM^{\eee_1\eee_2\cdots \eee_{i-1}}\wedge \BM^{(\aaa_i)}\ses \BM^{\eee_1\eee_2\cdots \eee_{i-1}1}
$$
in the same manner we obtain that
$$
\eqalign{
\CF\ssp \Bigl(\BM^{\eee_1\eee_2\cdots \eee_{i-1}}\wedge \BM^{(\aaa_i)}\Bigr)&\ses
\ssp  \Bigl(\prod_{j=1}^{i-1}\Bigl({-\nabla\over T_j}\Bigr)^{1-\eee_j}\ssp \Bigr)\Bigl({-\nabla\over T_i}\Bigr)^{1-1}
\prod_{j=i+1}^{m}\Bigl(1-{\nabla\over T_j}\Bigr) 
\ess\ssp \Phi_\mu
\cr
&\ses 
\left( \prod_{j=1}^{i-1}\Bigl({-\nabla\over T_j}\Bigr)^{1-\eee_j}
\prod_{j=i}^{m}\Bigl(1-{\nabla\over T_j}\Bigr) \Big /\Bigl(1-{\nabla\over T_i}\Bigr)\right)
\ess  \Phi_\mu\ess .
}
$$
Comparing with 4.18 we see that 4.16 with 4.17 does hold in this case. This completes our proof.
\sas
 
If $\BM$ is a direct sum of the modules $\BM^{\eee_1\eee_2\cdots \eee_m}$, by a slight abuse of
notation, we let $\Pi_{\BM}(z)$ denote the unique polynomial of degree $<m$ such that
$\CF\ssp \BM\ses \Pi_\BM(\nabla)\Phi_\mu$. This given, Proposition 1.1 has the following
beautiful corollary.

\Prop{4.2}
If $B$ is the row that contains the $i^{th}$ corner of $\mu$ then
$$
\Pi_B(z) \ses \Pi_A(z)\ssp {q^c-z/T_i\over 1-z/T_i}\sps
{q^c-1\over q-1} \prod_{j=1,j\neq i}^m\bigl(1-{z\over T_j}\bigr)
\eqno 4.19
$$

\Prop From 4.14 and 4.15 we derive that
$$
\Pi_B(\nabla)\Phi_\mu\ses \Pi_A(\nabla)\Phi_\mu\sps {q^c-1\over q-1}\ssp \Pi_\BC(\nabla)\Phi_\mu
\sps (q^c-1)\ssp \sum_{i=1}^a q^{s-1}  \Pi_{\BA_s\cap \BC}(\nabla)\Phi_\mu\ess .
\eqno 4.20
$$
Since by assumption  $\BC=\BM_{\aaa^{(i)}}$, formula 3.19 a) gives
$$
\Pi_\BC(z)\ses  \prod_{j=1,j\neq i}^m\bigl(1-{z\over T_j}\bigr)\ess .
$$
Recalling that the BH algorithm constructs $A_s$ as a direct sum of the modules $\BM^{\eee_1\eee_2\cdot\eee_{i-1}}$,
we can apply Proposition 4.1 and obtain that
$$
\Pi_{\BA_s\cap \BC}(z)\ses \Pi_{\BA_s}(z)\ess {1\over 1-  z/T_i}\ess . 
$$
Thus 4.20 yields
$$
\Pi_B(z)\ses \Pi_A(z)\sps {q^c-1\over q-1}\prod_{j=1,j\neq i}^m\bigl(1-{z\over T_j}\bigr)
\sps {q^c-1  \over 1-  z/T_i}\ssp \biggl(\sum_{i=1}^a q^{s-1}  \Pi_{\BA_s}(z)\biggr) 
\ess ,
$$
and 4.19 follows because from 4.12 we get that
$$
\sum_{i=1}^a q^{s-1}  \Pi_{\BA_s}(z)\ses \Pi_A(z)\ess .
$$

Note that, since the the BH algorithm keeps the same cell assignments when the row
length doesn't change, the polynomial $\Pi_A(z)$ will be the same as that constructed from
the row that contains the $(i-1)^{st}$ corner cell. Recalling our definitions of corner
weights given in 3.13 and 3.38, we see that the monomial $q^c$ appearing in 4.19 
is none other than the ratio  $x_i/u_{i-1}$. Moreover, since                                                                      
Thus if we change our notation
and denote $\Pi_A$ by $\Pi_{i-1}$ and $\Pi_B$ by $\Pi_i$, the recursion in 4.19 becomes
$$
\Pi_i(z)\ses \Pi_{i-1}(z)\ssp  {x_i/u_{i-1}-z/T_i\over 1- z/T_i}
\sps {x_i/u_{i-1}-1\over q-1}\ssp \prod_{j=1,j\neq i}^m\bigl(1-{z\over T_j}\bigr)\ess .
$$ 
Recalling that we have $T_\mu/T_i=x_i$, we may rewrite this as 
$$
\Pi_i(z)\ses \Pi_{i-1}(z)\ssp {x_i\over u_{i-1}} {1- u_{i-1}\ssp z/T_\mu\over 1- x_i\ssp z/T_\mu}
\sps {1\over u_{i-1}}{x_i-u_{i-1}\over q-1}\ssp \prod_{j=1,j\neq i}^m\bigl(1- x_i\ssp z/ T_\mu\bigr)\ess .
\eqno 4.21
$$
To solve this recursion, it is convenient to work with the expressions 
$$
\GGG_i(z)\ses {\Pi_i\bigl(\ssp z\ssp  T_\mu\bigr)\over \prod_{j=1}^m\bigl(1- x_j\ssp z\bigr)}
\bigsp(\ssp for \ess\ess i=1\ldots m\ess)\ess ,
\eqno 4.22
$$ 
so that 4.21 reduces to
$$
\GGG_i(z)\ses \GGG_{i-1}(z)\ssp {x_i\over u_{i-1}}\ssp {1- u_{i-1}\ssp z \over 1- x_i\ssp z }
\sps {1\over u_{i-1}}{x_i-u_{i-1}\over q-1}\ssp  {1 \over 1- x_i\ssp z }\ess .
\eqno 4.23
$$
Standard methods of the theory of difference equations yield that the solution of this
recurrence may be expressed in the form
$$
\GGG_i(z)\ses a\ssp  \GGG_i^{(1)}(z)\sps \GGG_i^{(2)}(z)\ess .
\eqno 4.24
$$
where $\GGG_i^{(1)}(z)$ and $\GGG_i^{(2)}(z)$  are respectively 
solutions of the homogeneous and non-homogeneous equations, and $a$ is a constant that is the determined by the
initial conditions. Now, it is easily verified that we can take
$$
\GGG_i^{(1)}(z)\ses {x_i\over u_{i-1}}\cdots { x_1\over  u_o}\ess {1-u_{i-1}z\over 1-x_i z}\cdots {1-u_oz \over 1-x_1 z}
\ess\ess\ess
{\rm and}
\ess\ess\ess \GGG_i^{(2)}(z)\ses {1\over 1-q} 
$$ 
Substituting this in 4.24 gives that
$$
\GGG_1(z)\ses a\ssp {x_1\over u_o}\ess {1-u_oz \over 1-x_1 z}\sps {1\over 1-q}\ess .
\eqno 4.25
$$
Now since the BH algorithm assigns $\BM_{\aaa^{(1)}}$ 
to every cell of the top row, and 
$$
\CF\BM_{\aaa^{(1)}}=\prod_{j=2}^m(1-\nabla/T_j)\Phi_\mu\ess , 
$$
we must have 
$$
\Pi_1(z)\ses (1+q+q^2+\cdots +q^{a_1})\ssp \prod_{j=2}^m(1-x_j\ssp z/T_\mu)\ess ,
$$
where $a_1$ is the coarm of the $1^{st}$ corner cell and therefore $a_1+1$ is the 
length of the top row. Thus, since $q^{a_1+1}=x_1/u_o\ssp ,$ 4.22 gives
$$
\GGG_1(z)\ses {1\over u_o}\ssp {u_o-x_1\over 1-q }{1 \over 1-x_1 z}\ess .
\eqno 4.26
$$
Equating the righthand sides of 4.25 and 4.26 and solving for $a$ gives 
$$
a\ses -{1\over 1-q}
$$
Using this in 4.24 we finally obtain that 
$$
\GGG_i(z)\ses -{1\over 1-q}\ess 
{x_i\over u_{i-1}}\cdots { x_1\over  u_o}\ess {1-u_{i-1}z\over 1-x_i z}\cdots {1-u_oz \over 1-x_1 z}
\sps {1\over 1-q}\ess .
\eqno 4.27
$$

Note that from our notation 
\footnote{$\dag$}
{Recall that the exponent $l_i$ in $x_i=t^{l_i}q^{a_i}$
gives the coleg of the $i^{th}$ corner cell
}
it follows that there are exactly $l_i-l_{i+1}$ rows that have the same length as that
which contains the $i^{th}$ corner cell.  Thus, under the BH algorithm,
each of these rows contributes to the right hand side of 4.4 a term of the form $t^s \Pi_i(\nabla)\Phi_\mu$
with  $\ssp l_{i+1}< s\leq l_i$. 

In summary the BH algorithm yields that
$$
\eqalign{
\delta_{p_1}\ssp \TH_\mu(x;q,t)&\ses 
\sum_{i=1}^m\ssp t^{1+l_{i+1}}\Bigl(1+t+\cdots +t^{l_i-l_{i+1}-1}\Bigr)\ess  \Pi_i(\nabla)\Phi_\mu
\cr
&\ses
{1\over 1-t}\ssp
\sum_{i=1}^m\ssp t^{1+l_{i+1}}\bigl(\ssp 1-t^{l_i-l_{i+1}}\bigr)\ess  \Pi_i(\nabla)\Phi_\mu
\cr
}
\eqno 4.28
$$ 
Note that since $t^{1+l_m}=x_m/u_m$ and $t^{l_j-l_{i+j}}=x_j/u_j$ we may rewrite 4.28 as
$$
\del_{p_1}\TH_\mu(x;q,t)\ses 
{1\over 1-t}\ssp
\sum_{i=1}^m\ssp {x_m\over u_m}\cdots {x_{i+1}\over u_{i+1}}\ssp 
\Bigl(  1-{x_i\over u_i} \ssp  \Bigr)\ssp  \Pi_i(\nabla)\Phi_\mu\ess .
\eqno 4.29
$$
Now, 4.27 gives that
$$
\eqalign{
(1-q)\ssp \sum_{i=1}^m\ssp {x_m\over u_m}\cdots &{x_{i+1}\over u_{i+1}}\ssp 
\Bigl(  1-{x_i\over u_i} \ssp  \Bigr)\ssp  \GGG_i(z) 
\ses \cr
&\ses \sms 
\sum_{i=1}^m\ssp {x_m\over u_m}\cdots {x_{i+1}\over u_{i+1}}\ssp \Bigl(  1-{x_i\over u_i} \Bigr)
{x_i\over u_{i-1}}\cdots { x_1\over  u_o}\ess {1-u_{i-1}z\over 1-x_i z}\cdots {1-u_oz \over 1-x_1 z} 
\cr
&\ess\ess\ess\ess\ess\ess
\bigsp\bigsp\bigsp\bigsp
\sps  \sum_{i=1}^m\ssp {x_m\over u_m}\cdots {x_{i+1}\over u_{i+1}}\ssp
\Bigl(  1-{x_i\over u_i} \ssp  \Bigr)\ess .
\cr
}\eqno{4.30}
$$ 
Calling, for a moment, $\bf S_1$ and $\bf S_2$ the first and second sums on the righthand side,
and recalling from 3.45 that $x_ox_1\cdots x_m=u_ou_1\cdots u_m$, we may write
$$
\eqalign{
{\bf S_1}&\ses \sms 
\sum_{i=1}^m\ssp {u_i\over x_o}\ssp  \Bigl(  1-{x_i\over u_i} \Bigr)
\ssp  {1-u_{i-1}z\over 1-x_i z}\cdots {1-u_oz \over 1-x_1 z}
\cr
&\ses 
\sms \sum_{i=1}^m\ssp {u_i\over x_o}\ssp 
\ssp  {1-u_{i-1}z\over 1-x_i z}\cdots {1-u_oz \over 1-x_1 z}
\sps 
\sum_{i=1}^m\ssp {x_i\over x_o}  
\ssp  {1-u_{i-1}z\over 1-x_i z}\cdots {1-u_oz \over 1-x_1 z}
\cr
&\ses
{1\over x_o z}
\biggl\{\ssp 
 \sum_{i=1}^m\ssp 
\ssp  { (1-u_{i}z) \cdots (1-u_oz )   \over (1-x_i z)\cdots (1-x_1 z)}
\sms 
\sum_{i=1}^m 
\ssp  { (1-u_{i-1}z) \cdots (1-u_oz )   \over (1-x_{i-1} z)\cdots (1-x_1 z)}
\ssp \biggr\}
\cr
&\ses
{1\over x_o z}
\biggl\{\ssp 
{ (1-u_{m}z) \cdots (1-u_oz )   \over (1-x_m z)\cdots (1-x_1 z)}\sms (1-u_o z)
\ssp \biggr\}\ess .
\cr
}
$$
At the same time we have (again because of 3.45)
$$
\eqalign{
{\bf S_2}&\ses  \sum_{i=1}^m\ssp {x_m\over u_m}\cdots {x_{i+1}\over u_{i+1}}\ssp
\sms \sum_{i=1}^m\ssp {x_m\over u_m}\cdots {x_{i }\over u_{i }}\ssp
\cr
&\ses
1 \sms {x_m\over u_m}\cdots {x_{ 1}\over u_{ 1}}\ses 1- {u_o\over x_o}\ses 
{(1-u_oz) \sms (1-x_o z)\over x_o z }
\cr 
}
$$ 
Using these expressions for $\bf S_1$ and $\bf S_2$ in 4.30 gives
$$
\sum_{i=1}^m\ssp {x_m\over u_m}\cdots {x_{i+1}\over u_{i+1}}\ssp 
\Bigl(  1-{x_i\over u_i} \ssp  \Bigr)\ssp  \Pi_i(z)\ses 
{1\over 1-q}\ssp {1\over x_o z}
\biggl\{\ssp 
{ (1-u_{m}z) \cdots (1-u_oz )   \over (1-x_m z)\cdots (1-x_1 z)}\sms (1-x_o z)
\ssp \biggr\}\ess .
$$
Thus from 4.22 we get (recalling that $x_o=1/qt$)
$$
\sum_{i=1}^m\ssp {x_m\over u_m}\cdots {x_{i+1}\over u_{i+1}}\ssp 
\Bigl(  1-{x_i\over u_i} \ssp  \Bigr)\ssp  \Pi_i(z)\ses
{1\over (1-1/t)(1-1/q)}\ssp {T_\mu\over z}\ssp
\biggl\{\ssp 
\prod_{i=0}^m(1-u_i\ssp  {z\over T_\mu})\sms \prod_{i=0}^m(1-x_i\ssp  {z\over T_\mu}) 
\biggr\} 
$$
So that 4.29 becomes
$$
\del_{p_1} \TH_\mu(x;q,t)\ses  
{1\over (1-1/t)(1-1/q)}\ssp {T_\mu\over \nabla}\ssp
\biggl\{\ssp 
\prod_{i=0}^m \Bigl(1-u_i\ssp  {\nabla\over T_\mu}\Bigr)\sms \prod_{i=0}^m\Bigl(1-x_i\ssp  {\nabla\over T_\mu}\Bigr) 
\biggr\}\ssp \Phi_\mu\ess , 
$$
which is in complete agreement with 3.43. 
It is indeed remarkable that these two quite
distinct paths lead exactly to  the same formula. In fact, in this section we have 
basically deduced it  only from the combinatorics of the BH algorithm.
On the other hand, in section 3. we made heavy use of Science Fiction identities
as well as various identities of the theory of Macdonald polynomials.
This yields simultaneous  support to the BH algorithm,  the SF heuristic and the $C=\TH$ conjecture. 
Since it all fits together so well, it is difficult to believe that it results from  purely accidental
circumsances. 
\sa

Needeless to say, the BH algorithm applied to any of the special cases we have studied in this
paper does produce a basis for the corresponding module $\BM_\mu$. It may be also successfully applied
to some infinite family of partitions. As an instance in point we shall treat the case
of $2$-corner partitions. So let $\mu$ have $l_a$ rows of length
$a$ and  $l_b$ rows of length $b$ with $c=b-a>0$. Let $\mu$ be a partition of $n+1$ and $\aaa$ and 
$\bbb$ denote the two partitions of $n$ obtained by removing one of the two corners of $\mu$.
More precisely, we set $\aaa=a^{l_a-1}b^{l_b}$ and $\bbb=a^{l_a}b^{l_b-1}$. Moreover, 
let us decompose the diagram of $\mu$ into $4$ parts as follows:
\sas

\itemitem{\bf For $\bf c>a$ } 
$$
\cases{
UPR\ses\Bigl\{ (r,s)\ssp :\ssp \matrix {l_b\leq r<l_b+l_a\ssp\cr \ssp 0\leq s<a\ssp\cr} \Bigr\}\cr\cr
INS\ses \Bigl\{ (r,s)\ssp :\ssp \matrix {0\leq r<l_b\ssp\cr \ssp 0\leq s<a\ssp\cr} \Bigr\}\cr\cr
MID\ses\Bigl\{ (r,s)\ssp :\ssp \matrix {0\leq r<l_b\ssp\cr \ssp a\leq s<c\ssp\cr} \Bigr\}\cr\cr
OUT\ses\Bigl\{ (r,s)\ssp :\ssp \matrix {0\leq r<l_b\ssp\cr \ssp c\leq s<b\ssp\cr} \Bigr\}\cr
}
$$
and
\itemitem{\bf For $\bf c\leq a$ } 
$$
\cases{
UPR\ses\Bigl\{ (r,s)\ssp :\ssp \matrix {l_b\leq r<l_b+l_a\ssp\cr \ssp 0\leq s<a\ssp\cr} \Bigr\}\cr\cr
INS\ses \Bigl\{ (r,s)\ssp :\ssp \matrix {0\leq r<l_b\ssp\cr \ssp 0\leq s<c\ssp\cr} \Bigr\}\cr\cr
MID\ses\Bigl\{ (r,s)\ssp :\ssp \matrix {0\leq r<l_b\ssp\cr \ssp c\leq s<a\ssp\cr} \Bigr\}\cr\cr
OUT\ses\Bigl\{ (r,s)\ssp :\ssp \matrix {0\leq r<l_b\ssp\cr \ssp a\leq s<b\ssp\cr} \Bigr\}\cr
}
$$

\noindent
Now let $\CB_{A\cup B}$ be a basis for $\BM_\aaa\vee \BM_\bbb$, let  $\CB_{A\cap B}$
be a basis for $\BM_\aaa\wedge \BM_\bbb$ and let $\CB_{A^\perp\cap B}^*$ denote a collection
whose flip with respect to $\Delta_\bbb$ forms a basis for the submodule $\BM_\aaa^\perp\wedge\BM_\bbb$.
This given, we have the following result.

\Thm{4.1}
If the $n!$ conjecture holds true for $\aaa$ and $\bbb$ and the intersection of the two modules 
$\BM_\aaa$ and $\BM_\bbb$ has dimension $n!/2$ then $\dim \BM_\mu$ has $(n+1)!$ dimensions and 
a basis for $\BM_\mu$ is given by the collection
$$
\CB_\mu\ses \sum_{(r,s)\in \mu}\ssp \Bigl\{ b(\del )\del_{x_{n+1}}^r\del_{y_{n+1}}^s  \Delta_\mu\ssp \Bigr\}_{b\in\CB_{r,s}} 
\eqno 4.31
$$
with the following choices:
$$
\CB_{r,s}=\CB_{A}\ess  for  \ess (r,s)\in UPR\ess ,
$$ 
$$
\CB_{r,s}=\CB_{A\cup B}\ess  for  \ess (r,s)\in INS
\ess\ess\ess and \ess\ess\ess \CB_{r,s}=\CB_{A^\perp\cap B}^*\ess  for \ess  (r,s)\in OUT 
$$
while for $(r,s)\in MID$
$$ 
\CB_{r,s}\ses 
\cases
{
\CB_B &    when $\ssp c>a$ 
\cr
\cr
 \CB_A &    when  $\ssp c\leq a$ 
\cr
}
$$

\Proof It is not difficult to see that the BH algorithm yields precisely the choices indicated above 
except that for $(r,s)\in OUT$ it gives $\CB_{r,s}=\CB_{A\cap B}$. Of course if we assume 
the validity of the SF heuristic then a choice is guaranteed to exist for a $\CB_{A\cap B}$ 
whose flip with respect to $\Delta_\bbb$ is in fact a basis for $\BM_\aaa^\perp\cap\BM_\bbb$. 
However, in the two corner case, we need not assume anything more than $\dim \BM_\aaa\cap\BM_\bbb=n!/2$.
Since a simple count shows that, with these choices, $\CB_\mu$ has $(n+1)!$ elements, we need only
show that for any choices of $b_{r,s}\in \CL[\CB_{r,s}]$ we cannot have
$$
\sum_{(r,s)\in \mu}\ssp b_{r,s}(\del )\del_{x_{n+1}}^r\del_{y_{n+1}}^s  \Delta_\mu\ses 0
\eqno 4.32
$$ 
without $b_{r,s}=0$ identically for all $(r,s)\in \mu$. Proceeding as we did in
the several examples we have previously considered, we can derive $n+1$ equations from 4.32
by equating to zero the coefficients of $x_{n+1}^{r_o}y_{n+1}^{s_o}$ for each $(r_o,s_o)\in \mu$.
These equations are all of the form
$$
\sum_{
\multi{(r',s')\in \mu\cr (r'+r_o,s'+s_o)\in \mu\cr }} \hskip -.2truein b_{(r',s')}(\del)D_{r'+r_o,s'+s_o}\ses 0
\eqno 4.33
$$     
where $D_{r,s}$ denotes a suitable multiple of the cofactor of the monomial $x_{n+1}^ry_{n+1}^s$ in the matrix
whose determinant gives $\Delta_\mu$. For convenience let us write $(r' ,s' )\prec (r,s)$ if and only if 
$(r',s')\neq (r,s)$ and 
$r'\leq r$ and $s'\leq s$.

Now, a close examination of few special cases reveals that the equations
in 4.33  may be successively solved to yield one of the following four types of conditions
\sas

\item {(1)} For a $b_{r,s}$ with $(r,s)\in INS$: 
\hfill(after showing that $b_{r',s'}=0$ for all $(r' ,s' )\prec (r,s)$) 
$$
b_{r,s}(\del)\Delta_\aaa= 0
\ess\ess\ess 
{\rm and}
\ess\ess\ess 
b_{r,s}(\del)\Delta_\bbb= 0
\eqno 4.34
$$

\item {(2)} For a $b_{r,s}$ with $(r,s)\in MID$: \hfill (after showing that $b_{r',s'}=0$ for all $(r' ,s' )\prec (r,s)$) 
$$
\cases{
b_{r,s}(\del)\Delta_\bbb= 0 & if $c>a$\cr\cr
b_{r,s}(\del)\Delta_\aaa= 0 & if $c\leq a$\cr}
\eqno 4.35
$$

\item {(3)} For some $b_{r,s}$ with $(r,s)\in UPR$:\hfill (after showing that $b_{r',s'}=0$ for all $(r' ,s' )\prec (r,s)$) 
$$
b_{r,s}(\del)\Delta_\aaa= 0
\eqno 4.36
$$

\item {(4)} Alternatively, for a pair $b_{r_1,s_1}$, $b_{r_2,s_2}$  with $(r_1,s_1)\in UPR$ and 
$(r_2,s_2)\in OUT$

\hfill (after showing that $b_{r',s'}=0$ for all $(r' ,s' )\prec (r_1,s_1)$ and for all $(r' ,s' )\prec (r_2,s_2)$)
$$
b_{r_1,s_1}(\del)\Delta_\aaa\sps b_{r_2,s_2}(\del)\Delta_\bbb = 0\ess .
\eqno 4.37
$$

Now in any of the cases (1), (2), and (3) the condition forces $b_{r,s}=0$. For instance in case (1), the choice 
$\CB_{r,s}=\CB_{A\cup B}$ yields that $b_{r,s}\in \BM_\aaa\vee\BM_\bbb$ and at the same time the condition in 
4.34 yields that $b_{r,s}\in \big(\BM_\aaa\vee\BM_\bbb\bigr)^\perp$. Thus $b_{r,s}$ must be orthogonal to itself
and therefore it must identically vanish. The conclusion in the other two cases are similarly obtained
since the corresponding choice of $\CB_{r,s}$ combined with 4.35 or 4.36, as the case may be, forces $b_{r,s}$
to be orthogonal to itself and therefore equal to zero. In case (4), since $\CB_{r_1,s_1}=\CB_A$,
the polynomial $b_{r_1,s_1}(\del)\Delta_\aaa$ may turn out be an arbitrary element of $\BM_\aaa$. So
the only way to force  $b_{r_1,s_1}$  to vanish is to assure that 4.37 forces $b_{r_1,s_1}(\del)\Delta_\aaa=0$,
and that will certainly be forced by 4.37 if  $b_{r_2,s_2}(\del)\Delta_\bbb\in \BM_\aaa^\perp$. 
Now the latter is guaranteed by our choice of $\CB_{r_2,s_2}=\CB_{A^\perp\cap B}^*$ for all $(r_2,s_2)\in OUT$.
This completes our proof.
\sas

We should note that the role played  by the condition $\dim(\BM_\aaa\wedge\BM_\bbb)=n!/2$ in the above proof
is simply to assure that the orthogonal complement of $\BM_\aaa\wedge\BM_\bbb$ in $\BM_\bbb$ has
also dimension $n!/2$ and this is the additional condition that needs to be satisfied so that in the final
count $\CB_\mu$ has exactly $(n+1)!$ elements.
\sas

There are a number of other identities and conjectures that may be derived from the SF heuristics,
these will be included in a forthcoming joint work with N. Bergeron,  M. Haiman. and G. Tessler [\BBGHM]
where many of the identities derived here are shown to have a remarkably suggestive combinatorial interpretation.
We should also refer the reader to the work of C. Chang [\chang] where Theorem 4.1 is used to give a unified
proof of the $n!$ conjecture for hook, two-row and two-column shapes.

\section{References}
\parindent=.25truein

\item{[\butler]}L. Buttler, 
(Personal communication).
\vskip .1truein

\item{[\chang]}  C. Chang, {\sl Doctoral Dissertation} UCSD (1997)
\vskip .1truein

\item {[\BBGHM]} F. Bergeron, N. Bergeron, A. M. Garsia,  M. Haiman and G. Tessler,
{\sl Ferrers diagrams with a missing cell and the $n!$ conjecture},
(in preparation)
\vskip .1truein

\item{[\bergeronhamel]}
F. Bergeron and S. Hamel, 
{\sl Intersection of Modules related to Macdonald's Polynomials\/},
(in preparation)
\vskip .1truein

\item {[\bergerongarsia]}   N. Bergeron and A. M. Garsia,
{\sl On Certain Spaces
of Harmonic Polynomials, Hypergeometric Functions on domains of
Positivity, Jack polynomials and Applicatioions\/}, Contemporary Mathematics, {\bf
138} (1992) 51-86..
\vskip .1truein

\item{[\garsiahaiman]}     
A. M. Garsia and M.  Haiman,
{\sl A graded representation module for Macdonald's polynomials},
Proc. Natl. Acad. Sci. USA V 90 (1993) 3607-3610. 
\vskip .1truein

\item {[\garsiahaimanorbit]}
A. M. Garsia and M.  Haiman, {\sl Orbit Harmonics and Graded Representations} 
(Research Monograph to appear as part of the Collection Published by the
Lab.~de Comb.~et Informatique Math\'ematique, edited by S. Brlek, 
U. du Qu\'ebec \`a Montr\'eal).
\vskip .1truein

\item{[\garsiahaimanpieri]}
A. M. Garsia and M. Haiman, {\sl Factorizations of Pieri rules for Macdonald polynomials},
Discrete Mathematics 139 (1995) 219-256.
\vskip .1truein

\item{[\garsiahaimancatalan]}
A. Garsia and M. Haiman, {\sl A Remarkable q,t-Catalan Sequence and q-Lagrange inversion},
J. of Alg. Comb. V. 5 (1996) pp. 191-244. 
\vskip .1truein

\item{[\garsiahaimankostka]}
A. Garsia and M. Haiman, {\sl Some bigraded $S_n$-modules and the Macdonald q,t-Kostka 
coefficients}, Electronic Journal of Alg. Comb. V. 3 \#2 (1996) pp. 561-620.  

(web site http://ejc.math.gatech.edu:8080/Journal/journalhome.html).
\vskip .1truein

\item {[\garsiaprocesi]}
A. M. Garsia and C. Procesi,
{\sl On certain graded $S_n$-modules and the q-Kostka polynomials}, 
Advances in Mathematics {\bf 94} (1992) 82-138.
\vskip .1truein

\item {[\grasiaatesler]}
A. M. Garsia and G. Tessler,
{\sl Plethystic Formulas for the Macdonald $q,t$-Kostka coefficients},
Advances in Mathematics V. 123 \#2 Nov. 1996 pp.  144-222.
\vskip .1truein

\item {[\garsiaremmel]}
A. M. Garsia and J. Remmel,
{\sl Plethystic Formulas and positivity for $q,t$-Kostka  Coefficients,
to appear in the Proceedings of Rotafest}  ( a referred Volume
in Honor of G. C. Rota)
\vskip .1truein

\item {[\garsiateslerhaiman]}
A. M. Garsia, G. Tessler and  M. Haiman,
{\sl Explicit Plethystic Formulas for the Macdonald  $q,t$-Kostka  Polynomials}
(in Preparation) 
\vskip .1truein

\item{[\haiman]}
M. Haiman, {\sl Macdonald Polynomials and Hilbert Schemes}, (preprint). 
\vskip .1truein

\item {[\kirillovnoumi]}
A. Kirillov and M. Noumi,
{\sl Raising operators for Macdonald Polynomials},
(preprint).
\vskip .1truein

\item {[\knop]}
F. Knop, {\sl Integrality of Two Variable Kostka Functions},
(preprint).
\vskip .1truein

\item {[\knopquantum]}
F. Knop, {\sl Symmetric and non-symmetric Quantum Capelli Polynomials},
(preprint).
\vskip .1truein

\item {[\lapointevinet]}
L. Lapointe and L. Vinet,
{\sl Creation Operators for Macdonald and Jack Polynomials},
(preprint)
\vskip .1truein

\item {[\macdonald]}
I. G. Macdonald, {\sl A new class of symmetric functions}, 
Actes du $20^e$ S\'eminaire Lotharingien, 
Publ. I.R.M.A. Strasbourg, (1988)
131-171.
\vskip .1truein

\item {[\macdonaldbook]}
I. G. Macdonald, {\sl Symmetric functions and Hall polynomials},
Second Edition, Clarendon Press, Oxford (1995).
\vskip .1truein

\item {[\reiner]} E. Reiner, {\sl A Proof of the $n!$ Conjecture for Generalized Hooks},
to appear in the Journal of Combinatorial Theory, Series A.  
\vskip .1truein

\item {[\sahi]}
S. Sahi, {\sl Interpolation and integrality for Macdonald's Polynomials},
(preprint)
\vskip .1truein

\item {[\youngsubst]}
A. Young, {\sl On quantitative substitutional analysis} (sixth paper),
The collected papers of A. Young, University of Toronto Press (1977)
pp. 434-435.  

\end

However, it can also be
defined in terms of the difference operator $\Delta$ which acts on a 
symmetric polynomial $P$ according to the formula
$$
\Delta\ssp P[X]\ses P[X]\sms P\bigl[X+{\textstyle {(1-t)(1-q)\over z}}\bigr]\Omega[-zX]\ssp
|_{z^0}\ess .
\eqno 3.12
$$
where for any expression $E$ we set
%\footnote{($\dag$)}
%{Here ``$|_{z^0} $'' extracts the constant term} 
$$
\Omega[E]\ses exp\biggl[\ssp \sum_{k\geq 1}\ssp {1\over k}\ssp p_k[E]\ssp \biggr]\ess .
$$
Using Macdonald theory it is shown in [\garsiahaimancatalan] that we have
$$
\Delta\ssp \TH_\mu\ses (1-t)(1-q)\ssp B_\mu(q,t)\ssp \TH_\mu\ess ,
\eqno 3.13
$$
where
$$
B_\mu(q,t)\ses \sum_{s\in \mu}\ssp t^{l'_\mu(s)}\ssp q^{a'_\mu(s)}
\eqno 3.14
$$
Since 3.13 shows that $\Delta$ and $\nabla$ have the same set of eigenfunctions
and the eigenvalues of $\Delta$ are distinct we see that $\nabla$ can be expressed
as a polynomial in $\Delta$.
\sas